\documentclass{siamltex}

\usepackage[
  hyperindex  = {true},
  colorlinks  = {true},
  linkcolor   = {blue},
  citecolor   = {blue}
]{hyperref}
\usepackage{multirow}
\usepackage[margin=0pt,labelsep=space,labelfont=bf,textfont=it,up,belowskip=-8pt,aboveskip=5pt]{caption}
\usepackage{latexsym,graphicx, amssymb, amsfonts, subfigure, setspace, geometry, amsmath}
\usepackage{pgfplots}
\pgfplotsset{compat=1.8}
\usepgfplotslibrary{statistics}
\usepackage{grffile}
\usepackage{pgfplotstable}
\usepackage{slashbox}
\usepackage{mathtools}
\usepackage[utf8]{inputenc}
\usepackage[english]{babel}
\usepackage{float}
\usepackage{esvect}
\usepackage{algorithm,algorithmic}
\usepackage{enumitem}

\pgfplotsset{compat=newest}

\makeatletter
\newcommand{\doublewidetilde}[1]{{%
  \mathpalette\double@widetilde{#1}%
}}
\newcommand{\double@widetilde}[2]{%
  \sbox\z@{$\m@th#1\widetilde{#2}$}%
  \ht\z@=.9\ht\z@
  \widetilde{\box\z@}%
}
\makeatother

\newtheorem{prob}{Problem}

\newtheorem{remark}{Remark}[section]

\newcommand{\xb}{\mathbf{x}}

\newcommand{\db}{\mathbf{d}}
\newcommand{\Fb}{\mathbf{F}}
\newcommand{\Jb}{\mathbf{J}}

\newcommand{\Tb}{\mathbf{T}}

\newcommand{\parB}{{\partial \mathcal{B}}}

\newcommand{\ipoint}[1]{\textit{\textbf{#1}}}

\DeclareMathOperator*{\argmin}{arg\,min}

\pgfmathdeclarefunction{gauss}{3}{%
  \pgfmathparse{1/(#3*sqrt(2*pi))*exp(-((#1-#2)^2)/(2*#3^2))}%
}

\title{Reconstruction of a compactly supported sound profile in the presence of a random background medium}

\author{Carlos Borges\thanks{Institute for Computation Engineering and Sciences, University of Texas, Austin, TX}
\and
George Biros\thanks{Department of Mechanical Engineering and Institute for Computation Engineering and Sciences, University of Texas, Austin, TX}
}

\begin{document}

\maketitle

\begin{abstract}
In this paper, we present algorithms for reconstructing an unknown compact scatterer embedded in a random noisy background medium, given measurements of the scattered field and information about the background medium and the sound profile. We present six different methods for the solution of this inverse problem using different amounts of scattered data and prior information about the random background medium and the scatterer. The different inversion algorithms are defined by a combination of stochastic programming methods and Bayesian formulation. Our basic results show that if we have data for just one instance of the random background medium the best strategy is to invert for both random medium and unknown scatterer with appropriate regularization. However, if we have data for multiple instances of the medium it may be worth solving a coupled set of multiple inverse problems. We present several numerical results for inverting for various scatterer geometries under different inversion scenarios. The main take-away of our study is that one should invert for both unknown scatterer and random medium, with appropriate, prior-information based regularization. Furthermore, if data from multiple realizations of the background medium is available, then it may be beneficial to combine results from multiple inversions. 
\end{abstract}


\section{Introduction}\label{s:intro}
Inverse scattering finds applications in  medical imaging  \cite{Hoskins,kuchment2014radon,
2006-MMRAMI-WM,nashed2002inverse,scherzer2010handbook,simonetti2008inverse,Simonetti2013,Tsui201683, Blackledge1985,pastorino2008medical}, remote sensing \cite{Ustinov2014,Wang01022012,jin2006theory}, ocean acoustics  \cite{chavent2012inverse,0266-5611-10-5-003}, nondestructive testing \cite{collins1995nondestructive,engl2012inverse,marklein2006inverse,marklein2005numerical,langenberg1993imaging,Bleistein1983},  geophysics \cite{aster2013parameter,tarantola2013inverse,zhdanov2002geophysical,snieder1999inverse,Eldad_04}, and sonar and radar \cite{cheney2009fundamentals,Colton,Borges2017}. In this work, we recover a compactly supported \ipoint{unknown scatterer}, denoted $q(\xb)$, in the presence of a compactly supported \ipoint{random background medium} denoted by $\eta(\xb)$, from far field acoustic scattering measurements.   We are not interested (and don't have enough data) to exactly recover $\eta(\xb)$ but we have some prior information on its statistics. What we're interested in is recovering $q(\xb)$ as accurately as possible, given the prior information on $\eta$.

\begin{figure}
\center
\begin{tikzpicture}[scale=0.60]

\draw[ultra thick,->] (-3-1,2.5) -- (-1-1,2.5);
\draw (-2.5-1,1.5) -- (-2.5-1,3.5);
\draw (-2-1,1.5) -- (-2-1,3.5);
\draw (-1.5-1,1.5) -- (-1.5-1,3.5);

\node at (0,-.1) { };
\node at (0,5.1) { };

\draw[dashed] (0,0) rectangle (5,5);

\draw[dashed] (2.5,2.5) circle [radius=4];

\draw [fill=blue!40!white] (1, 1)
  .. controls ++(165:-1) and ++(240: 1) .. ( 4, 1)
  .. controls ++(240:-1) and ++(165:-1) .. ( 4, 4)
  .. controls ++(165: 1) and ++(175:-2) .. (2, 4.5)
  .. controls ++(175: 2) and ++(165: 1) .. ( 1, 1);

\draw [fill=black!40!white] (1.5, 1.5)
  .. controls ++(165:-1) and ++(240: 1) .. (4 , 1.5)
  .. controls ++(240:-1) and ++(165:-1) .. ( 3.5, 3)
  .. controls ++(165: 1) and ++(175:-2) .. (1.8, 3)
  .. controls ++(175: 2) and ++(165: 1) .. ( 1.5,1.5 );

\node at (0.5,0.5) {$\Omega$};
\node at (-2-1,1) {$u^{\emph{inc}}$};
\node at (8,1) {$u^{\emph{scat}}$};
\node at (2.5,2.2) {$\text{supp}(q)$};
\node at (2.5,3.5) {$\text{supp}(\eta)$};
\node at (6,-0.5) {$\parB$};

\draw[ultra thick,->] (6+1,2) -- (7.75+1,2.75);
\draw (6.5+1,1.25) -- (6+1,3);
\draw (7+1,1.5) -- (6.5+1,3.25);
\draw (7.5+1,1.75) -- (7+1,3.5);


\end{tikzpicture}

\caption{Scattering from a compactly supported inhomogeneity with scatterer
 $q(\xb)$ in the presence of a background medium represented by a random variable with realization $\eta(\xb)$.
In the {\em forward scattering problem}, $q(\xb)$ and $\eta(\xb)$ are known and one seeks to compute the scattered field,
either within $\Omega$ or on the boundary $\partial\mathcal{B}$ of a disk.
In the {\em inverse scattering problem}, $q(\xb)$ is unknown and we seek
to determine it from measurements of the scattered field on $\partial\mathcal{B}$ given prior knowledge of general characteristics  $q(\xb)$ and $\eta(\xb)$.}\label{fig1}
\end{figure}
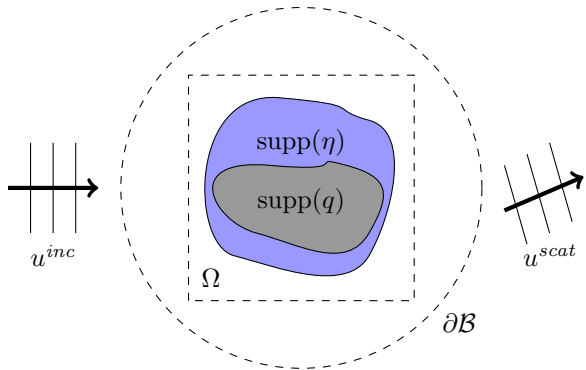

{\bf Problem Statement:} We consider the scattering problem in two dimensions. We assume that both $q(\xb)$ and  $\eta(\xb)$ are in $C^2_0(\Omega)$, where $\Omega\in\mathbb{R}^2$. We define the forward scattering operator $\mathcal{F}:C^2_{0}(\Omega)\rightarrow L^2(\partial \mathcal{B})$ by
\begin{equation}\label{e:forward}
 u_\eta^{\emph{scat}} = \mathcal{F} (q+\eta; u^{\emph{inc}}),
\end{equation}
where $u^{\emph{inc}}$ is the incident field, $\partial\mathcal{B}$ is the boundary of a disk centered at the origin and ${\rm supp}(q)$ is in the interior of $\mathcal{B}$.  The operator $\mathcal{F}$ is well defined since the forward scattering problem is well posed \cite{Colton}. To obtain the value of the scattered data $u_\eta^{\emph{scat}}$ at a point $\xb$, we must solve the variable coefficient Helmholtz equation or its integral equation counterpart, the Lippmann-Schwinger equation \cite{Colton,Nedelec}. (In the rest of the paper, unless it is otherwise stated, we suppress the notational dependence of $\mathcal{F}$ on $u^{\emph{inc}}$.)

Here our main interest is the inverse problem of recovering an approximation $\tilde{q}$ of the scatterer $q$ given \ipoint{measurements} $\db_\eta$ of the acoustic field  scattered by $q$ and $\eta$ by solving a nonlinear least squares problem given by
\begin{equation}
\tilde{q}=\argmin_q \frac{1}{2}\|\db_\eta-\Fb(q+\eta)\|^2,
\label{eq:invprob}
\end{equation}
where the $i^\text{th}$ component of $\Fb$ is the evaluation of $\mathcal{F}$ at the $i^\text{th}$ location at which field measurements are taken. In general, this problem is nonlinear and ill posed, as is the deterministic inverse scattering problem. In the most part of the paper, we will assume that we have high-quality, noise-free data from multiple frequencies, and multiple illuminations. In the absence of $\eta$,  after appropriate discretization of $q$, \eqref{eq:invprob} becomes well posed~\cite{Borges2017}. The main complication will be due to the presence of $\eta$, which will be the main source of the difficulty of reconstructing $q$, since scattering from $q$ and $\eta$ is mixed by the forward operator.

{\bf Notation:} We present the most common symbols used in this paper in Table \ref{table:symbol}. 

\begin{table}
\caption{List of main symbols used in this article.}\label{table:symbol}
{\small
\begin{center}
\begin{tabular}{ll}
\hline
Symbol & Description \\
\hline\hline
$k$                           & Wavenumber (or frequency) \\
$q$                           & Scatterer of the domain \\
$\eta$                       & Noisy background medium \\
$\mathcal{F}$           & Analytical forward scattering operator \\
$\Fb$                        & Operator $\mathcal{F}$ calculated at points where measurements are taken\\
$\mathcal{J}$           & Fr\'{e}chet derivative of $\mathcal{F}$ \\
$\Jb$                        & Operator $\mathcal{J}$ calculated at points where measurements are taken\\
$\mathbb{E}(\eta)$   & Expected value of $\eta$ \\
$\mathcal{T}_\eta$   & Isotropy tensor for $\eta$ \\
$\mathcal{T}_q$       & Isotropy tensor for $q$ \\
$p(\eta)$                   & Probability distribution of $\eta$ \\
$p(q)$                       & Probability distribution of $q$ \\
$\alpha$                    & Regularization parameter for $\eta$ \\
$\beta$                     & Regularization parameter for $q$ \\
$\theta$                     & Vector with incidence direction of plane wave \\
$\parB$                     & Circle around the domain where measurements are obtained \\
$u^{\emph{inc}}$          & Incident plane wave \\
$u^{\emph{scat}}$         & Scattered field off of $q$ \\
$u_\eta^{\emph{scat}}$ & Scattered field off of $q+\eta$ \\
$u_s^{\emph{scat}}$     & Scattered field off of $q+\eta_s$, $\eta_s$ is the $s^\text{th}$ realization of $\eta$ \\
$\db$                         & Measurements of $u^{\emph{scat}}$ on $\partial\mathcal{B}$ \\
$\db_\eta$                 & Measurements of $u_\eta^{\emph{scat}}$ on $\partial\mathcal{B}$ \\
$\db_s$                     & Measurements of $u_s^{\emph{scat}}$ on $\partial\mathcal{B}$ \\
$\xb$                         & Position in $\mathbb{R}^2$ space \\
$C_I$                        & Cost of the solution of one instance of the inverse problem \\
$N_s$                       & Number of samples (realizations) of $\eta$ \\
$N$                           & Number of points in the discretization of the domain \\
\hline 
\end{tabular}
\end{center}
}
\end{table}

{\bf Methodology:} We present six methods for the solution of \eqref{eq:invprob} that require different amounts of scattered data and different information about the prior knowledge of the probability distribution of the object and the background medium.

a) {\bf Single Inversion Single Data No Prior (SISDNP): }
In this scenario, we assume that the data have been generated by a single realization of $\eta$ but we only have access to the expected value of $\eta$. We  use the measurements $\db_\eta$ and the expected value $\mathbb{E}(\eta)$ to obtain the approximation $q_\mathrm{SISDNP}$ of the scatterer by solving
\begin{equation*}
q_\mathrm{SISDNP}=\left(\arg \min_z \frac{1}{2} \|\db_\eta-\Fb(z) \|^2\right) - \mathbb{E}(\eta).
\end{equation*}
In other words, here we apply our favorite inversion algorithm to reconstruct $z=q+\eta$ and then we subtract the expectation of $\eta$.

b) {\bf Multiple Inversion Single Data No Prior (MISDNP): } In this scenario, we assume that the data have been generated by a single realization of $\eta$ and that we can sample (or we're given samples) from the distribution of $\eta$. We solve a series of independent inverse problems for $q_s$ given a sample $\eta_s$ of $\eta$. Then we average $q_s$ and we subtract the mean of $\eta$. 
That is, we solve for $s=1,\ldots,N_s$ the problems
\begin{equation*}
q_s=\argmin_q  \frac{1}{2} \|\db_\eta-\Fb(q+\eta_s) \|^2.
\end{equation*}
Then,
\begin{equation*}
q_\mathrm{MISDNP}=\frac{1}{N_s}\sum_{s=1}^{N_s}q_s - \mathbb{E}(\eta).
\end{equation*}

c) {\bf Multiple Inversion Multiple Data No Prior (MIMDNP): } This scenario is very different in that we assume that we are given multiple measurements $\db_s$, $s=1,\ldots,N_s$, each for different (unknown) $\eta_s$ and only the expected value $\mathbb{E}(\eta)$. To obtain the approximation $q_\mathrm{MISDNP}$ of the scatterer we first solve  $s=1,\ldots,N_s$  independent inverse  problems
\begin{equation*}
q_s=\argmin_q \frac{1}{2} \|\db_s-\Fb(q) \|^2.
\end{equation*}
Then, the solution is obtained by 
\begin{equation*}
q_\mathrm{MIMDNP}=\frac{1}{N_s}\sum_{s=1}^{N_s}q_s - \mathbb{E}(\eta).
\end{equation*}

d) {\bf Single Inversion Multiple Data No Prior (SIMDNP): } The scenario is exactly the same as for problem (c). But now we first average the data and we solve a single inverse problem.  That is, we use  $\db_s$, $s=1,\ldots,N_s$, and the expected value $\mathbb{E}(\eta)$ to obtain the approximation $q_\mathrm{SIMDNP}$ of the scatterer by solving 
\begin{equation*}
q_\mathrm{SIMDNP}=\argmin_q \frac{1}{2} \left\lVert\frac{1}{N_s}\sum_{s=1}^{N_s} \db_s -\Fb(q) \right\rVert^2 - \mathbb{E}(\eta).
\end{equation*}
Method (d) is much cheaper than method (c), but is it effective?

e) {\bf Algorithm Single Inversion Single Data with Prior (SISDP): } In all remaining methods, we will assume that we have access to a probability distribution for $q$ and $\eta$ (e.g., empirical distributions obtained from samples).  In SISDP, we assume that we have data for a single (unknown) realization of $\eta$. We use $\db_\eta$, and the prior probability distributions $p(q)$ and $p(\eta)$ and solve for both $q$ and $\eta$:
\begin{equation*}
\min_{q,\eta} \frac{1}{2} \left\lVert \db_\eta-\Fb(q,\eta) \right\rVert^2 -\alpha \log p(\eta) -\beta \log p(q),
\end{equation*}
where the probability distributions are evaluated at the points of  measurements. The solution $q_\mathrm{SISDP}$ is the $q$ component of the solution. 

f) {\bf Algorithm Multiple Inversion Multiple Data with Prior (MIMDP): }
Finally, here we assume that we have multiple  measurements $\db_s$ from different (unknown) samples $\eta_s$,  $s=1,\ldots,N_s$, as well as the prior probability distributions $p(q)$ and $p(\eta)$. First, we solve $s=1,\ldots,N_s$ inverse problems
\begin{equation*}
\min_{q,\eta} \frac{1}{2} \|\db_s-\Fb(q,\eta) \|^2 -\alpha \log p(\eta) -\beta \log p(q),
\end{equation*}
and obtain $q_s$, the $q$-minimizer. The solution is then obtained by calculating
\begin{equation*}
q_\mathrm{MIMDP}=\frac{1}{N_s}\sum_{s=1}^{N_s} \tilde{q}_s - \mathbb{E}(\eta).
\end{equation*}

Methods (a)--(f) cover a rich spectrum of scenarios in which we may have or not access to probability distributions and we can have data for multiple realizations of $\eta$. An easy way to see the similarities and differences of the different methods is to linearize $\mathcal{F}$, and this is the route we pursue in our analysis. Then we conduct numerical experiments to showcase our results.  {\bf SISDNP} is the simplest algorithm. {\bf MISDNP} and {\bf MIMDNP} are variants of the sample average approximation method commonly used in stochastic programming~\cite{ruszczynski2003stochastic}, where as {\bf SIMDNP} attempts to construct a cheaper alternative to {\bf MIMDNP}. {\bf SISDP} and {\bf MIMDP} use standard regularized formulations that correspond to a point estimate within a Bayesian framework. The variance or the Conditional-Value-at-Risk statistics can also be used on the objective function, but we did not explore these schemes in this paper, for details we direct the interested reader to \cite{MHeinken17,DPKouri}.

In Table \ref{table:methods}, we present a summary of the methods, each with its respective objective function, the amount of scattered data needed, the information used, and the computational work as a function of the cost $C_I$  for solving \eqref{eq:invprob} and the number of inverse problems $N_s$. 
\begin{table}
{\small
\caption{Summary of inversion algorithms.}\label{table:methods}
\begin{center}

\bgroup
\def\arraystretch{1.2}
\begin{tabular}{|c|c|c|c|c|}
\hline
Method & Objective function & Inversion & Information & Work \\
\hline\hline
{\bf SISDNP}    & $\argmin_q \frac{1}{2} \|\db_\eta-\Fb(q) \|^2$                                                     & Single    & $\mathbb{E}(\eta)$ & $C_I$ \\
{\bf MISDNP}   & $\argmin_q  \frac{1}{2} \|\db_\eta-\Fb(q+\eta_s) \|^2$                                        & Multiple  & $\mathbb{E}(\eta)$ & $N_sC_I$ \\
{\bf MIMDNP}  & $\argmin_q \frac{1}{2} \|\db_s-\Fb(q) \|^2$                                                          & Multiple  & $\mathbb{E}(\eta)$ & $N_sC_I$ \\
{\bf SIMDNP}   & $\argmin_q \frac{1}{2} \|\frac{1}{N_s}\sum_{s=1}^{N_s} \db_s -\Fb(q) \|^2$       & Single    & $\mathbb{E}(\eta)$ & $C_I$ \\
{\bf SISDP}      & $\min_{q,\eta}\frac{1}{2} \|\db_\eta-\Fb(q,\eta) \|^2 -\alpha \log p(\eta) -\beta \log p(q)$ & Single    & $p(\eta)$, $p(q)$ & $C_I$ \\
{\bf MIMDP}     & $\min_{q,\eta}\frac{1}{2} \|\db_s-\Fb(q,\eta) \|^2 -\alpha \log p(\eta) -\beta \log p(q)$     & Multiple  & $p(\eta)$, $p(q)$, $\mathbb{E}(\eta)$ & $N_sC_I$ \\
\hline 
\end{tabular}
\egroup

\end{center}
}
\end{table}

\begin{figure}[!htp]
\centering
\subfigure[domain$+$background]{
\includegraphics[width=0.30\textwidth]{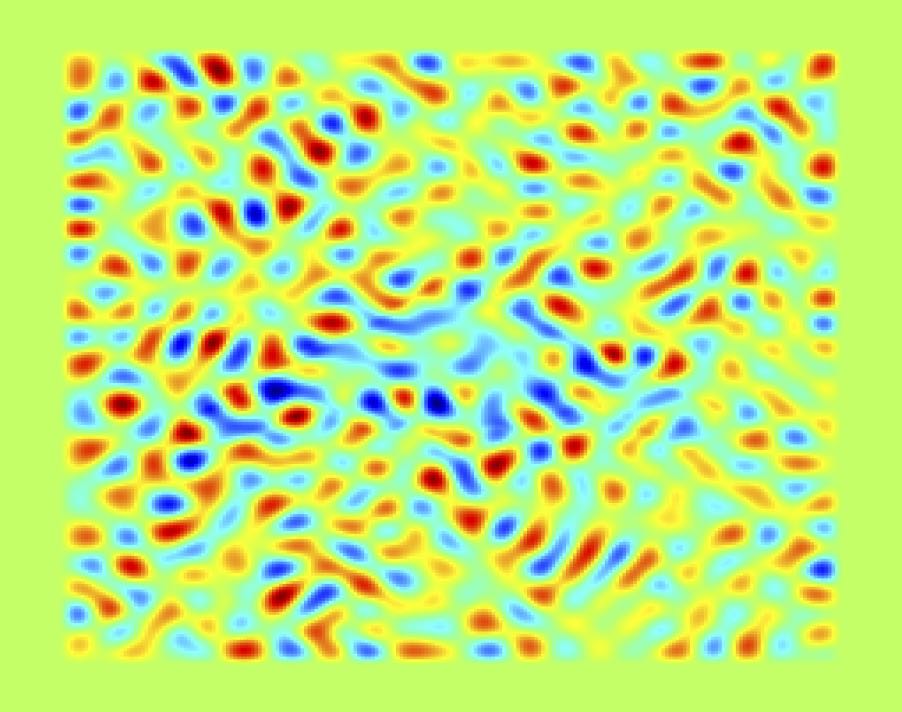}
}
\subfigure[{\bf RLA} solution]{
\includegraphics[width=0.30\textwidth]{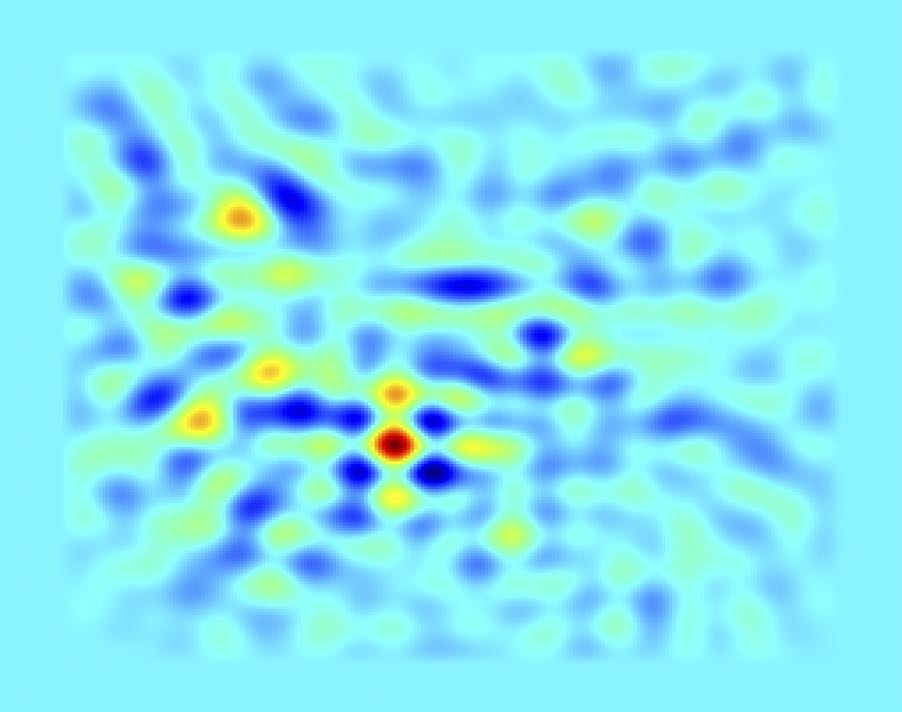}
}
\subfigure[{\bf SISDP} solution]{
\includegraphics[width=0.30\textwidth]{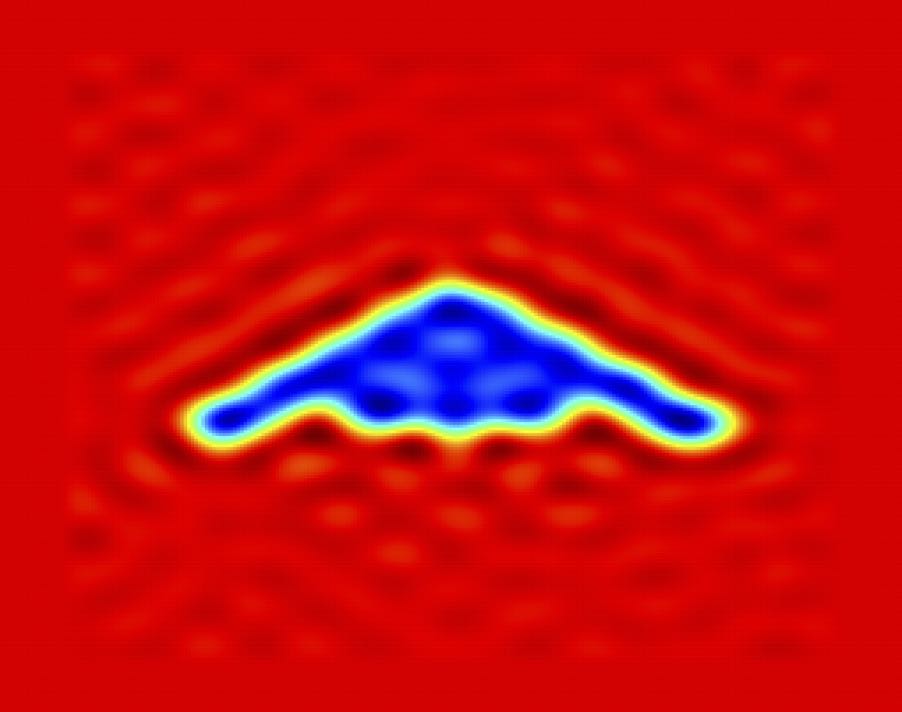}
}

\caption{Example of reconstruction of a scatterer embedded in a noisy unknown medium. In this scenario we have data from a single realization of $\eta$ and we assume we have priors for both $q$ and $\eta$. We make no assumptions on the magnitude of $\eta$ relative to $q$. In {\bf (a)}, we depict ground truth $q+\eta$.  In {\bf (b)}, we simply solve \eqref{eq:invprob} for $q+\eta$ using the RLA algorithm~\cite{Borges2017}. In (c) we use the method {\bf SISDP}. We depict the reconstructed $q$, which turns out to be very close to the ground truth $q$.}\label{pic:teaser_pic}
\end{figure}

{\bf Contributions:} We present six methods for the reconstruction of a scatterer $q$ in the presence of a random background medium $\eta$. These methods use different amounts of data and information regarding $q$ and $\eta$. We have the following remarks regarding those methods:
\begin{itemize}
\item We show that given all the information necessary and all the samples necessary, asymptotically, the methods {\bf SISDNP} and {\bf MISDNP} are equivalent and produce similar results; 
\item Considering methods that use data from only one realization of the background medium, we show that if the probability distributions of $q$ and $\eta$ are distinct enough, we obtain better results using {\bf SISDP} instead of {\bf SISDNP} or {\bf MISDNP}; 
\item When the probability distributions of $q$ and $\eta$ are not distinct enough to provide good reconstructions, it is necessary to have more data scattered off of $q$ in the presence of more realizations of $\eta$. This explains why {\bf MIMDP} has better results than {\bf SISDP}; 
\item Regarding methods that use data generated by multiple realizations of $\eta$ but no prior information, {\bf SIMDNP} is much faster than {\bf MIMDNP} and provides relatively accurate reconstructions of $q$; however, the latter method provides more accurate reconstructions at higher frequencies;
\item We show that {\bf SIMDNP} is equivalent to solving the problem using the Born approximation;
\item Our numerical examples show that, as expected, {\bf MIMDP} produce better results than {\bf MIMDNP}, since it uses more information about the probability distributions of $q$ and $\eta$; and
\item Regarding the methods that do not have prior information about the probability distribution of $q$ and $\eta$, {\bf MIMDNP} provides more accurate reconstructions of $q$ than {\bf SISDNP} or {\bf MISDNP}. This result is expected since {\bf MIMDNP} requires more data from several realizations of $\eta$.
\end{itemize}
All the inversions are done using the \ipoint{recursive linearization algorithm} (RLA) presented in~\cite{Borges2017}.  In Figure \ref{pic:teaser_pic}, we present an example of the reconstruction using SISDP. 

{\bf Limitations:} Our approach has the following limitations:
\begin{itemize}
\item Due to the lack of information in the provided data measurements and probability distributions of $\eta$ and $q$, we are not able to obtain accurate reconstructions of $q$ (at least not without contamination from $\eta$) using {\bf SISDNP} and {\bf MISDNP};
\item {\bf SIMDNP} is much faster than {\bf MIMDNP}, but {\bf SIMDNP} loses accuracy at higher frequencies, because the forward operator becomes more nonlinear;
\item {\bf SIMDNP} loses accuracy as the variance of $\eta$ increases, this is also due to the fact that the forward operator is nonlinear for larger perturbations of the domain; 
\item If the probability distributions of $q$ and $\eta$ are similar, {\bf SISDP} does not provide good approximations of the scatterer;
\item We have not considered the case of noisy data, which is typical in practice, but we do not expect any different results  other than a deterioration on the quality of the reconstruction. This is because we always regularize, either by discretization or by using some prior information. Our reconstructions are stable;
\item A result,  which is perhaps possible but missing here is to directly connect the overlap of the priors for $q$ and $\eta$ to the quality of the reconstruction. We provide empirical results but not a priori quantitative statement; and
\item There exist alternative formulations that we have not considered here,  e.g., either penalizing the variance of the posterior of the reconstructed scatterer or imposing chance constraints as a function of the reconstructed scatterer. 
\end{itemize}

{\bf Related Work:} The topic of deterministic inverse scattering problems has been extensively studied \cite{Colton,BaoLi2015,kaipio2010,stuart2010inverse}. 
There is less work on stochastic inverse scattering problems. Some work has been done regarding inverse scattering random source problems, which identify random sources from scattered data \cite{Lirandom2017,Baorandom2014,Lirandom2011,Baorandom2010,Baorandom2016,Baorandom2013,Helin2014}. Regarding the problem of recovering a domain in the presence of cluttered environments there has been extensive work in coherent interferometry~\cite{BorceaPapanicolau2005,BorceaParanicolau2006,BorceaPapanicolau2011}, and the work in the reconstruction of a domain that has very different features than the random cluttered background environment  in ~\cite{alberti2017,Alberti2017319,solna2012,bal2009physics}.  Most of the existing work assumes certain statistics on the background medium, for example small point scatterers. In our work, we make no assumptions other than smoothness and we present results in which the $L^2$-norm of the background medium is significantly larger than the target scatterer.

{\bf Article Outline:} In Section~\ref{s:forward}, we briefly describe the forward operator for the scattering medium problem and the solution of the inverse scattering medium problem. In Section~\ref{s:noprior}, we discuss the reconstruction of the scatterer of the domain in the presence of a background medium using no prior knowledge of the probability distribution of the background medium and the domain, while, in Section~\ref{s:prior}, we  consider the reconstruction of the scatterer when prior knowledge of the probability distributions is given. In Section~\ref{s:results}, we present the results of numerical experiments using the proposed methods. Concluding remarks and a short discussion of future work appear in Section~\ref{s:conclusion}.

\section{Forward and inverse scattering medium problem}\label{s:forward}
In this section, we summarize the mathematical and numerical tools used to solve the forward and inverse scattering medium problems without consideration to the presence of the background medium. Only the basic formulation and the algorithms used in our experiments are presented. For more details see Appendix A.

\subsection{Forward scattering medium problem}
In the forward scattering problem the objective is to find $u^{\emph{scat}}$, the field scattered off of the scatterer $q(\xb)$, by an incoming plane wave function $u^{\emph{inc}}(\xb)=\exp(i k \, \xb \cdot \theta)$, where $k$ is the wavenumber and $\theta$ is the incidence direction. The scattered field is the solution of the free space Helmholtz equation:
\begin{equation}
\Delta u^{\emph{scat}}(\xb) + k^2 (1+q(\xb)) u^{\text{\emph{scat}}}(\xb)
= -k^2 q(\xb) u^{\text{\emph{inc}}}(\xb),
\label{eq:fscat}
\end{equation}
where $u^{\emph{scat}}$ must satisfy the Sommerfeld radiation condition
\begin{equation}
\lim_{r \rightarrow \infty}
\sqrt{r} \, \left( \frac{\partial u^{\emph{scat}}}{\partial r} -
ik u^{\emph{scat}} \right) = 0,
\label{eq:somerrad}
\end{equation}
where $r=\|\xb\|$. 

Several authors have presented numerical methods for the solution of \eqref{eq:fscat} \cite{BorgesLS,ambikasaran2013fastdirect, borm2003hierarchical,borm2003introduction,chandrasekaran2006fast1,chen2002fast,corona2015,ifmm_darve,hackbusch2001introduction,ho2012fast,martinsson2013,xia2010fast,leonardo2016fast}. In this article, we use the HPS fast direct solver for the Helmholtz equation \cite{Gillman2014}. The computational cost of this solver is $O(N^{3/2})$ for the factorization step and $O(N)$ for each new right-hand-side, where $N$ is the number of points used to discretize the domain. We present an outline of the solver in Appendix A, but for a more detailed description, we direct the reader to \cite{Gillman2014}. The details of the forward solver are not critical for our discussion. The most important aspect is the way that $q$ is discretized, since in several of our methods, we regularize by coarsening the discretization. 

\begin{figure}[!htp]
\centering
\subfigure[Scatterer]{
\includegraphics[width=0.25\textwidth]{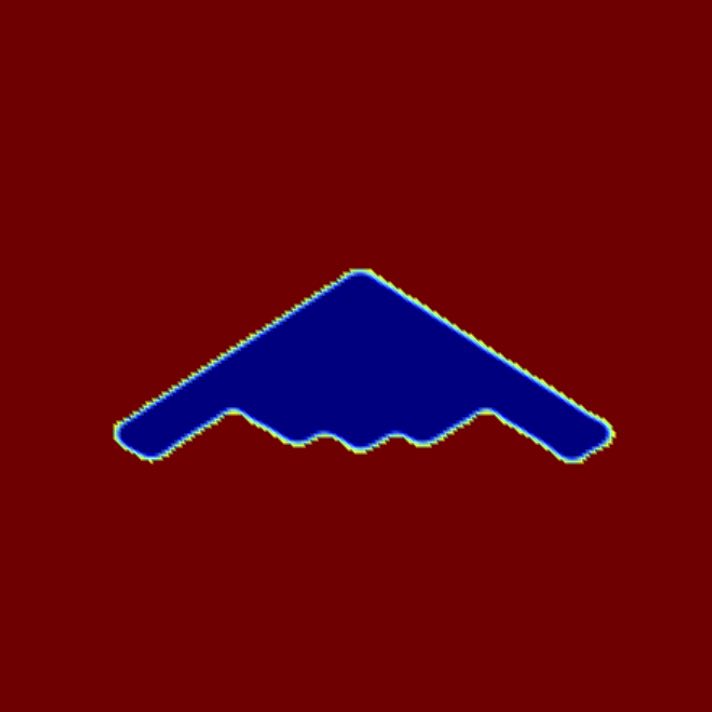}\label{pic:forward_bomber_nonoise}
}
\subfigure[$u^\emph{scat}$ for $\theta=(1,0)$]{
\includegraphics[width=0.25\textwidth]{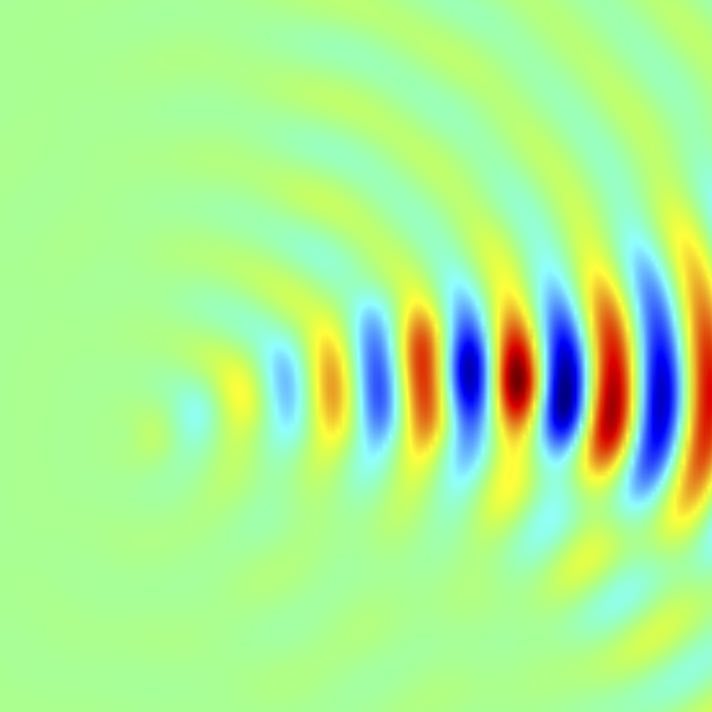}\label{pic:forward_bomber_nonoise_d0}
}
\subfigure[$u^\emph{scat}$ for $\theta=(0,1)$]{
\includegraphics[width=0.25\textwidth]{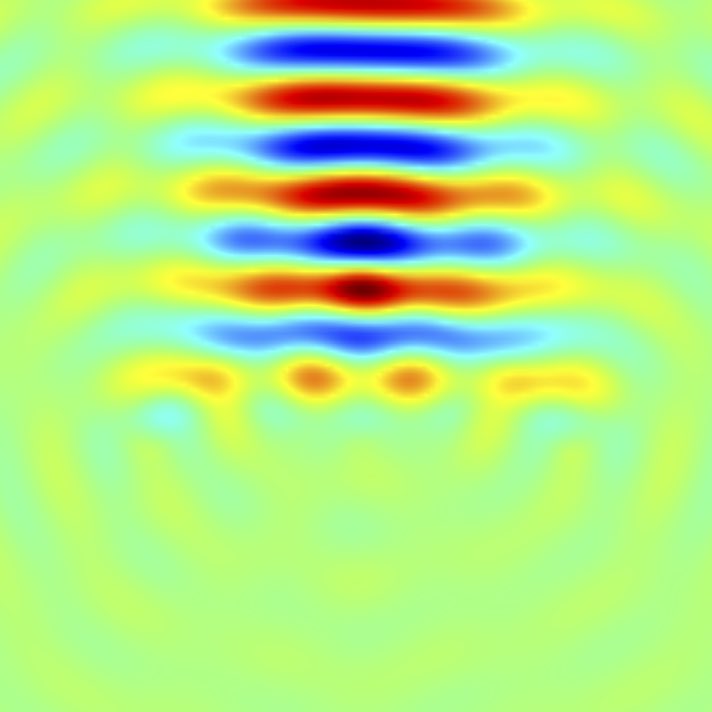}\label{pic:forward_bomber_nonoise_d90}
}

\subfigure[Scatterer]{
\includegraphics[width=0.25\textwidth]{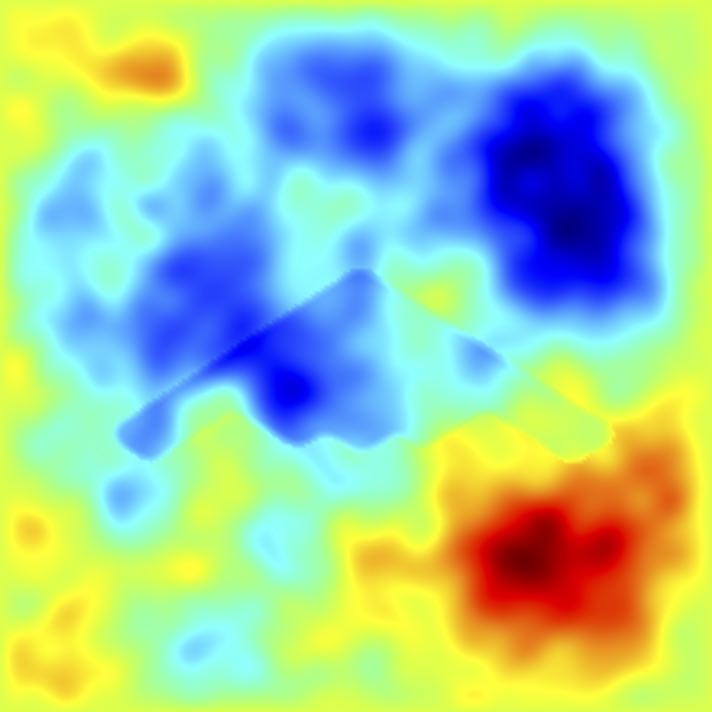}\label{pic:forward_bomber_noise}
}
\subfigure[$u^\emph{scat}$ for $\theta=(1,0)$]{
\includegraphics[width=0.25\textwidth]{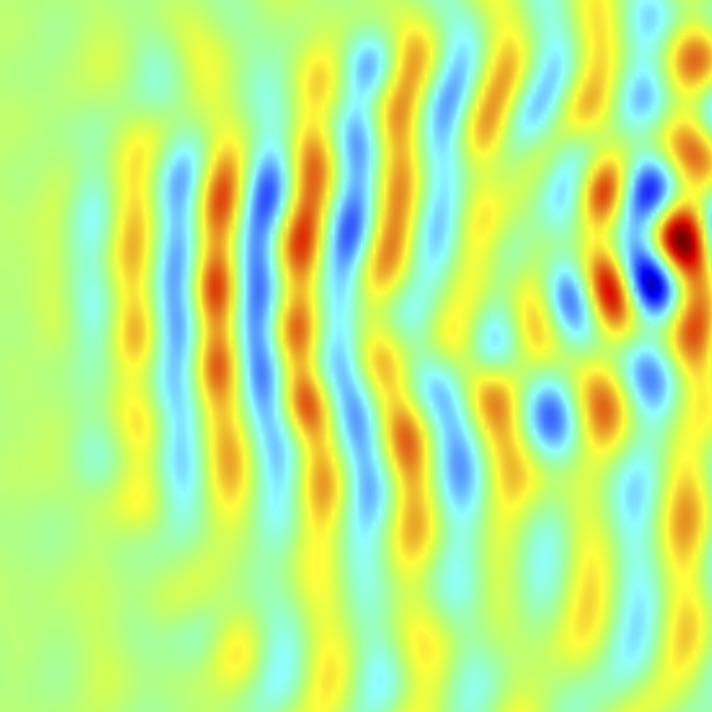}\label{pic:forward_bomber_noise_d0}
}
\subfigure[$u^\emph{scat}$ for $\theta=(0,1)$]{
\includegraphics[width=0.25\textwidth]{./forward_bomber_noise_d0_k15.jpg}\label{pic:forward_bomber_noise_d90}
}

\subfigure[Scatterer]{
\includegraphics[width=0.25\textwidth]{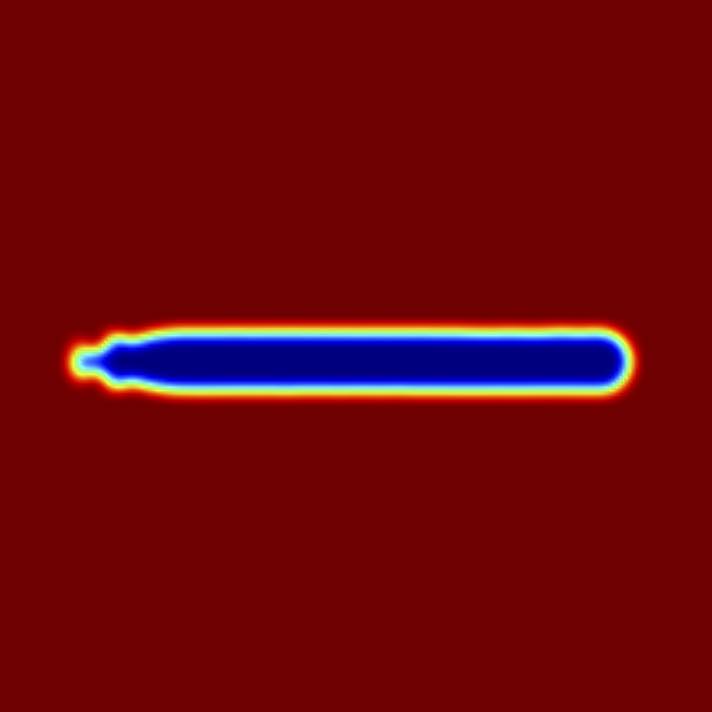}\label{pic:forward_submarine_nonoise}
}
\subfigure[$u^\emph{scat}$ for $\theta=(1,0)$]{
\includegraphics[width=0.25\textwidth]{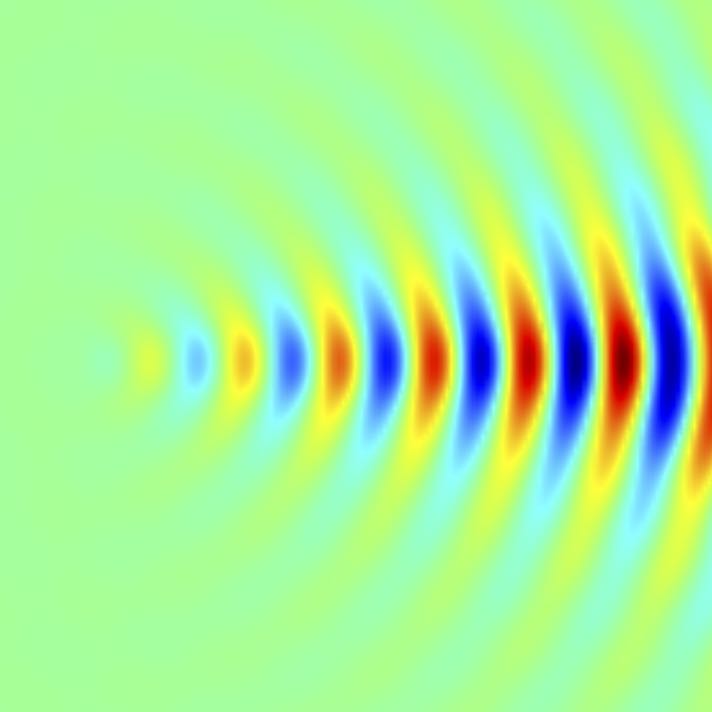}\label{pic:forward_submarine_nonoise_d0}
}
\subfigure[$u^\emph{scat}$ for $\theta=(0,1)$]{
\includegraphics[width=0.25\textwidth]{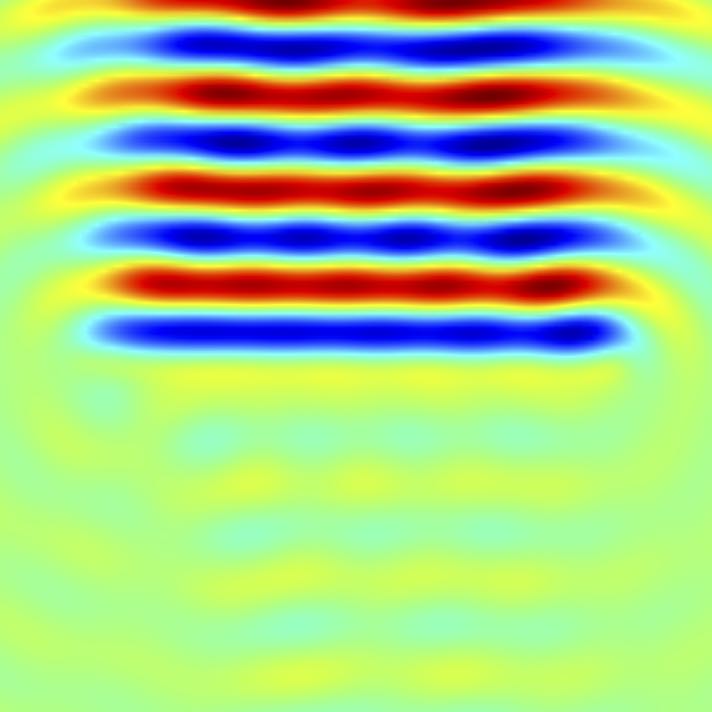}\label{pic:forward_submarine_nonoise_d90}
}

\subfigure[Scatterer]{
\includegraphics[width=0.25\textwidth]{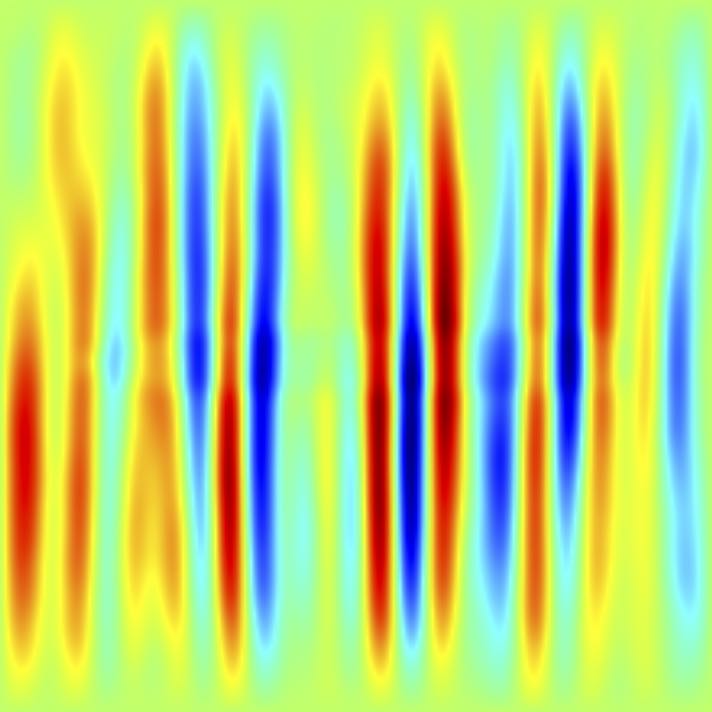}\label{pic:forward_submarine_noise}
}
\subfigure[$u^\emph{scat}$ for $\theta=(1,0)$]{
\includegraphics[width=0.25\textwidth]{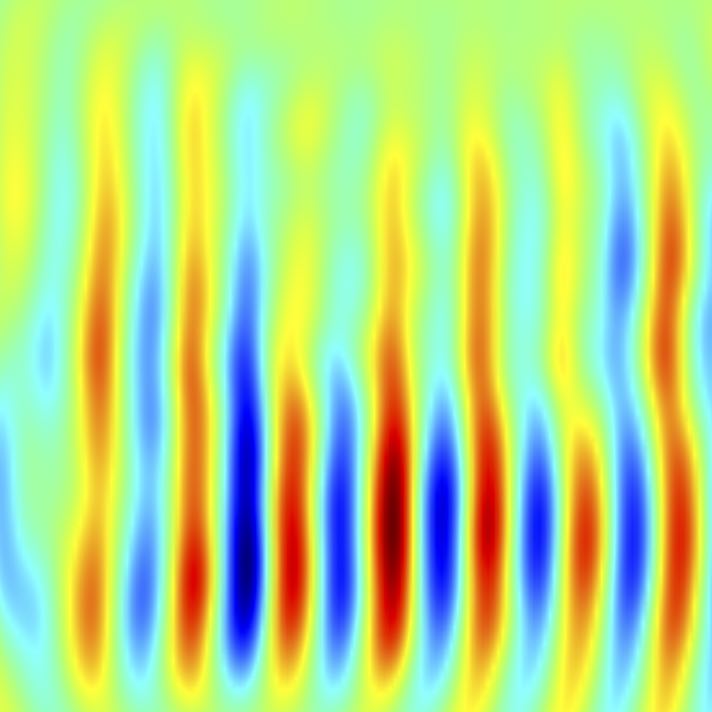}\label{pic:forward_submarine_noise_d0}
}
\subfigure[$u^\emph{scat}$ for $\theta=(0,1)$]{
\includegraphics[width=0.25\textwidth]{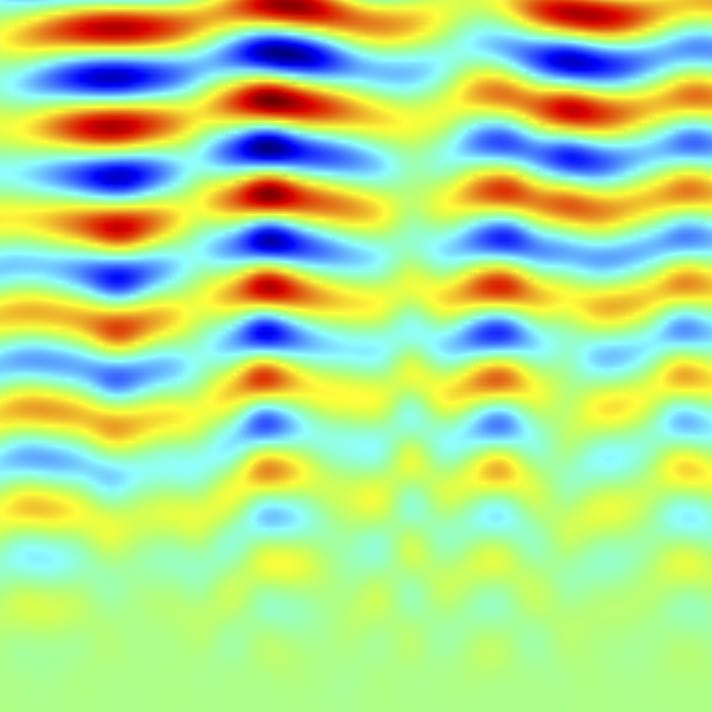}\label{pic:forward_submarine_noise_d90}
}

\caption{In each row, we present the scatterer of the medium followed by the fields scattered off of the scatterer by an incident plane wave with incident directions $\theta=(1,0)$ and $(0,1)$. }\label{pic:forward_field}
\end{figure}

In Figure \ref{pic:forward_field}, we present examples of the scattered field for several scatterer geometries. In \ref{pic:forward_bomber_nonoise_d0} and \ref{pic:forward_bomber_nonoise_d90}, we have the field scattered off of the domain represented in \ref{pic:forward_bomber_nonoise} from plane waves with $k=15$ and incidence directions $\theta=(1,0)$ and $(0,1)$, respectively.  In \ref{pic:forward_bomber_noise_d0} and \ref{pic:forward_bomber_noise_d90}, we have the field scattered off of the domain represented in \ref{pic:forward_bomber_noise} from plane waves with $k=15$ and incidence directions $\theta=(1,0)$ and $(0,1)$, respectively. In \ref{pic:forward_submarine_nonoise_d0} and \ref{pic:forward_submarine_nonoise_d90}, we have the field scattered off of the domain represented in \ref{pic:forward_submarine_nonoise} from plane waves with $k=15$ and incidence directions $\theta=(1,0)$ and $(0,1)$, respectively.  In \ref{pic:forward_submarine_noise_d0} and \ref{pic:forward_submarine_noise_d90}, we have the field scattered off of the domain represented in \ref{pic:forward_submarine_noise} from plane waves with $k=15$ and incidence directions $\theta=(1,0)$ and $(0,1)$, respectively.

\subsection{Inverse scattering medium problem}
The inverse scattering medium problem can be characterized as follows:
\begin{prob}\label{problem:invscat}
Given the measurements $\db(\xb_p)$, at the points $\xb_p$, $p=1,\ldots,N_p$ uniformly distributed on the circle $\partial B$, of the far field produced by the scattering off of $q$ of an incident plane wave $u^{\emph{inc}}=\exp(ik\xb\cdot\theta)$, recover the scatterer $q$.
\end{prob}

We solve Problem \ref{problem:invscat} by minimizing the objective function
\begin{equation}
f(q)=\frac{1}{2}\|\db-\Fb(q)\|^2,
\label{eq:obfunc}
\end{equation}
using the Gauss-Newton method. We start with an initial guess $q_0$ and iteratively update the domain $q_{n}=q_{n-1}+\delta q$, $n\in\mathbb{N}$, by solving the linear system
\begin{equation}
\Jb \, \delta q =\db - {\bf F}(q_{n-1}),
\label{linsingfreqsys}
\end{equation}
where $\Jb$ is the evaluation of the operator $\mathcal{J}$ at the measurement points, and $\mathcal{J}$ is the Fr\'{e}chet derivative of the operator $\mathcal{F}$ at $q_{n-1}$, as described in Theorem \ref{frechet_thm}: 

\begin{theorem} \label{frechet_thm} 
Let $u^\emph{inc}$ be an incoming field and let $u = u^\emph{inc} + u^\emph{scat}$ denote the total field solving the scattering problem
\begin{equation*}
\Delta u(\xb) +k^2(1+q(\xb)) u(\xb) = 0,
\end{equation*}
where $u^\emph{scat}$ satisfies the Sommerfeld radiation condition. Let $\delta q$ be a perturbation of $q$ and let $\mathcal{F}(q)$ denote the forward scattering operator. Then
\begin{equation*}
\mathcal{J}\,  \delta q(\xb) = v^\emph{scat}(\xb)
\end{equation*}
where $v^\emph{scat}(\xb)$ is the scattered field $v(\xb)$ calculated at $\partial \mathcal{B}$ and $v(\xb)$ denotes the solution to 
\begin{equation*}
\Delta v(\xb) +k^2(1+q(\xb)) v(\xb) = k^2 \, \delta q(\xb) \, u(\xb)
\end{equation*}
satisfying the Sommerfeld radiation condition.
\end{theorem}

The proof of this theorem may be found in \cite{Colton}

\begin{remark}
The extension of the Gauss-Newton method using data from $N_\theta$ impinging incident plane waves is straightforward. Consider that $\db(\theta_j)$ and $\Fb_{\theta_j}$ are the data measurements and forward operator referent to the plane wave with incident direction $\theta_j$, $j=1,\ldots,N_\theta$. The objective function becomes
\begin{equation*}
f(q)=\frac{1}{2}\sum_{j=1}^{N_\theta}\|\db(\theta_j)-\Fb_{\theta_j}(q)\|^2.
\end{equation*}
From this point forward, we will suppress the notation indicating the direction of the incident plane wave when referring to the measurements and the forward operator. We consider that the operators and measurements are being used for multiple incoming waves.
\end{remark}

The system \eqref{linsingfreqsys} is ill posed and its conditioning depends on several factors, such as $k$, the incidence directions of the incident waves, the initial guess, and the representation of $\delta q$. By Heisenberg's uncertainty principle, at wavenumber $k$, only $O(k^2)$ independent stable measurements can be made at finite precision, and those measurements only convey information regarding the lower frequencies of the domain. Using this fact, we approximate $q$ by the function $q_k(\xb)$ with ${\rm supp} (q_k) \in  [-\frac{\pi}{2},\frac{\pi}{2}]^2$, defined as follows:
\begin{equation}
q_k(\xb=(x_1,x_2))=
\sum_{m_1,m_2=1}^{M(k)} q_{m_1,m_2} B_{m_1,m_2}(x_1,x_2),
\label{eq:sineseries}
\end{equation}
where $B_{m_1,m_2}(x_1,x_2)=\sin\left( m_1 \left(x_1 + \frac{\pi}{2} \right) \right) \sin\left( m_2 \left(x_2 + \frac{\pi}{2} \right) \right)$ and the maximum frequency $M(k)$ depends on the wavenumber $k$. Representing $q$ using $q_k$ is equivalent to projecting $q$ onto a sine series with a fixed number of modes, that is, regularizing the solution of the inverse problem by discretization. Denoting this projection onto the sine modes by $\mathcal{P}_{M(k)}$ \eqref{linsingfreqsys} yields:
\begin{equation*}
\Jb \mathcal{P}_{M(k)} \delta q= \db - \Fb(q_{n}).
\end{equation*}
To ease the notation, we henceforth use $\delta q$ for both the coefficients of its approximation in the sine series and its values at points, and we suppress the projection operator $\mathcal{P}_{M(k)}$, so when we see $\Jb$, it should be considered as $\Jb \mathcal{P}_{M(k)}$. We provide specific notation when it is necessary to avoid confusion.

At low frequencies, the inverse scattering problem is uniquely solvable; however, it presents poor stability, meaning that it is difficult to obtain high resolution of the contrast function. On the other hand, at higher frequencies, the objective function presents multiple minima but is very stable. This trade-off between frequency and stability of the problem forms the basis of the \ipoint{Recursive Linearization Algorithm}. The {\bf RLA} uses standard frequency continuation to solve a sequence of inverse single-frequency scattering problems at increasing frequencies, using the solution of each problem as the initial guess for the subsequent problem. Of course, in a real application certain frequencies may not be available but application-specific details are beyond the scope of this paper. 

\begin{algorithm}
\caption{Recursive Linearization Algorithm with Gauss-Newton method (RLA).}
\label{alg:rla}
\begin{algorithmic}[1]
\STATE{{\bf Input:} data $\db(k_j)$ for $j=1,\ldots,Q$ with $k_1  < \dots < k_Q$, initial guess $q_0$, tolerances $\epsilon_1$,$\epsilon_2$ and maximum number of iterations $N_{it}$.}
\FOR{$j=1,\ldots,Q$}
\STATE{Set $q\coloneqq q_{j-1}$, $\delta q\coloneqq 0$ and $it\coloneqq 0$.}
\WHILE{$\|\db(k)-\Fb(q)\|\geq\epsilon_1$ and $it<N_{it}$ and $\delta q\geq\epsilon_2$} 
\STATE{Solve $\Jb\, \delta q=\db(k)-\Fb(q)$}
\STATE{Update $q\leftarrow q+\delta q$}
\STATE{Update $it\leftarrow it+1$}
\ENDWHILE
\STATE{Set $q_j\coloneqq q$.}
\ENDFOR
\end{algorithmic}
\end{algorithm}

\begin{remark}
In this article, the initial guess provided for the inverse problem is henceforth assumed to be in the basin of attraction of the Gauss-Newton method.
\end{remark}

\section{Inverse random medium scattering problem with no prior knowledge}\label{s:noprior}
In this section, we consider the problem of reconstructing the scatterer $q^\ast$ from measurements of the scattered field in the presence of the background noisy medium $\eta^\ast$ with no prior knowledge of the probability distribution of $q^\ast$ and $\eta^\ast$. We consider two cases. In the first case, the data measurements are generated using one realization of the background medium, while in the second case, we use data generate by multiple realizations.

\subsection{Single realization of the background medium}
In this subsection, we present two algorithms for the solution of the following problem:
\begin{prob} \label{prob:SDNP}
Given the measurements $\db_{\eta^\ast}$ of the field scattered off of the domain $q^\ast+\eta^\ast$, and given $\mathbb{E}(\eta)$, find an approximation for the unknown scatterer $q^\ast$.
\end{prob}

In the first algorithm, called {\bf SISDNP}, we use the {\bf RLA} with the Gauss-Newton method to solve the problem
\begin{equation*}
\tilde{q}=\argmin_q \frac{1}{2}\|\db_{\eta^\ast} - \Fb(q)\|^2.
\end{equation*}

At step $j$ of the Gauss-Newton method, we solve for $\delta q$ the equation
\begin{equation}
\Jb \, \delta q =\db_{\eta^\ast} - {\bf F}(q_{j-1}),
\end{equation}
after which we update the domain according to $q_j=q_{j-1}+\delta q$. 

The solution of the method is given by
\begin{equation*}
q_\mathrm{SISDNP}=\tilde{q}-\mathbb{E}(\eta).
\end{equation*}
The computational complexity of this algorithm is the same as that of applying the {\bf RLA} with the Gauss-Newton method once. 

\begin{algorithm}
\caption{SISDNP}
\label{alg:SISDNP}
\begin{algorithmic}[1]
\STATE{{\bf Input:} data $\db_{\eta^\ast}(k_j)$ for $j=1,\ldots,Q$ with $k_1  < \dots < k_Q$ and $\mathbb{E}(\eta)$.}
\STATE{Use the {\bf RLA} with Gauss-Newton method to solve $\tilde{q}=\argmin_q \frac{1}{2}\|\db_{\eta^\ast} - \Fb(q) \|^2$.}
\STATE{Calculate $q_\mathrm{SISDNP}=\tilde{q}-\mathbb{E}(\eta)$.}
\end{algorithmic}
\end{algorithm}

The second algorithm, {\bf MISDNP}, assumes that $N_s$ realizations of the background medium $\eta_s$, $s=1,\ldots,N_s$, are available. For each realization $\eta_s$, we use the {\bf RLA} with the Gauss-Newton method to solve
\begin{equation*}
\tilde{q}_s=\argmin_q \frac{1}{2}\|\db_{\eta^\ast} - \Fb(q+\eta_s)\|^2.
\end{equation*}

At the $j^\text{th}$ iteration of the Gauss-Newton method, we solve for $\delta q$ the equation
\begin{equation}
\Jb \, \delta q =\db_{\eta^\ast} - {\bf F}(q_{j-1}+\eta_s),
\end{equation}
after which we update the domain according to $q_j=q_{j-1}+\delta q$. 

The solution obtained using {\bf MISDNP} is given by
\begin{equation*}
q_\mathrm{MISDNP}=\frac{1}{N_s}\sum_{s=1}^{N_s}\tilde{q}_s-\mathbb{E}(\eta).
\end{equation*}
The computational complexity of this algorithm is equal to $N_s$ times the complexity of {\bf RLA} with the Gauss-Newton method. 

\begin{algorithm}
\caption{MISDNP}
\label{alg:MISDNP}
\begin{algorithmic}[1]
\STATE{{\bf Input:} data $\db_{\eta^\ast}(k_j)$ for $j=1,\ldots,Q$ with $k_1  < \dots < k_Q$, $\eta_s$, with $s=1,\ldots,N_s$ and $\mathbb{E}(\eta)$.}
\FOR{$s=1,\ldots,N_s$}
\STATE{Use the {\bf RLA} with Gauss-Newton method to solve $\tilde{q}_s=\argmin_q \frac{1}{2}\|\db_{\eta^\ast} - \Fb(q+\eta_s) \|^2$.}
\ENDFOR
\STATE{Calculate $q_\mathrm{MISDNP}=\sum_{s=1}^{N_s}\tilde{q}_s-\mathbb{E}(\eta)$.}
\end{algorithmic}
\end{algorithm}

Under certain (strong) assumptions in Lemma \ref{thm:SDNP} and in Example \ref{example:SDNP} we show that the solutions obtained by {\bf MISDNP} and {\bf SISDNP} are equivallent and are approximations of $q^\ast+\eta^\ast-\mathbb{E}(\eta)$, and not $q^\ast$. 

\begin{lemma} \label{thm:SDNP}
  Assume that the nonlinear least-squares problem \ref{eq:invprob}  is strictly convex  for the discretized $q$ and that  $\Fb$ is differentiable with a full-rank Jacobian $\Jb$. Moreover, let $\db_{\eta^\ast}=\Fb(q^\ast+\eta^\ast)$ be the data measurements with no noise. That is, $\db_\eta^\ast$ is in the range of $\Fb$. Further assume that we know $\mathbb{E}(\eta)$ exactly and that all numerical calculations are in exact arithmetic.
  
Let 
\begin{equation*}
\tilde{q}_1=\argmin_q\frac{1}{2}\|\db_{\eta^\ast}-\Fb(q)\|^2-\mathbb{E}(\eta),
\end{equation*}
and
\begin{equation*}
\tilde{q}_2=\frac{1}{N_s}\sum_{s=1}^{N_s}\argmin_q\frac{1}{2}\|\db_{\eta^\ast}-\Fb(q+\eta_s)\|^2-\mathbb{E}(\eta). 
\end{equation*}
Then 
\begin{equation*}
q_1 = q_2 = q^\ast+\eta^\ast-\mathbb{E}(\eta).
\end{equation*}
\end{lemma}
{\bf Proof:}
Since $\db_{\eta^\ast}=\Fb(q^\ast+\eta^\ast)$ and $\Fb$ is differentiable and $\Jb^\ast\Jb$ is invertible, we can apply the Gauss-Newton method. Since the least squares functional is strictly convex
\begin{equation*}
q_{sol}=\argmin_q\frac{1}{2}\|\db_{\eta^\ast}-\Fb(q)\|
\end{equation*}
results in  $q_{sol}=q^\ast+\eta^\ast$. So we have that if we use the {\bf SISDNP} algorithm, we obtain the solution $\tilde{q}_1=q^\ast+\eta^\ast-\mathbb{E}(\eta)$.

Similarly if we solve multiple inverse problems for each realization $\eta_s$ of $\eta$,
\begin{equation*}
q_s=\argmin_q\frac{1}{2}\|\db_{\eta^\ast}-\Fb(q+\eta_s)\|,
\end{equation*}
we obtain $q_s=q^\ast+\eta^\ast-\eta_s$. Averaging over all samples, we obtain the solution for {\bf MISDNP}
\begin{equation*}
\tilde{q}_2=q^\ast+\eta^\ast-\mathbb{E}(\eta).
\end{equation*}
$~\square$

Table \ref{table:SDNP} summarizes the algorithms {\bf SISDNP} and {\bf MISDNP}. In Example \ref{example:SDNP}, we present numerical experiments exemplifying and comparing {\bf SISDNP} and {\bf MISDNP}.

\begin{table}[!htbp]
{\small
\begin{center}

\bgroup
\def\arraystretch{1.2}
\begin{tabular}{|c|c|c|c|}
\hline
Algorithm & Function & Complexity & Approximation \\
\hline\hline
 {\bf SISDNP} & $\min_q \frac{1}{2}\|\db_{\eta^\ast} - \Fb(q)\|^2$ & $C_I$ & $q^\ast+\eta^\ast-\mathbb{E}(\eta)$ \\
\hline
 {\bf MISDNP} & $\min_q \frac{1}{2}\|\db_{\eta^\ast} - \Fb(q+\eta_s)\|^2$ & $N_sC_I$ & $q^\ast+\eta^\ast-\mathbb{E}(\eta)$ \\
\hline
\end{tabular}
\egroup

\end{center}
}
\caption{Summary of algorithms {\bf SISDNP} and {\bf MISDNP}.} \label{table:SDNP}
\end{table}

Although the assumptions of \ref{thm:SIMDNP} are not in general true, assuming no noise and appropriate discretization of $q$, the Hessian of the nonlinear least squares is the Gauss-Newton Jacobian and the problem is locally convex. So, at least locally, SISDNP and MISDNP are equivalent. But of course, SISDNP is much cheaper. We conclude that for data from a single realization of $\eta$, using multiple inversions doesn't provide any additional information. Moreover, both methods result in significant error and fail to reconstruct $q$.

\subsection{Multiple realizations of the background medium}
Another scenario we consider here is the case in which we have data measurements from different realizations of the background medium $\eta$. The precise statement of the problem reads as follows.
\begin{prob} \label{prob:MDNP}
Given the measurements $\db_s$, $s=1,\ldots,N_s$, of the field scattered off of the domain $q^\ast+\eta_s^\ast$, and given $\mathbb{E}(\eta)$, find an approximation for the unknown scatterer $q^\ast$.
\end{prob}

In the first algorithm, {\bf MIMDNP}, we use the {\bf RLA} with the Gauss-Newton method to solve the problem  
\begin{equation*}
\tilde{q}_s=\argmin_q \frac{1}{2}\|\db_s - \Fb(q)\|^2.
\end{equation*}

At step $j$ of the Gauss-Newton method, we solve for $\delta q$ the equation
\begin{equation}
\Jb \, \delta q =\db_s - {\bf F}(q_{j-1}),
\end{equation}
after which we update the domain according to $q_j=q_{j-1}+\delta q$. 

The solution $q_\mathrm{MIMDNP}$ is obtained by averaging the $\tilde{q}_s$ and subtracting the average of $\eta$. The result is
\begin{equation*}
q_\mathrm{MIMDNP}=\frac{1}{N_s}\sum_{s=1}^{N_s}\tilde{q}_s - \mathbb{E}(\eta).
\label{eq:MIMDNP}
\end{equation*}

Applying the {\bf RLA} at $\db_s$ gives an approximation $\tilde{q}_s\approx q+\eta_s$, and the average provides an approximation $\frac{1}{N_s}\sum_{s=1}^{N_s}\tilde{q}_s\approx q^\ast+\mathbb{E}(\eta)$. Since $\mathbb{E}(\eta)=\frac{1}{N_s}\sum_{s=1}^{N_s}\eta_s$ is given, assuming that we have enough samples to calculate the expected value of the background noise medium, $q_\mathrm{MIMDNP}$ is the best band-limited approximation possible for the scatterer $q^\ast$.

The computational complexity of this algorithm is equal to $N_s$ times the computational complexity of the {\bf RLA} applied to each data set. 

\begin{remark}
Equation \eqref{eq:MIMDNP} is equivalent to
\begin{equation*}
q_\mathrm{MIMDNP}=\mathbb{E}_s\left(\argmin_q \frac{1}{2}\|\db_s - \Fb(q)\|^2\right)-\mathbb{E}(\eta),
\end{equation*}
where the first average is taken over the index $s$.
\end{remark}

\begin{algorithm}
\caption{MIMDNP}
\label{alg:MIMDNP}
\begin{algorithmic}[1]
\STATE{{\bf Input:} data $\db_s(k_j)$ for $j=1,\ldots,Q$ with $k_1  < \dots < k_Q$ and $s=1,\ldots,N_s$, and $\mathbb{E}(\eta)$.}
\FOR{$s=1,\ldots,N_s$}
\STATE{Use the {\bf RLA} with Gauss-Newton method to solve $\tilde{q}_s=\argmin_q \frac{1}{2}\|\db_s - \Fb(q) \|^2.$.}
\ENDFOR
\STATE{Calculate $q_\mathrm{MIMDNP}=\frac{1}{N_s}\sum_{s=1}^{N_s}\tilde{q}_s-\mathbb{E}(\eta)$.}
\end{algorithmic}
\end{algorithm}

In the second algorithm, {\bf SIMDNP}, we first calculate the average of the data measurements and subtract the field scattered by $\mathbb{E}(\eta)$:
\begin{equation*}
\bar{\db}=\frac{1}{N_s}\sum_{s=1}^{N_s}\db_s-\Fb(\mathbb{E}(\eta)).
\end{equation*}
Next, we apply the {\bf RLA} with the Gauss-Newton method to solve the problem
\begin{equation*}
\tilde{q}=\argmin_q \frac{1}{2}\|\bar{\db} - \Fb(q)\|^2.
\end{equation*}

At step $j$ of the Gauss-Newton method, we solve for $\delta q$ the equation
\begin{equation}
\Jb \, \delta q =\bar{\db} - {\bf F}(q_{j-1}),
\end{equation}
after which we update the domain according to $q_j=q_{j-1}+\delta q$. 

Finally, we calculate the solution:
\begin{equation*}
q_\mathrm{SIMDNP}=\tilde{q}-\mathbb{E}(\eta).
\end{equation*}

The computational complexity of this algorithm is equal to the computational complexity of the {\bf RLA} applied to $\bar{\db}$. 

\begin{algorithm}
\caption{SIMDNP}
\label{alg:SIMDNP}
\begin{algorithmic}[1]
\STATE{{\bf Input:} data $\db_s(k_j)$ for $j=1,\ldots,Q$ with $k_1  < \dots < k_Q$ and $s=1,\ldots,N_s$, and $\mathbb{E}(\eta)$.}
\STATE{Calculate $\bar{\db}=\frac{1}{N_s}\sum_{s=1}^{N_s}\db_s-\Fb(\mathbb{E}(\eta))$.}
\STATE{Use the {\bf RLA} with Gauss-Newton method to solve $q_\mathrm{SIMDNP}=\argmin_q \frac{1}{2}\|\bar{\db} - \Fb(q) \|^2$.}
\end{algorithmic}
\end{algorithm}

On one hand, the clear advantage of {\bf SIMDNP} over {\bf MIMDNP} is that instead of solving $N_s$ problems simultaneously, only one application of the {\bf RLA} in the average of the data is necessary. On the other hand, as we can see in Example \ref{example:MDNP}, {\bf SIMDNP} is less accurate at higher frequencies and for background medium functions with large variance. This behavior is easily explained by Lemma \ref{thm:SIMDNP}. As the wavenumber and the variance of the background medium increase the convexity of the inverse problem is lost  and the two methods are no longer equivalent.

\begin{lemma}\label{thm:SIMDNP}
Let $k^2\|q^\ast+\eta_s\|_{\infty}\ll 1$, $k^2\|\eta_s\|_{\infty}\ll 1$, $\db=\Fb(q^\ast)$, $\mathbb{E}(\eta)\equiv0$, $\db_s=\Fb(q^\ast+\eta_s)$, for $s\in \mathbb{N}$, and 
\begin{equation}
q_N=\argmin_q \frac{1}{2}\|d-\Fb(q)\|^2
\label{eq:newton}
\end{equation} 
be the solution given by the {\bf RLA} with the Gauss-Newton method. If we apply the {\bf SIMDNP} algorithm using an arbitrarily large number of samples, we have that 
\begin{equation*}
q_\mathrm{SIMDNP}\rightarrow q_N.
\end{equation*}
\end{lemma}

{\bf Proof:} To solve Problem \eqref{eq:newton}, we apply the Gauss-Newton method. Starting from the initial guess $q_0$, at the $j^\text{th}$ iteration we solve
\begin{equation}
\Jb\, \delta q=\db-\Fb(q_j)
\label{eq:newtonstep}
\end{equation}
and update the domain according to $q_{j+1}=q_j+\delta q$. As $j\rightarrow\infty$, our series converges to the Newton solution, $q_j\rightarrow q_N$.

If we use the {\bf SIMDNP} algorithm, starting from the initial guess $q_0$, at each step, we solve
\begin{equation}
\Jb\, \delta q=\bar{\db}-\Fb(q_j)
\label{eq:SIMDNPstep}
\end{equation}
and update the domain according to $q_{j+1}=q_j+\delta q$. 

Regarding the data for {\bf SIMDNP}, we have that 
\begin{equation}
\bar{\db}=\frac{1}{N_s}\sum_{s=1}^{N_s}\db_s=\frac{1}{N_s}\sum_{s=1}^{N_s}\Fb(q^\ast+\eta_s).
\label{eq:dataSIMDNP}
\end{equation}
Since $k^2\|q^\ast+\eta_s\|_{\infty}\ll 1$, $q^\ast+\eta_s$ is a small perturbation of the domain. We are in the Born approximation regime, which means that the forward operator becomes linear \cite{Colton}. Using this fact on \eqref{eq:dataSIMDNP}, we obtain a sequence
\begin{equation*}
\bar{\db}_s=\Fb(q^\ast)+\frac{1}{N_s}\sum_{s=1}^{N_s}\Fb(\eta_s).
\end{equation*}
Since $k^2\|\eta_s\|_{\infty}\ll 1$, we can use again the linearity of the forward operator and obtain
\begin{equation*}
\bar{\db}_s=\Fb(q^\ast)+\Fb(\frac{1}{N_s}\sum_{s=1}^{N_s}\eta_s).
\end{equation*}
It is clear that $\bar{\db}_s\rightarrow \Fb(q^\ast) + \Fb(\mathbb{E}(\eta))=\db$ as $s\rightarrow\infty$, since $\mathbb{E}(\eta)\equiv0$.

Using the limit of the sequence $\bar{\db}_s$ in Equation \eqref{eq:SIMDNPstep}, we get the same equation as \eqref{eq:newtonstep}. This means that both methods provide the same sequence of approximations of the domain, which converges to the same solution.$\square$

It is worth noticing, in Lemma \ref{thm:SIMDNP1}, that using the {\bf SIMDNP} algorithm is equivalent to applying the {\bf RLA} with the Gauss-Newton method in the problem
\begin{equation*}
\argmin_q\mathbb{E}\left(\frac{1}{2}\|\db-\Fb(q)\|^2\right).
\end{equation*}

\begin{lemma}\label{thm:SIMDNP1}
Suppose we have the data measurements $\db_s=\Fb(q^\ast+\eta_s^\ast)$ for $s=1,\ldots,N_s$ and without loss of generality $\mathbb{E}(\eta)=0$. If we use the Gauss-Newton method to solve the problem
\begin{equation}
\argmin_q\mathbb{E}\left(\frac{1}{2}\|\db-\Fb(q)\|^2\right)=\argmin_q \frac{1}{N_s}\sum_{s=1}^{N_s}\frac{1}{2}\|\db_s-\Fb(q)\|^2,
\label{method_2}
\end{equation}
we obtain the same result as using the Gauss-Newton method to solve the problem
\begin{equation}
\argmin_q \frac{1}{2}\left\lVert\frac{1}{N_s}\sum_{s=1}^{N_s}\db_s-\Fb(q)\right\rVert^2,
\label{method_3}
\end{equation}
when using the same initial guess $q_0$.
\end{lemma}

{\bf Proof:} Suppose we want to solve \eqref{method_3} using the Gauss-Newton method with $q_0$ as the initial guess. At the $j^\text{th}$ iteration, we solve the system
\begin{equation}
\Jb^\ast \Jb \delta q= \Jb^\ast \left( \frac{1}{N_s}\sum_{s=1}^{N_s} \db_s \right)- \Jb^\ast \Fb(q_j),
\label{normal_method_3}
\end{equation}
and update $q_{j+1}=q_{j}+\delta q$. Analogously, to solve \eqref{method_2}, we apply the Gauss-Newton method and obtain
\[
 \begin{bmatrix}
\Jb \\
\vdots \\
\Jb
\end{bmatrix} \delta q=
\begin{bmatrix}
\db_1-\Fb(q_j) \\
\vdots \\
\db_{N_s}-\Fb(q_j)
\end{bmatrix}.
\]
It is straightforward to see that the normal equations to solve the system above are the same as \eqref{normal_method_3}. Since at each step the normal equations to be solved are the same and the initial guess for both methods is the same, then both methods yield the same solution at each step.$\square$

Table \ref{table:MDNP} summarizes the {\bf MIMDNP} and {\bf SIMDNP} algorithms. In Example \ref{example:MDNP}, we present numerical experiments exemplifying and comparing the methods.

\begin{table}[!htbp]
{\small
\begin{center}

\bgroup
\def\arraystretch{1.2}
\begin{tabular}{|c|c|c|c|}
\hline
Algorithm & Function & Complexity \\
\hline\hline
 {\bf MIMDNP} & $\min_q \frac{1}{2}\|\db_s - \Fb(q)\|^2$ & $N_sC_I$ \\

\hline

 {\bf SIMDNP} & $\min_q \frac{1}{2}\|\bar{\db} - \Fb(q)\|^2$ & $C_I$ \\
\hline
\end{tabular}
\egroup

\end{center}
}
\caption{Summary of the {\bf MIMDNP} and {\bf SIMDNP} algorithms.} \label{table:MDNP}
\end{table}

\section{Inverse random medium scattering problem with priors}\label{s:prior}
Now we consider the problem of reconstructing from measurements of the scattered field the scatterer $q$ in the presence of the background noisy medium $\eta$ given the prior probability distributions of $q$ and $\eta$. We present two schemes, {\bf SISDP} (when we have data from one realization of $\eta$)  and {\bf MIMDP} (when we have data from multiple realizations of $\eta$). Unlike the previous schemes, now we will be \ipoint{solving the inverse problem for both} $\eta$ and $q$. We assume that $q$ and $\eta$ are statistically independent. Another difference with the schemes in Section~\ref{s:noprior} is that we regularize using the priors and we do not truncate the discretization. 

\subsection{Single realization of the background medium}
In this subsection, we present an algorithm for the solution of the following problem:
\begin{prob} \label{prob:SDP}
Given measurements $\db_{\eta^\ast}$ of the scattered field off $q^\ast+\eta^\ast$, and given the prior probability distributions $p(\eta)$ and $p(q)$ of the background medium and scatterer, respectively, we seek to find an approximation for the unknown scatterer $q^\ast$.
\end{prob}

In particular, we assume that the background noisy medium has the probability distribution
\begin{equation}
p(\eta)=\exp(-\eta^T\Tb_\eta \eta), 
\label{eq:prob_eta}
\end{equation}
while the scatterer has the probability distribution
\begin{equation}
p(q)=\exp(-q^T\Tb_q q),
\label{eq:prob_q}
\end{equation}
where $\Tb_\eta$ and $\Tb_q$ are the inverse covariance operators for the respective distributions. In our numerical experiments we will construct these operators from standard anisotropic smoothness priors for $q$ and $\eta$. In particular,
$\Tb_q=\nabla\cdot \mathcal{T}_q \nabla$  and $\Tb_\eta = \nabla\cdot \mathcal{T}_\eta \nabla$. Since we use the pseudo-spectral discretization for $q$ and $\eta$ these operators become diagonal. 

In {\bf SISDP}, we apply the {\bf RLA} using the Gauss-Newton method to solve 
\begin{equation*}
\min_{q,\eta}  \frac{1}{2}\|\db_{\eta^\ast} - \Fb(q+\eta)\|^2 -\frac{\alpha}{2} \log(p(\eta)) -\frac{\beta}{2} \log(p(q)) ,
\end{equation*}
where $\alpha$ and $\beta$ are regularization parameters for $\eta$ and $q$, respectively. Using the probability functions \eqref{eq:prob_eta} and \eqref{eq:prob_q} as the prior, we obtain
\begin{equation}
\min_{q,\eta}  \frac{1}{2}\|\db_{\eta^\ast} - \Fb(q+\eta)\|^2 + \frac{\alpha}{2}\|\eta\|^2_{\Tb_\eta} + \frac{\beta}{2}\|q\|^2_{\Tb_q}. 
\label{eq:functionSISDP}
\end{equation}
At step $j$, we solve the system
\begin{equation}
\begin{bmatrix} \Jb & \Jb \\ 0 & \Tb_\eta^{1/2} \\ \Tb_q^{1/2} & 0 \end{bmatrix}  \left[ \begin{array}{c} \delta q \\ \delta \eta \end{array} \right] = \left[ \begin{array}{c}\db_{\eta^\ast} - \Fb(q+\eta) \\ 0 \\ 0 \end{array} \right],
\label{eq:SISDP:system}
\end{equation}
where $\Jb$ is the Fr\'{e}chet derivative of $\Fb$ at $q+\eta$ . The updated solutions are $q\leftarrow q+\delta q$ and $\eta\leftarrow \eta+\delta \eta$.

The computational complexity of this algorithm is equal to the computational complexity of the {\bf RLA}. 

\begin{algorithm}
\caption{SISDP}
\label{alg:SISDP}
\begin{algorithmic}[1]
\STATE{{\bf Input:} data $\db(k_j)$ for $j=1,\ldots,Q$ with $k_1  < \dots < k_Q$, and the isotropy tensors $\mathcal{T}_\eta$ and $\mathcal{T}_q$.}
\STATE{Calculate $\Tb_\eta$ and $\Tb_q$ using $\mathcal{T}_\eta$ and $\mathcal{T}_q$ respectively.}
\STATE{Use the {\bf RLA} with Gauss-Newton method to solve 
$$\min_{q,\eta}  \frac{1}{2}\|\db_{\eta^\ast} - \Fb(q+\eta)\| ^2 + \frac{\alpha}{2}\|\eta\|_{\Tb_\eta}^2 + \frac{\beta}{2}\|q\|_{\Tb_\eta}^2.$$}
\STATE{Set $q_\mathrm{SISDP}\leftarrow q$.}
\end{algorithmic}
\end{algorithm}

\subsection{Multiple realizations of the background medium}
In this subsection we present an algorithm for the solution of the following problem:
\begin{prob} \label{prob:MDP}
Given measurements $\db_s$ of the field scattered off of the domain $q^\ast+\eta^\ast_s$, $s=1,\ldots,N_s$, given the prior  probability distributions $p(\eta)$ and $p(q)$ of the background medium and scatterer, respectively, and given the expected value $\mathbb{E}(\eta)$, find an approximation for the unknown scatterer $q^\ast$.
\end{prob}

In the {\bf MIMDP} algorithm, first, for each data measurement $\db_s$, we apply the {\bf SISDP} algorithm to obtain a solution, which we denote $\tilde{q}_s$. Next, we average the results and subtract the expected value of the background medium, obtaining the solution:
\begin{equation*}
q_\mathrm{MIMDP}=\frac{1}{N_s}\sum_{s=1}^{N_s}\tilde{q}_s-\mathbb{E}(\eta).
\end{equation*}

The computational complexity of this algorithm is equal to $N_s$ times the computational complexity of the {\bf RLA} applied to each data set $\db_s$. 

\begin{algorithm}
\caption{MIMDP}
\label{alg:MIMDP}
\begin{algorithmic}[1]
\STATE{{\bf Input:} data $\db(k_j)$ for $j=1,\ldots,Q$ with $k_1  < \dots < k_Q$, and the isotropy tensors $\mathcal{T}_\eta$ and $\mathcal{T}_q$.}
\STATE{Calculate $\Tb_\eta$ and $\Tb_q$ using $\mathcal{T}_\eta$ and $\mathcal{T}_q$ respectively.}
\FOR{$s=1,\ldots,N_s$}
\STATE{Use the {\bf SISDP} algorithm with input data $\db_s$, obtaining the solution $\tilde{q}_s$ for the scatterer.}
\ENDFOR
\STATE{Calculate $q_\mathrm{MIMDP}=\frac{1}{N_s}\sum_{s=1}^{N_s} \tilde{q}_s - \mathbb{E}(\eta)$.}
\end{algorithmic}
\end{algorithm}

A key issue here is the choice of $\mathcal{T}_\eta$ and $\mathcal{T}_q$ as well as the regularization parameters $\alpha$ and $\beta$. These are related to the specific application and the prior distributions for $q$ and $\eta$. In our numerical experiments, we give details in the context of specific examples. As we will see, depending on the prior probability distributions we may be able to reconstruct $q^\ast$ accurately. However, if the prior distributions are similar, {\bf SISDP} may fail to reconstruct. As the number of realizations for $\eta$ increases {\bf MIMDP} does better.

\begin{figure}[h!]
  \centering
\subfigure[similar priors]{
\begin{tikzpicture}[scale=0.6]
\begin{axis}[
  no markers, 
  domain=0:6, 
  samples=100,
  ymin=0,
  axis lines*=left, 
  every axis y label/.style={at=(current axis.above origin),anchor=south},
  every axis x label/.style={at=(current axis.right of origin),anchor=west},
  legend pos=north west,
  legend entries={$q_i$,$\eta_i$},
  height=5cm, 
  width=12cm,
  xtick=\empty, 
  ytick=\empty,
  enlargelimits=false, 
  clip=false, 
  axis on top,
  grid = major,
  hide y axis
  ]

 \addplot [very thick,cyan!50!black] {gauss(x, 3, 1)};
 \addplot [very thick,red!50!black] {gauss(x, 3, .8)};

\pgfmathsetmacro\valueA{gauss(1.5,3,1)}
\pgfmathsetmacro\valueB{gauss(1.5,3,.8)}
\pgfmathsetmacro\valueAA{gauss(3,3,1)}
\pgfmathsetmacro\valueBB{gauss(3,3,.8)}



\draw [gray] (axis cs:1.5,0) -- (axis cs:1.5,\valueA);
\draw [gray] (axis cs:1.5,0) -- (axis cs:1.5,\valueB);
\draw [gray] (axis cs:3,0) -- (axis cs:3,\valueAA);
\draw [gray] (axis cs:3,0) -- (axis cs:3,\valueBB);


\node[below] at (axis cs:1.5, 0)  {$z_i$}; 
\node[below] at (axis cs:3, 0)  {$0$}; 
\end{axis}

\end{tikzpicture}\label{fig:dist_features_a}
}
\subfigure[Diffent priors]{
\begin{tikzpicture}[scale=0.6]
\begin{axis}[
  no markers, 
  domain=0:6, 
  samples=100,
  ymin=0,
  axis lines*=left, 
  every axis y label/.style={at=(current axis.above origin),anchor=south},
  every axis x label/.style={at=(current axis.right of origin),anchor=west},
  legend pos=north west,
  legend entries={$q_i$,$\eta_i$},
  height=5cm, 
  width=12cm,
  xtick=\empty, 
  ytick=\empty,
  enlargelimits=false, 
  clip=false, 
  axis on top,
  grid = major,
  hide y axis
  ]

 \addplot [very thick,cyan!50!black] {gauss(x, 3, 2)};
 \addplot [very thick,red!50!black] {gauss(x, 3, .3)};

\pgfmathsetmacro\valueA{gauss(1.5,3,2)}
\pgfmathsetmacro\valueB{gauss(1.5,3,.3)}
\pgfmathsetmacro\valueAA{gauss(3,3,2)}
\pgfmathsetmacro\valueBB{gauss(3,3,.3)}

\draw [gray] (axis cs:1.5,0) -- (axis cs:1.5,\valueA);
\draw [gray] (axis cs:1.5,0) -- (axis cs:1.5,\valueB);
\draw [gray] (axis cs:3,0) -- (axis cs:3,\valueAA);
\draw [gray] (axis cs:3,0) -- (axis cs:3,\valueBB);

\node[below] at (axis cs:1.5, 0)  {$z_i$}; 
\node[below] at (axis cs:3, 0)  {$0$}; 
\end{axis}

\end{tikzpicture} \label{fig:dist_features_b} 
}
\caption{Probability distributions for $q$ and $\eta$ when they have: (a) similar, and (b) different priors.}\label{fig:dist_features}
\end{figure}
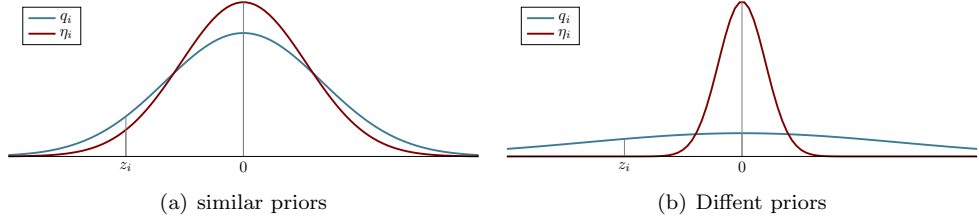

We use a simple example to discuss the interplay between prior for $\eta$ and $q$ and our ability to reconstruct $q$. We consider the linearization of the nonlinear inversion to a simple least squares problem, which we assume that is well posed for the sum $\eta + q$, denoted by $z$. Then we try to reconstruct $q$ and $\eta$ given $z$, assuming that the prior covariance operators for $\eta$ and $q$ can be diagonalized by the same basis (in our case spectral). 

So first we solve
$$
z^\ast = \argmin_z \frac{1}{2}\|\Jb z-\db\|^2
$$
where we assume that the solution is stable. 

To obtain the value of $q$ and $\eta$ from the solution $z^\ast$ we solve
\begin{equation}
\argmin_{q,\eta} \frac{1}{2}\|q+\eta-z^\ast\|^2+\frac{1}{2}\|q\|^2_{\Tb_q}+\frac{1}{2}\|\eta\|^2_{\Tb_\eta}. \label{eq:inverse_stat}
\end{equation}
Notice that here and in our experiments both $q$ and $\eta$ are normally distributed with zero mean. The extension to the non-zero mean case is straighfortward.

First, we assume that $q$ and $\eta$ are statistically independent. Further, assume that  $\Tb_q$ and $\Tb_\eta$ are diagonal with entries  $\sigma_{q_i}^{-2}$ and $\sigma_{\eta_i}^{-2}$ respectively. That is components $q_i$ and $\eta_i$ are independent with priors $N(0,\sigma_{q_i})$ and $N(0,\sigma_{\eta_i})$  respectively. Therefore to determine $q$ and $\eta$ we solve  for each component $i$,
\begin{equation}
\argmin_{q_i,\eta_i} \frac{1}{2}(q_i+\eta_i-z_i)^2+\frac{1}{2\sigma_{q_i}^2}q_i^2+\frac{1}{2\sigma_{\eta_i}^2}\eta_i^2, \label{eq:comp_inverse}
\end{equation}
whose analytic solution is given by
\begin{equation}
q_i=\frac{\sigma_{q_i}^2 z_i}{1+\sigma_{q_i}^2+\sigma_{\eta_i}^2} \label{eq:q_i_sol}
\end{equation}
and
\begin{equation}
\eta_i=\frac{\sigma_{\eta_i}^2 z_i}{1+\sigma_{q_i}^2+\sigma_{\eta_i}^2}. \label{eq:eta_i_sol}
\end{equation}

Equations \ref{eq:q_i_sol} and \ref{eq:eta_i_sol} give us a good idea of how the spectral decomposition of the priors relate to the reconstruction of the features of $q$. As we can see, if $\sigma_{\eta_i} \gg \sigma_{q_i}$ then $q_i \rightarrow 0$ and $\eta_i \rightarrow z_i$. If $\sigma_{\eta_i} \ll \sigma_{q_i}$ then $q_i \rightarrow z_i$ and $\eta_i \rightarrow 0$. Figure \ref{fig:dist_features} illustrates the priors for  $q_i$ and $\eta_i$. When they have similar priors it is harder to reconstruct $q$ accurately. If we use spectral truncation for the regularization of $q$, then this is roughly equivalent to have a very large $\sigma_{q_i}$.  This analysis can be easily extended to the case when $\Tb_\eta$ is not diagonalizable by the same bases as $\Tb_q$ by some computations to calculate the variance of $\eta_i$.

\newcommand{\barsigma}{\delta}
\newcommand{\qp}{q_p}
\newcommand{\qb}{q_b}

\section{Numerical experiments}\label{s:results}
We conduct several numerical experiments to compare the methods introduced in Sections~\ref{s:noprior} and \ref{s:prior}.
In all our results, we use the spectral representation (sine series) for both $q$ and $\eta$  but we  consider different noise levels,  regularization, and observation scenarios.

Below, we summarize the experiments we conducted.
\begin{itemize}
\item In Example \ref{example:SDNP}, we compare {\bf SISDNP} and {\bf MISDNP} algorithms. We use spectral truncation for $q$ so that the problem is well posed. We're just sampling $\eta$, not inverting for it.  The $\eta$ prior used for sampling is a Mat\'{e}rn-like isotropic prior. As mentioned before the data was produced by a single realization of $\eta$.  
\item In Example \ref{example:MDNP}, we compare the {\bf MIMDNP} and {\bf SIMDNP} algorithms. Here we assume we have data (in all receiver locations for all illumination angles) from multiple realizations of $\eta$. The priors for $\eta$ and $q$ are the same as in \ref{example:SDNP}. 
\item In Example \ref{example:SDP}, we use the {\bf SISDP} algorithm to solve the inverse problem in three different cases for several noise levels. The data is observed for a single realization of $\eta$.
\begin{enumerate}[label=(\alph*)]
\item $q$ and $\eta$ have exactly the same prior; 
\item $q$ and $\eta$ have very different priors: $q$ is supported only in low-frequencies and $\eta$ only in high frequencies (high-frequency noise). Since the inverse problem is not stable in high frequencies we need regularization for $\eta$. 
\item $q$ and $\eta$ have partially overlapping priors. 
\end{enumerate}
\item Example \ref{example:MDP}, we apply the {\bf MIMDP}  algorithm in two different cases: 
\begin{enumerate}[label=(\alph*)]
\item $q$ and $\eta$ have similar power spectrum and the same prior; and
\item $q$ and $\eta$ have similar power spectrum and different priors.
\end{enumerate}
In this example, the data is observed for multiple realizations of $\eta$. 
\end{itemize}

In Table \ref{table:examples}, we summarize the examples we investigated with their respective goals and figures with its results. Table \ref{table:figures} summarizes the contents of the figures with their descriptions.
\begin{table}[h!]
\caption{List of Examples, their goals and respective figures.}\label{table:examples}
{\small
\begin{center}

\bgroup
\def\arraystretch{1.2}
\begin{tabular}{| c | c | c |}
\hline
Example & Goal & Figures \\
\hline\hline
\ref{example:SDNP} & Compare {\bf SISDNP} to {\bf MISDNP} & \ref{ex:SDNP_results} \\
\hline
\ref{example:MDNP} & Compare {\bf SIMDNP}  to {\bf MIMDNP} & \ref{error_exp1_MIMDNP}, \ref{error_exp1_SIMDNP}, \ref{ex1_MIMDNP_results}, \ref{ex1_SIMDNP_results}  \\
\hline
\ref{example:SDP} part (a) & \begin{tabular}{@{}c@{}}Show limitations of the {\bf SISDP} \\ and compare to {\bf RLA} \end{tabular} & \ref{ex:SISDP:fig1}, \ref{ex:SISDP:qb:af:results}, \ref{ex:SISDP:qp:af:results}\\
\hline
\ref{example:SDP} part (b) & \begin{tabular}{@{}c@{}} Show advantages of {\bf SISDP} \\ for isotropic background medium \end{tabular} & \ref{fig_boxplots_bump_isotropic_hf_pld1}, \ref{fig_boxplots_bomber_isotropic_hf_pld1},\ref{fig_bumps_isotropic_hf_pld1},\ref{fig_plane_isotropic_hf_pld1} \\
\hline
\ref{example:SDP} part (c) & \begin{tabular}{@{}c@{}} Show advantages and limitations of {\bf SISDP} \\ for anisotropic background medium \end{tabular}& \ref{fig_anisotropic_af}  \\
\hline
\ref{example:MDP} part (a)  & \begin{tabular}{@{}c@{}}Show improvements of {\bf MIMDP} over \\ {\bf SISDP} for isotropic background medium \end{tabular} & \ref{fig_pid_plane_samples} \\
\hline
\ref{example:MDP} part (b)  & \begin{tabular}{@{}c@{}} Show improvements of {\bf MIMDP} over \\ {\bf SISDP} for anisotropic background medium \end{tabular} & \ref{fig_pid_submarine_samples}  \\
\hline 
\end{tabular}
\egroup

\end{center}
}
\end{table}

\begin{table}[h!]
\caption{List of Figures and their short descriptions.}\label{table:figures}
{\small
\begin{center}

\bgroup
\def\arraystretch{1.2}
\begin{tabular}{| c | l |}
\hline
Figure & Short description \\
\hline\hline
\ref{fig:domains} & Domains $\qb$, $\qp$ and $q_\mathrm{sub}$ used in our experiments.  \\
\hline
\ref{ex:SDNP_results} & \begin{tabular}{@{}l@{}} Comparison of the reconstructions of $\qb$ using {\bf SISDNP}  to {\bf MISDNP}   with \\ background medium with different noise-levels. \end{tabular} \\
\hline
\ref{error_exp1_MIMDNP} & \begin{tabular}{@{}l@{}} Plot of the error $E_\mathrm{MIMDNP}$ using different number of samples for background \\ medium with different noise levels. \end{tabular} \\
\hline
\ref{error_exp1_SIMDNP} & \begin{tabular}{@{}l@{}} Plot of the error $E_\mathrm{SIMDNP}$ using different number of samples for background \\ medium with different noise levels. \end{tabular} \\
\hline
\ref{ex1_MIMDNP_results} & \begin{tabular}{@{}l@{}} Approximation of $\qb$ obtained by {\bf MIMDNP} using different number of samples  \\ for background medium with different noise levels. \end{tabular} \\
\hline 
\ref{ex1_SIMDNP_results}  & \begin{tabular}{@{}l@{}} Approximation of $\qb$ obtained by {\bf SIMDNP} using different number of samples  \\ for background medium with different noise levels. \end{tabular} \\
\hline
\ref{ex:SISDP:fig1} & \begin{tabular}{@{}l@{}} Approximations of $\qb$ and $\qp$ obtained using the {\bf RLA} for an isotropic background \\ medium, and $q$ and $\eta$ have similar prior. \end{tabular} \\ 
\hline
\ref{ex:SISDP:qb:af:results} & \begin{tabular}{@{}l@{}} Approximation of $\qb$ obtained by {\bf SISDP} using different values of the regularization \\ parameter for an isotropic background medium and $q$ and $\eta$ have similar prior.  \end{tabular} \\
\hline
\ref{ex:SISDP:qp:af:results} & \begin{tabular}{@{}l@{}} Approximation of $\qp$ obtained by {\bf SISDP} using different values of the regularization \\ parameter for an isotropic background medium, and $q$ and $\eta$ have similar prior.  \end{tabular} \\
\hline
\ref{fig_boxplots_bump_isotropic_hf_pld1} & \begin{tabular}{@{}l@{}} Box plots of the error $E_\mathrm{SISDP}$ in the reconstruction of $\qb$ \\ at different wavenumbers when the background medium has different noise levels. \end{tabular} \\
\hline
\ref{fig_boxplots_bomber_isotropic_hf_pld1} & \begin{tabular}{@{}l@{}} Box plots of the error $E_\mathrm{SISDP}$ in the reconstruction of $\qp$ \\ at different wavenumbers when the background medium has different noise levels. \end{tabular} \\
\hline
\ref{fig_bumps_isotropic_hf_pld1} & \begin{tabular}{@{}l@{}} Approximation of $\qb$ obtained by {\bf SISDP} for an isotropic background medium with\\ different noise levels when $\qb$ and $\eta$ have different prior in the frequency domain. \end{tabular} \\
\hline
\ref{fig_plane_isotropic_hf_pld1} & \begin{tabular}{@{}l@{}} Approximation of $\qp$ obtained by {\bf SISDP} for an isotropic background medium with\\ different noise levels when $\qp$ and $\eta$ have different prior in the frequency domain.\end{tabular} \\
\hline
\ref{fig_anisotropic_af} & \begin{tabular}{@{}l@{}} Approximation of $q_\mathrm{sub}$ obtained by {\bf SISDP} for an anisotropic background medium with\\ different noise levels, when $q_\mathrm{sub}$ and $\eta$ are generated by different priors. \end{tabular} \\
\hline
\ref{fig_pid_plane_samples} & \begin{tabular}{@{}l@{}} Approximation of $\qp$ obtained by {\bf MIMDP} using different number of samples for an \\ isotropic background medium when $\qp$ and $\eta$ have similar prior. \end{tabular} \\
\hline
\ref{fig_pid_submarine_samples}  & \begin{tabular}{@{}l@{}} Approximation of $q_\mathrm{sub}$ obtained by {\bf MIMDP} using different number of samples for an \\ anisotropic background medium when $q_\mathrm{sub}$ and $\eta$ are generated by different priors. \end{tabular} \\
\hline 
\end{tabular}
\egroup

\end{center}
}
\end{table}

\begin{figure}[h!]
  \centering
\subfigure[Isometric view of $-\qb$]{
\includegraphics[width=0.3\textwidth]{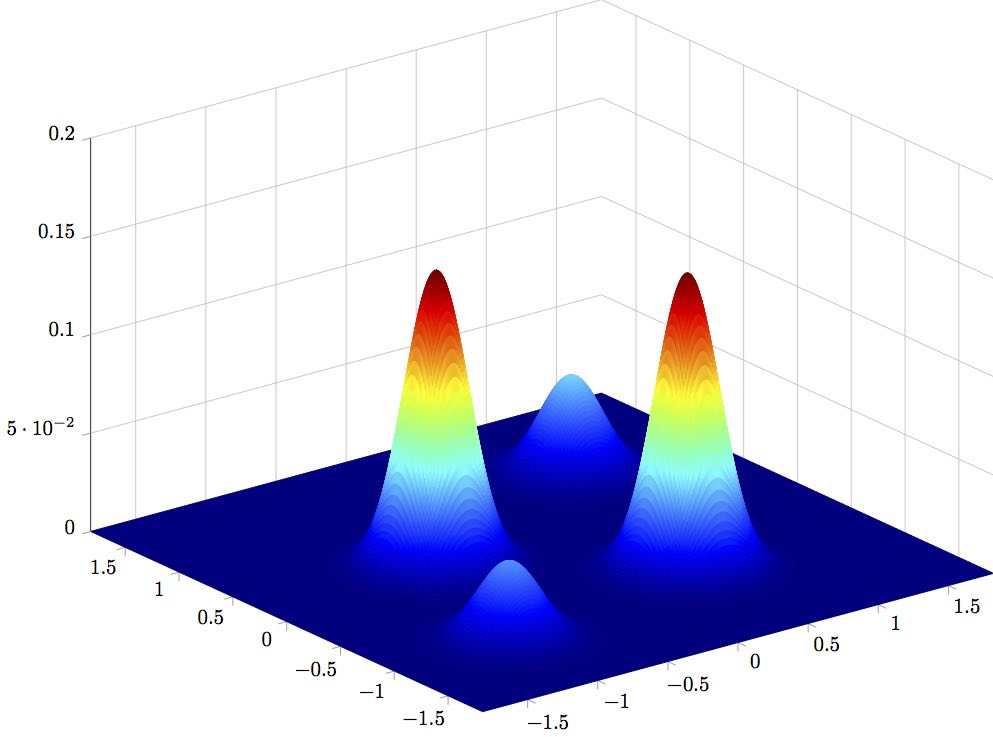} \label{contrast_qbumps_iso}
}
\subfigure[Isometric view of $-\qp$]{
\includegraphics[width=0.3\textwidth]{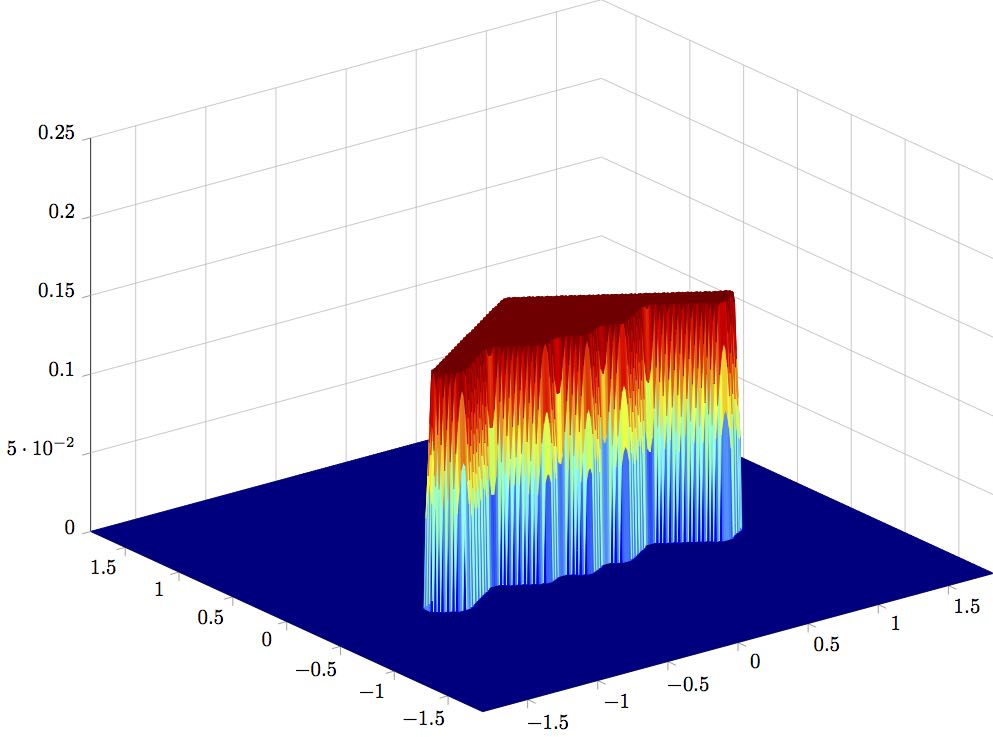} \label{contrast_qplane_iso}
}
\subfigure[Isometric view of $-q_\mathrm{sub}$]{
\includegraphics[width=0.3\textwidth]{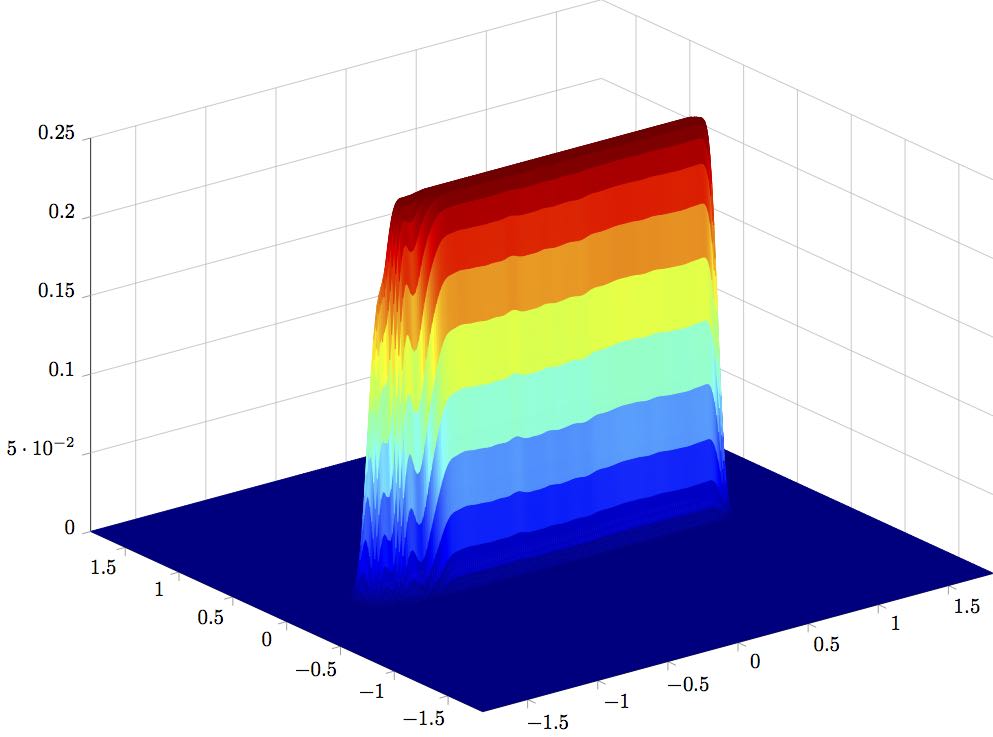} \label{contrast_qsub_iso}
}

\subfigure[Top view of $\qb$]{
\includegraphics[width=0.3\textwidth]{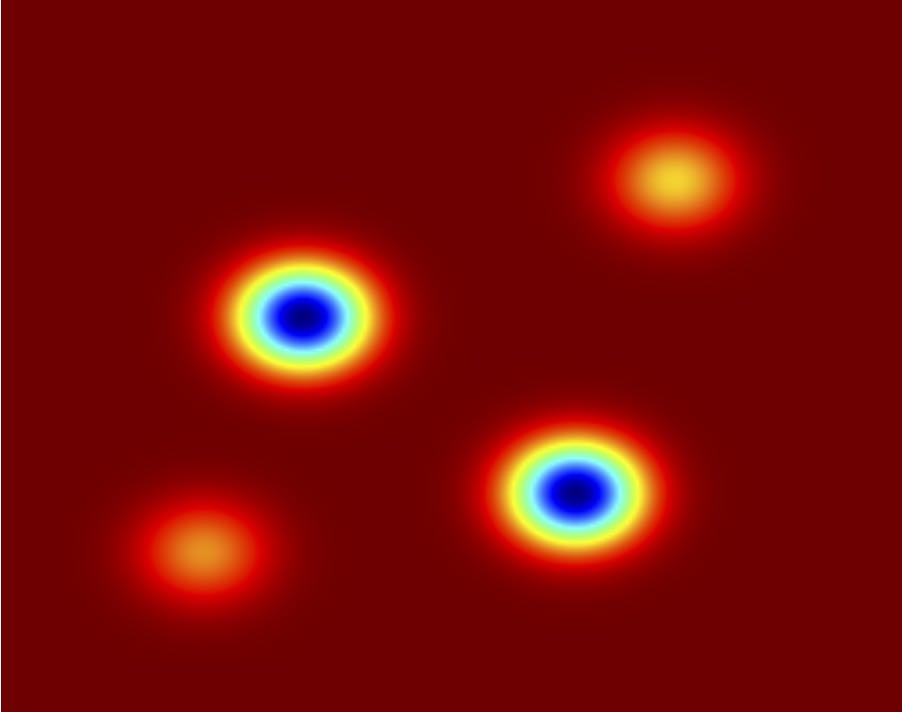} \label{contrast_qbumps_top}
}
\subfigure[Top view of $\qp$]{
\includegraphics[width=0.3\textwidth]{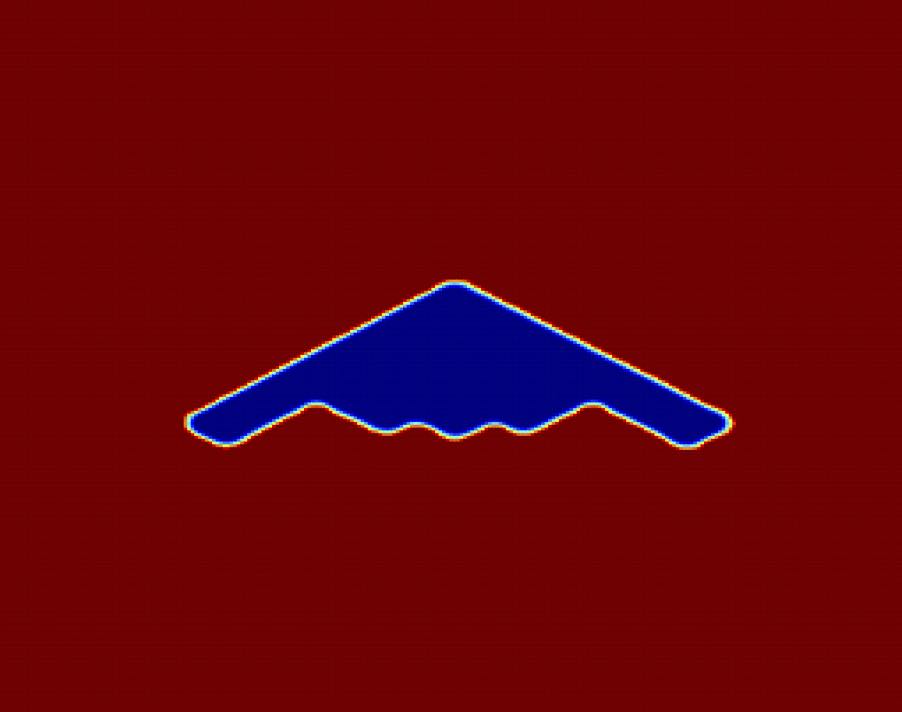} \label{contrast_qplane_top}
}
\subfigure[Top view of $q_\mathrm{sub}$]{
\includegraphics[width=0.3\textwidth]{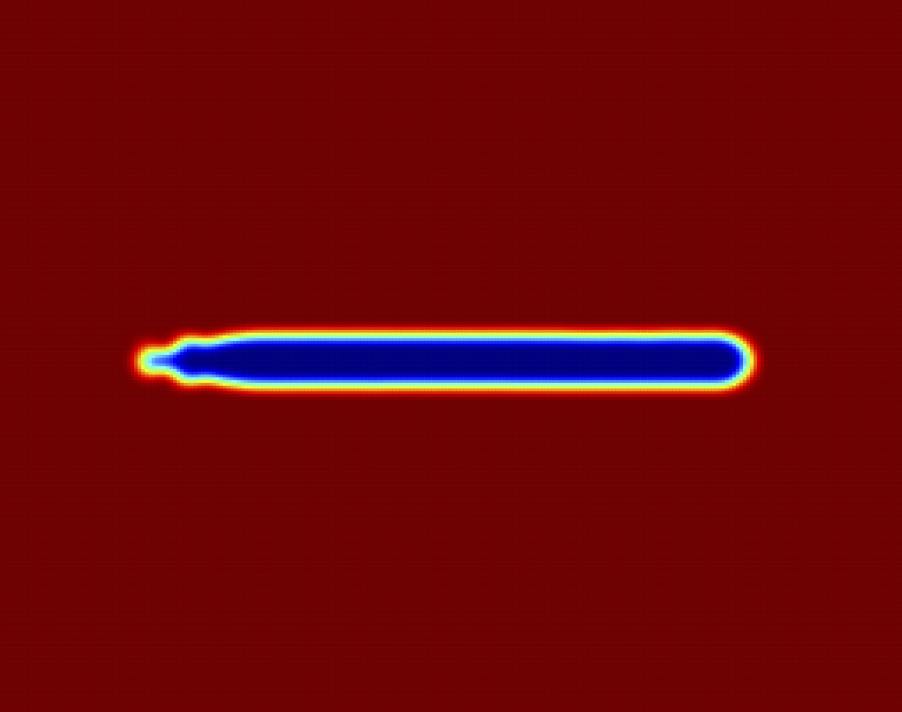} \label{contrast_qsub_top}
}
\caption{\footnotesize Plots of the scatterers used in our examples. }\label{fig:domains}
\end{figure}

{\bf Scatterer functions:} Three different scatterers are used in our experiments. The first scatterer considered can be seen in Figures \ref{contrast_qbumps_iso} and \ref{contrast_qbumps_top} and is composed of 4 Gaussian bumps with compact support on the domain $\Omega=\left[-\pi/2,\pi/2\right]^2$. Its analytical representation is
\begin{equation*}
\qb(\xb)=\sum_{j=1}^{4}g_j(\xb),
\end{equation*}
where $g_1(\xb)=-0.15\exp(-15\|\xb-(-0.6,0.2)\|^2)$, $g_2(\xb)=-0.15\exp(-15\|\xb-(0.5,-0.7)\|^2)$, $g_3(\xb)=-0.05\exp(-15\|\xb-(0.9,0.9)\|^2)$, and $g_4(\xb)=-0.04\exp(-15\|\xb-(-1.0,-1.0)\|^2)$. This function can be satisfactorily reconstructed with low-frequency modes of the sine series.

The second function considered is the scatterer $\qp$ that resembles the shape of a plane. The surface plots of the isometric and top views of this function are, respectively, in Figures \ref{contrast_qplane_iso} and \ref{contrast_qplane_top}. Due to the steep derivative on its boundary, the reconstruction of this function requires more frequency modes than that of $\qb$.

Finally, for the experiment where there is prior knowledge of the probability distribution of the scatterer and the background medium and these probabilities are regulated by anisotropic medium priors, we use a scatterer $q_\mathrm{sub}$ that looks like a submarine. The surface plots of the isometric and top views of this function are, respectively, in Figures \ref{contrast_qsub_iso} and \ref{contrast_qsub_top}.

To improve the clarity of the details of the function, in the isometric view of all examples presented, the scatterer is multiplied by $-1$.

\subsection{Forward and inverse solver configuration}\label{s:setup}
The HPS solver is used to generate data for the forward problem with at least 10 points per wavelength, which gives 5 digits of accuracy. Synthetic data measurements are generated for frequencies $k_j=k_1+(j-1)\delta k$, with $j=1,\ldots,Q$. At each frequency, the data is generated at $N_p=120$ points uniformly distributed in the circle of radius $R=3$ for $N_\theta=30$ incident plane waves with incidence direction $\theta_m=(\cos(2\pi m/N_\theta),\sin(2\pi m/N_\theta))$, $m=1,\ldots,N_\theta$. We do not add any noise to the data.  This is a classical ``inverse crime'' but the focus of the paper is not the particular inversion algorithm but the effect of the randomness of the background medium.

We represent the approximate solution of the domain using  \eqref{eq:sineseries}, with $M(k)=30$ for all $k$, which amounts to about 900 coefficients. For lower frequencies, the problem is still ill posed, which requires some regularization technique depending on the method that is used. We specify the regularization used in each example.

The main algorithmic component is a Gauss-Newton method for the solution of the inverse problem at a single frequency.  At each Gauss-Newton iteration $j$, the HPS solver is used to calculate the forward operator applied at the guess domain $q_j$ and its Fr\'{e}chet derivative. We use the HPS solver for the inverse problem with 8 points per wavelength, which gives 4 digits of accuracy, with the intent of introducing a small model error.
In our examples, we set $q_0$ to zero.

We use three stopping criteria for the Gauss-Newton iteration: 
\begin{enumerate}[label=(\alph*)]
\item we reach a  maximum number of iterations $N_\text{it}=20$;
\item the relative data-mismatch norm is below $\epsilon_\text{res}=10^{-7}$; or
\item the norm of the step update divided by the number of unknowns is smaller than $\epsilon_\text{step}=10^{-7}$.
\end{enumerate}

\subsection{Data from single $\eta$ and spectral truncation for $q$} \label{example:SDNP}

We present the results obtained using the {\bf SISDNP} and {\bf MISDNP} algorithms on one data set. The scatterer chosen to be recovered is $\qb$.
As mentioned we use a spectral truncation regularization. All high-frequency coefficients that satisfy $m_1+m_2>2k$ are filtered out. The data measurements are obtained at the minimum frequency $k_1=1$, with $\delta k=0.5$, and the number of frequencies used is $Q=18$, giving data measurements up to the maximum frequency $k_{18}=9.5$.

To generate the smooth-background noisy medium samples $\eta_s$, $s=1,\ldots,N_s$, we create a $N_\Omega\times N_\Omega$ uniform grid of points $\tilde{\bf x}$ in the domain $\Omega$. The background noisy medium is obtained by solving
\begin{equation}\label{e:ex1-eta-reg}
{\bf \eta}_s= (\mu I+D)^{-1} \sigma,
\end{equation}
where $\sigma\in\mathbb{R}^{N_\Omega^2}$ is a vector whose elements are obtained from the normal distribution $\mathcal{N}(0,\delta)$ with mean $0$ and variance $\delta$, $\mu=10^3$, $I$ is the $N_\Omega^2\times N_\Omega^2$ identity matrix, and $D$ is the $N_\Omega^2\times N_\Omega^2$ matrix obtained from discretizing the Neumann problem with the usual five-point operator on the $N_\Omega\times N_\Omega$ uniform grid.

\begin{figure}[!htp]
\centering
\subfigure[$\qb+\eta^\ast$ for $\delta=10$]{
\includegraphics[width=0.29\textwidth]{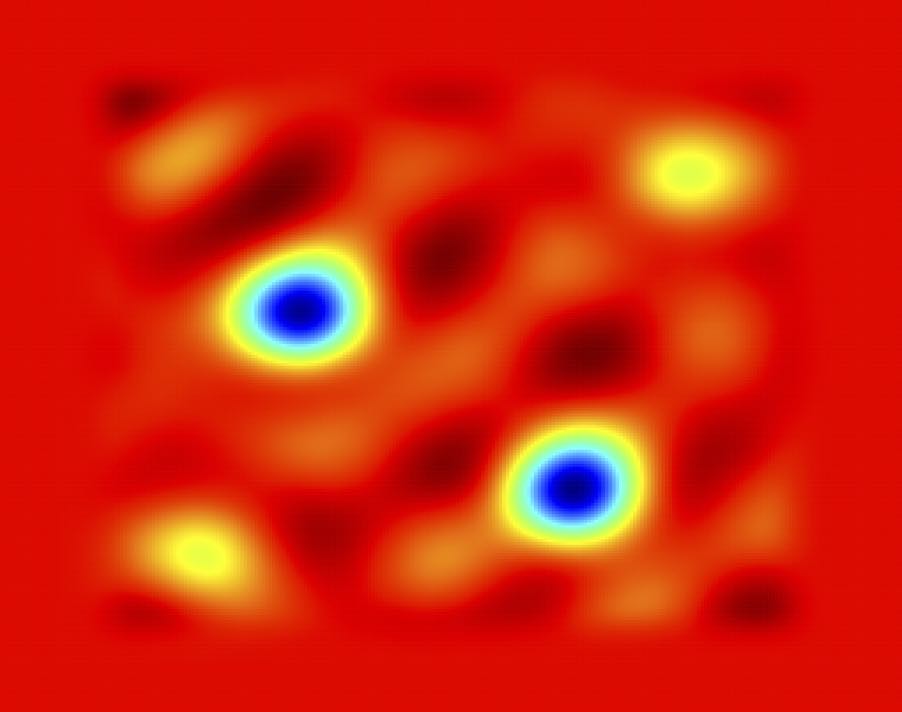}\label{ex1:domain_d10}
}
\subfigure[$q_\mathrm{SISDNP}$ for $\delta=10$]{
\includegraphics[width=0.29\textwidth]{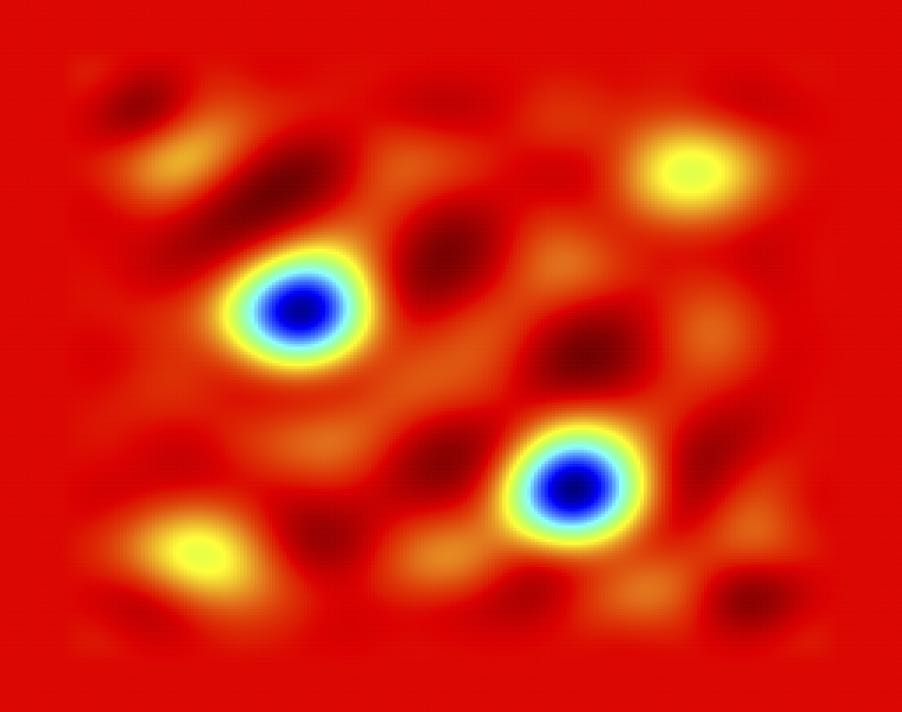}\label{ex1:SISDNP_d10}
}
\subfigure[$q_\mathrm{MISDNP}$ for $\delta=10$]{
\includegraphics[width=0.29\textwidth]{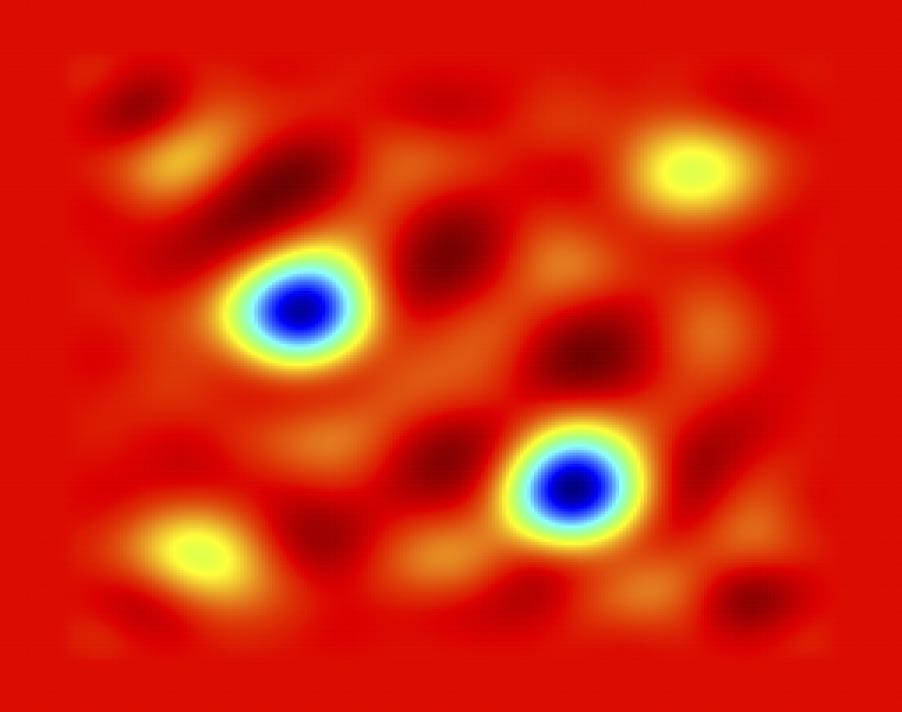}\label{ex1:=MISDNP_d10}
}

\subfigure[$\qb+\eta^\ast$ for $\delta=20$]{
\includegraphics[width=0.29\textwidth]{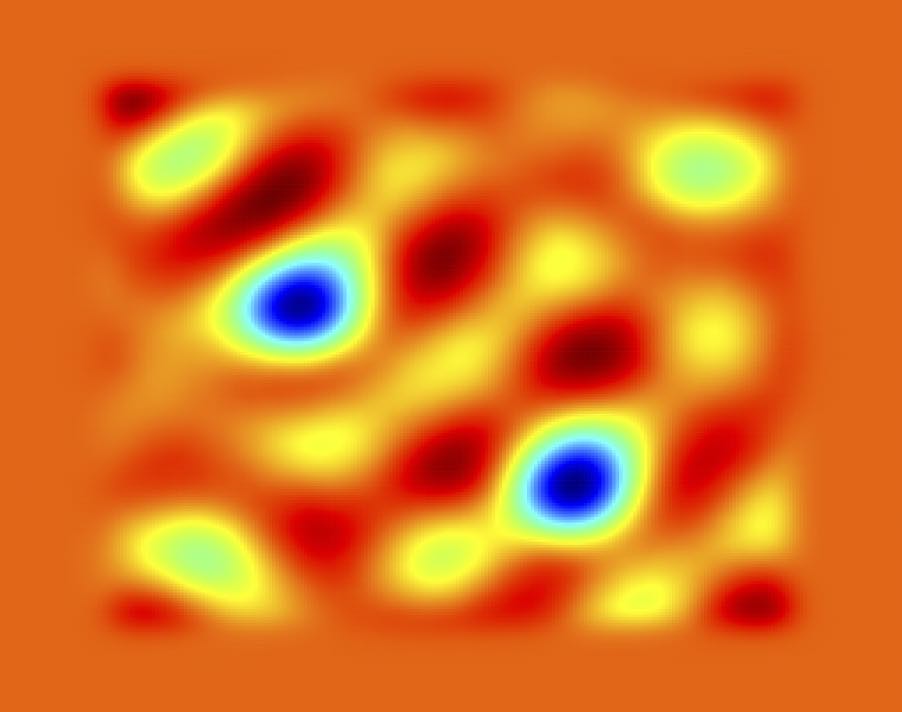}\label{ex1:domain_d20}
}
\subfigure[$q_\mathrm{SISDNP}$ for $\delta=20$]{
\includegraphics[width=0.29\textwidth]{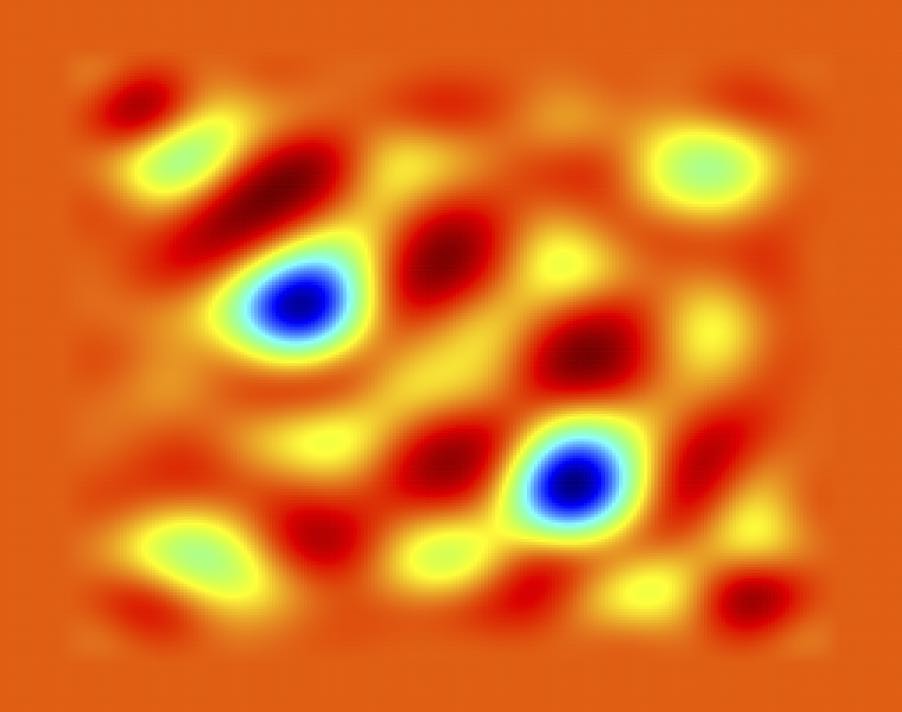}\label{ex1:SISDNP_d20}
}
\subfigure[$q_\mathrm{MISDNP}$ for $\delta=20$]{
\includegraphics[width=0.29\textwidth]{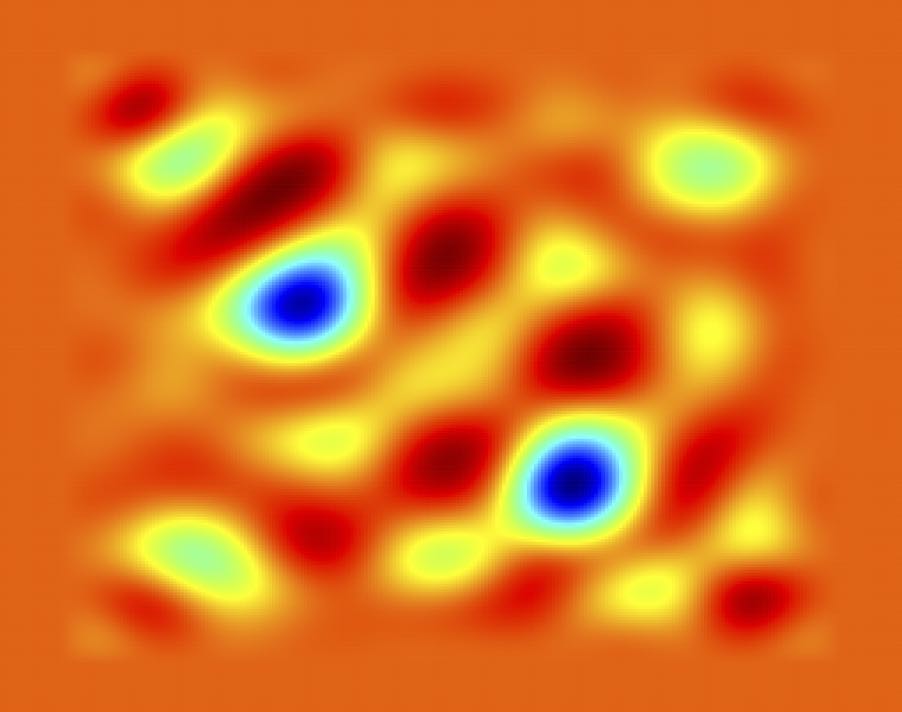}\label{ex1:MISDNP_d20}
}

\subfigure[$\qb+\eta^\ast$ for $\delta=40$]{
\includegraphics[width=0.29\textwidth]{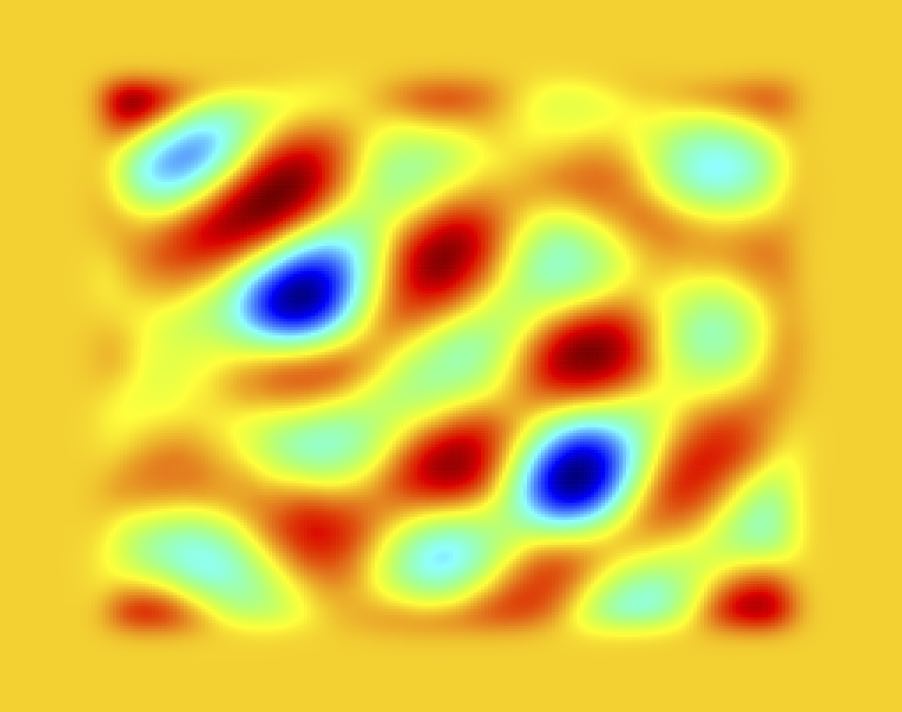}\label{ex1:domain_d40}
}
\subfigure[$q_\mathrm{SISDNP}$ for $\delta=40$]{
\includegraphics[width=0.29\textwidth]{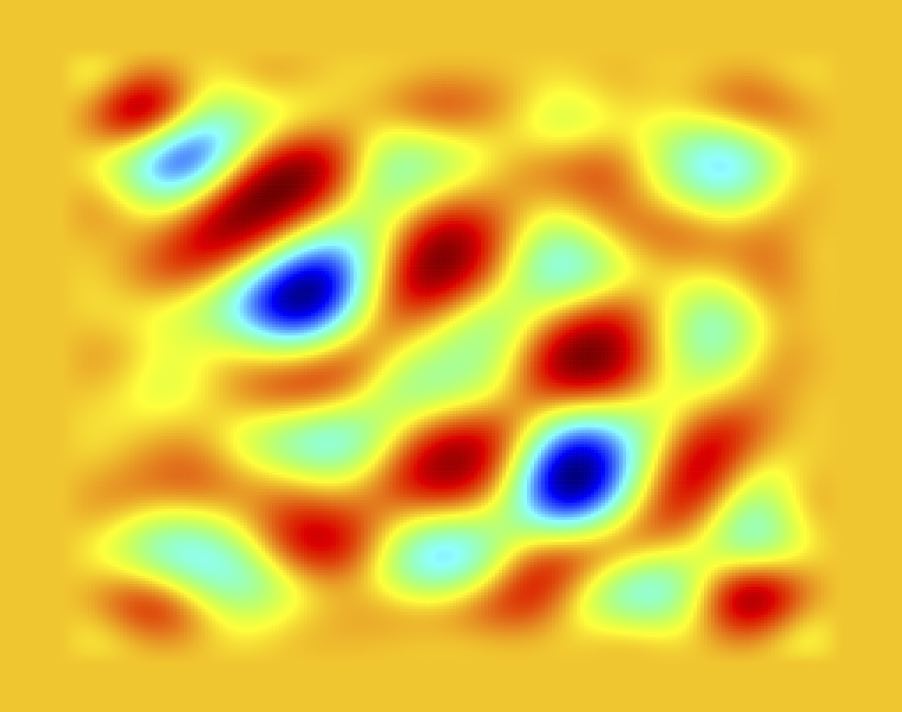}\label{ex1:SISDNP_d40}
}
\subfigure[$q_\mathrm{MISDNP}$ for $\delta=40$]{
\includegraphics[width=0.29\textwidth]{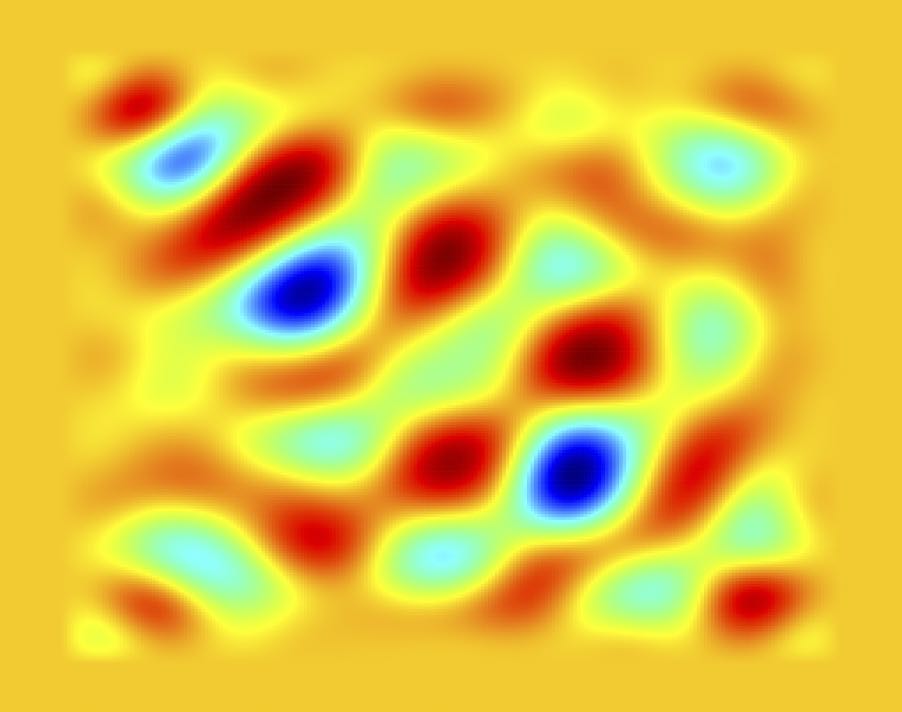}\label{ex1:MISDNP_d40}
}

\subfigure[$\qb+\eta^\ast$ for $\delta=80$]{
\includegraphics[width=0.29\textwidth]{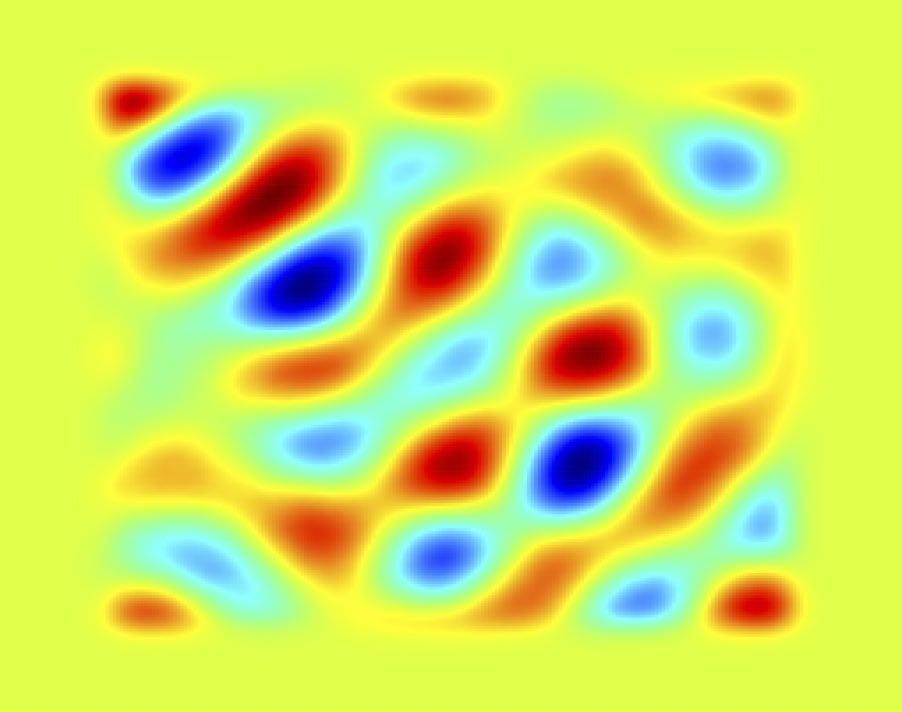}\label{ex1:domain_d80}
}
\subfigure[$q_\mathrm{SISDNP}$ for $\delta=80$]{
\includegraphics[width=0.29\textwidth]{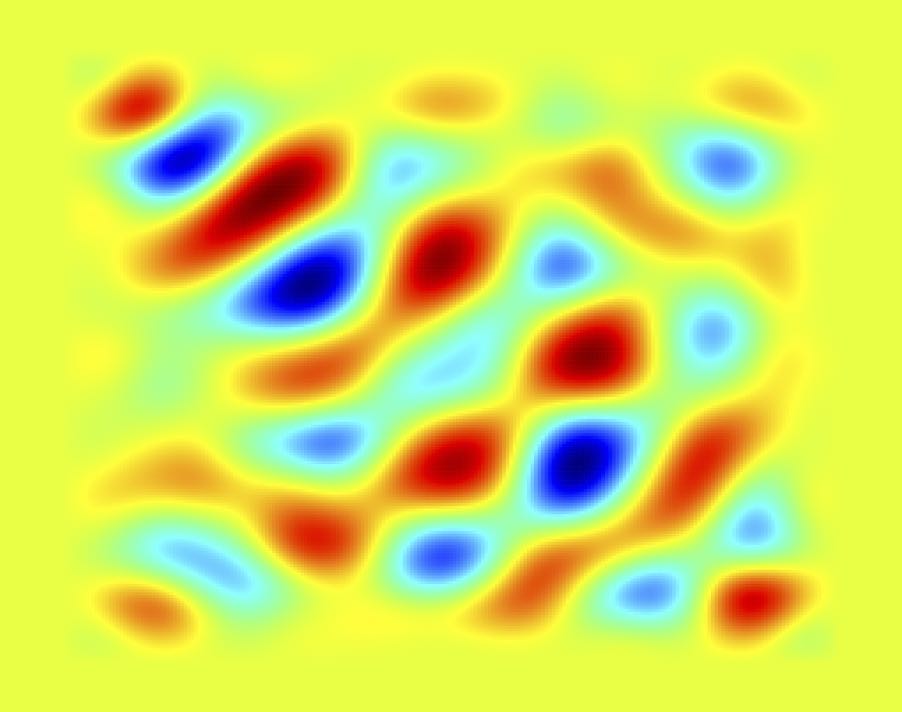}\label{ex1:SISDNP_d80}
}
\subfigure[$q_\mathrm{MISDNP}$ for $\delta=80$]{
\includegraphics[width=0.29\textwidth]{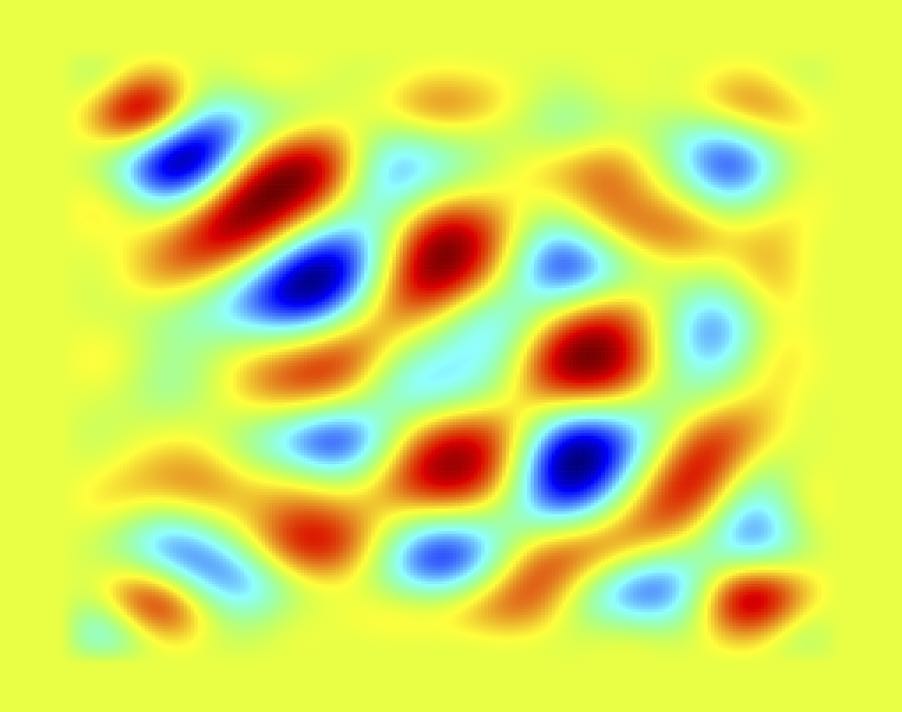}\label{ex1:MISDNP_d80}
}

\caption{Reconstruction of $\qb$ for Example \ref{example:SDNP}. The {\bf SISDNP} and {\bf MISDNP} algorithms are used to reconstruct the domain $\qb$ in the presence of a background medium. From top to bottom, we present the domain $\qb+\eta^\ast$, and the solutions $q_\mathrm{SIMDNP}$ and $q_\mathrm{MISDNP}$,  when the background medium is generated using the parameter $\delta=10$, $20$, $40$ and $80$. The solution $q_\mathrm{MISDNP}$ is always obtained using $N_s=100$ samples for the background domain.}\label{ex:SDNP_results}
\end{figure}

The data is generated for the domain $\qb+\eta^\ast$. In Figure \ref{ex:SDNP_results}, we present the reconstruction using the {\bf SISDNP} and {\bf MISDNP} algorithms for the background medium with $\delta=10$, $20$, $40$ and $80$. In this example, this translates to ratios of $\|\eta^\ast\|/\|q^\ast\|=0.36$, $0.72$, $1.43$, and $2.86$. For the solution with the {\bf MISDNP} algorithm, we use $N_s=100$ samples of $\eta$. Since $\eta$ has a significant number of components on the energy spectrum of $q$ it is very hard to separate them: our reconstructions have information from both $\qb$ and the background noisy medium $\eta^\ast$, and it is not possible to separate this information.  It is easy to show that if the inverse problem is linear the point estimate solutions to {\bf SISDNP} and {\bf MISDNP} are identical. Experimentally, for this example, we reach the same conclusion for specific nonlinear scattering problem. 

\subsection{Data from multiple $\eta$ and spectral truncation for $q$}\label{example:MDNP}

This example is an extension of the Example \ref{example:SDNP}. We use the same scatterer and generate $\eta$ using the same procedure; but here we assume that we have  scattered data for $\qb+\eta_s$, for $s=1,\ldots,N_s$. Of course the values of $\eta_s$ are not known but are independently sampled  from~\eqref{e:ex1-eta-reg}. 

As in the previous example, we use noise levels $\delta=10$, $20$, $40$ and $80$  to generate samples for $\eta$. We apply the {\bf MIMDNP} and {\bf SIMDNP} algorithms  to recover the scatterer with $N_s=10$, $50$, $100$, $500$ and $1000$ samples for a background noise function with $\delta=10$ and $20$, and with $N_s=10$, $50$, $100$, $500$, $1000$, $2000$ and $4000$ for a background noise function with $\delta=40$ and $80$. 

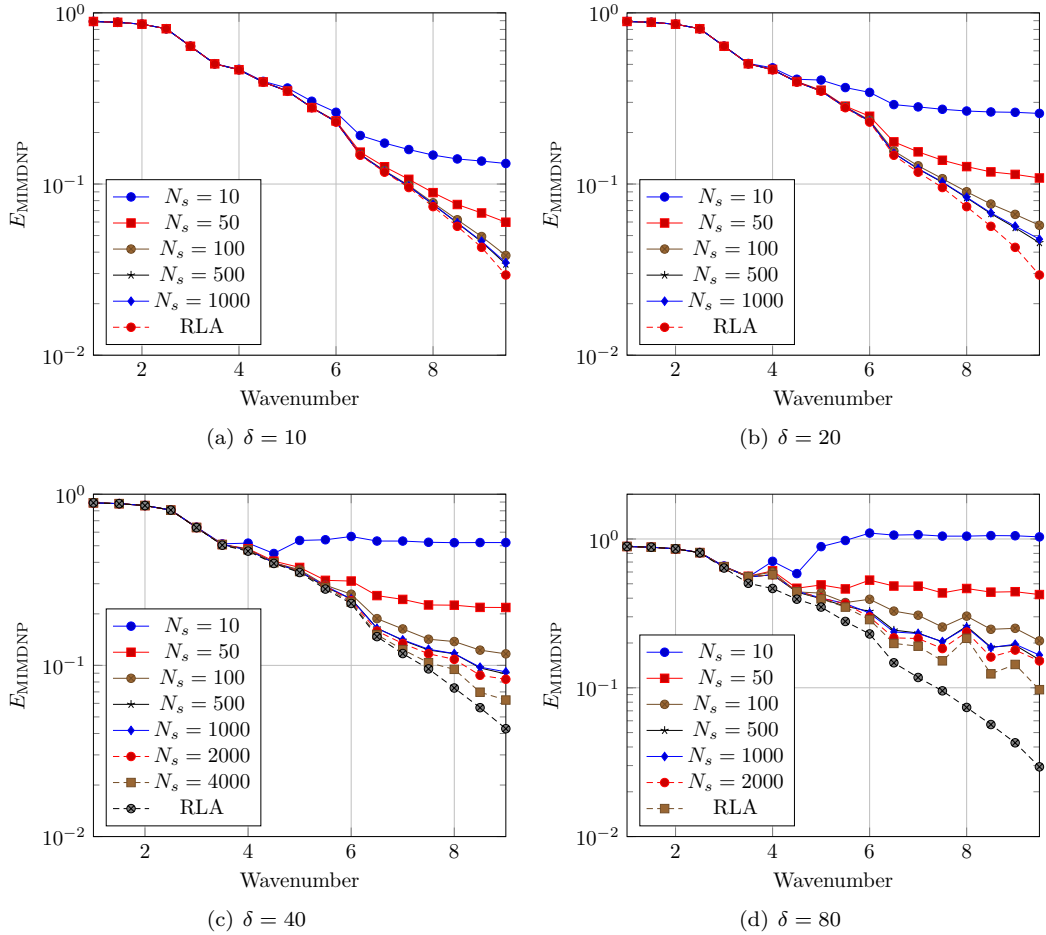
\begin{figure}[h!]
  \centering
\subfigure[$\delta=10$]{
\begin{tikzpicture}[scale=0.8]
\begin{semilogyaxis}[xmin=1, xmax=9.5,domain=1:9.5,ymin=1e-2, ymax=1,domain=1e-4:1,
    xlabel={Wavenumber},
    ylabel={$E_\mathrm{MIMDNP}$},
    grid=major,
    legend pos=south west,
    legend entries={$N_s=10$,$N_s=50$,$N_s=100$,$N_s=500$,$N_s=1000$, RLA},
]
\addplot coordinates {
(1,8.905749e-01) (1.5,8.807977e-01) (2,8.588176e-01) (2.5,8.081662e-01) (3,6.382204e-01) (3.5,5.043617e-01)
(4,4.679236e-01) (4.5,3.979904e-01) (5,3.639482e-01) (5.5,3.045967e-01) (6,2.625970e-01) (6.5,1.921515e-01)
(7,1.733987e-01) (7.5,1.590457e-01) (8,1.476157e-01) (8.5,1.401419e-01) (9,1.362228e-01) (9.5,1.318421e-01)
};
\addplot coordinates {
(1,8.905099e-01) (1.5,8.807231e-01) (2,8.586115e-01) (2.5,8.078946e-01) (3,6.379181e-01) (3.5,5.036669e-01)
(4,4.654713e-01) (4.5,3.949823e-01) (5,3.495916e-01) (5.5,2.795889e-01) (6,2.337395e-01) (6.5,1.535954e-01)
(7,1.262142e-01) (7.5,1.065686e-01) (8,8.906962e-02) (8.5,7.592751e-02) (9,6.775805e-02) (9.5,5.984426e-02)
};
\addplot coordinates {
(1,8.905120e-01) (1.5,8.807283e-01) (2,8.586087e-01) (2.5,8.078616e-01) (3,6.379882e-01) (3.5,5.036927e-01)
(4,4.650981e-01) (4.5,3.942259e-01) (5,3.495636e-01) (5.5,2.803177e-01) (6,2.313263e-01) (6.5,1.491596e-01)
(7,1.196586e-01) (7.5,9.827026e-02) (8,7.773157e-02) (8.5,6.184938e-02) (9,4.936052e-02) (9.5,3.825456e-02)
};
\addplot coordinates {
(1,8.905019e-01) (1.5,8.807194e-01) (2,8.586020e-01) (2.5,8.078364e-01) (3,6.379982e-01) (3.5,5.037122e-01)
(4,4.649259e-01) (4.5,3.940283e-01) (5,3.487649e-01) (5.5,2.794725e-01) (6,2.308917e-01) (6.5,1.483731e-01 )
(7,1.186726e-01) (7.5,9.716863e-02) (8,7.590270e-02) (8.5,5.919425e-02) (9,4.594865e-02) (9.5,3.401364e-02)
};
\addplot coordinates {
(1,8.905003e-01) (1.5,8.807192e-01) (2,8.586053e-01) (2.5,8.078218e-01) (3,6.380188e-01) (3.5,5.037078e-01)
(4,4.649193e-01) (4.5,3.940718e-01) (5,3.486128e-01) (5.5,2.795123e-01) (6,2.308761e-01) (6.5,1.482288e-01)
(7,1.187622e-01) (7.5,9.729271e-02) (8,7.615305e-02) (8.5,5.942253e-02) (9,4.629814e-02) (9.5,3.476329e-02)
};
\addplot coordinates {
(1,8.905047e-01) (1.5,8.807221e-01) (2,8.586277e-01) (2.5,8.078779e-01) (3,6.379677e-01) (3.5,5.036224e-01)
(4,4.648982e-01) (4.5,3.940108e-01) (5,3.487298e-01) (5.5,2.791646e-01) (6,2.300323e-01) (6.5,1.472535e-01)
(7,1.172597e-01) (7.5,9.538726e-02) (8,7.377173e-02) (8.5,5.662153e-02) (9,4.267085e-02) (9.5,2.941934e-02)
};
\end{semilogyaxis}
\end{tikzpicture}
}
\subfigure[$\delta=20$]{
\begin{tikzpicture}[scale=0.8]
\begin{semilogyaxis}[xmin=1, xmax=9.5,domain=1:9.5,ymin=1e-2, ymax=1,domain=1e-4:1,
    xlabel={Wavenumber},
    ylabel={$E_\mathrm{MIMDNP}$},
    grid=major,
    legend pos=south west,
    legend entries={$N_s=10$,$N_s=50$,$N_s=100$,$N_s=500$,$N_s=1000$, RLA},
]
\addplot coordinates {
(1,8.906773e-01) (1.5,8.809039e-01) (2,8.590748e-01) (2.5,8.084981e-01) (3,6.391197e-01) (3.5,5.060909e-01)
(4,4.774077e-01) (4.5,4.094603e-01) (5,4.049270e-01) (5.5,3.660846e-01) (6,3.431286e-01) (6.5,2.910742e-01)
(7,2.822404e-01) (7.5,2.733601e-01) (8,2.670062e-01) (8.5,2.636856e-01) (9,2.624232e-01) (9.5,2.585492e-01)
};
\addplot coordinates {
(1,8.905122e-01) (1.5,8.807179e-01) (2,8.585625e-01) (2.5,8.078492e-01) (3,6.379846e-01) (3.5,5.040756e-01)
(4,4.672417e-01) (4.5,3.971513e-01) (5,3.534562e-01) (5.5,2.851481e-01) (6,2.488874e-01) (6.5,1.763224e-01)
(7,1.541781e-01) (7.5,1.377475e-01) (8,1.264202e-01) (8.5,1.177334e-01) (9,1.138507e-01) (9.5,1.084674e-01)
};
\addplot coordinates {
(1,8.905165e-01) (1.5,8.807278e-01) (2,8.585540e-01) (2.5,8.077752e-01) (3,6.380891e-01) (3.5,5.040628e-01)
(4,4.660601e-01) (4.5,3.947309e-01) (5,3.515575e-01) (5.5,2.827502e-01) (6,2.360314e-01) (6.5,1.559455e-01)
(7,1.278030e-01) (7.5,1.075077e-01) (8,9.016266e-02) (8.5,7.641191e-02) (9,6.647632e-02) (9.5,5.741441e-02)
};
\addplot coordinates {
(1,8.904967e-01) (1.5,8.807125e-01) (2,8.585458e-01) (2.5,8.077398e-01) (3,6.381082e-01) (3.5,5.040920e-01)
(4,4.654432e-01) (4.5,3.942801e-01) (5,3.491303e-01) (5.5,2.803209e-01) (6,2.331489e-01) (6.5,1.514003e-01)
(7,1.225504e-01) (7.5,1.022590e-01) (8,8.311162e-02) (8.5,6.712796e-02) (9,5.555798e-02) (9.5,4.541681e-02)
};
\addplot coordinates {
(1,8.904937e-01) (1.5,8.807130e-01) (2,8.585552e-01) (2.5,8.077163e-01) (3,6.381620e-01) (3.5,5.041056e-01)
(4,4.654770e-01) (4.5,3.945464e-01) (5,3.491709e-01) (5.5,2.811721e-01) (6,2.331377e-01) (6.5,1.512377e-01)
(7,1.229341e-01) (7.5,1.028065e-01) (8,8.373992e-02) (8.5,6.784373e-02) (9,5.675609e-02) (9.5,4.763799e-02)
};
\addplot coordinates {
(1,8.905047e-01) (1.5,8.807221e-01) (2,8.586277e-01) (2.5,8.078779e-01) (3,6.379677e-01) (3.5,5.036224e-01)
(4,4.648982e-01) (4.5,3.940108e-01) (5,3.487298e-01) (5.5,2.791646e-01) (6,2.300323e-01) (6.5,1.472535e-01)
(7,1.172597e-01) (7.5,9.538726e-02) (8,7.377173e-02) (8.5,5.662153e-02) (9,4.267085e-02) (9.5,2.941934e-02)
};
\end{semilogyaxis}
\end{tikzpicture}
}
\subfigure[$\delta=40$]{
\begin{tikzpicture}[scale=0.8]
\begin{semilogyaxis}[xmin=1, xmax=9,domain=1:9.5,ymin=1e-2, ymax=1,domain=1e-4:1,
    xlabel={Wavenumber},
    ylabel={$E_\mathrm{MIMDNP}$},
    grid=major,
    legend pos=south west,
    legend entries={$N_s=10$, $N_s=50$, $N_s=100$, $N_s=500$, $N_s=1000$, $N_s=2000$, $N_s=4000$, RLA},
]
\addplot coordinates {
(1,8.909578e-01) (1.5,8.811614e-01) (2,8.596812e-01) (2.5,8.091385e-01) (3,6.426982e-01) (3.5,5.129399e-01)
(4,5.179803e-01) (4.5,4.505336e-01) (5,5.369553e-01) (5.5,5.418553e-01) (6,5.662233e-01) (6.5,5.326168e-01)
(7,5.325150e-01) (7.5,5.239969e-01) (8,5.208914e-01) (8.5,5.219307e-01) (9,5.217784e-01) (9.5,5.146145e-01)
};
\addplot coordinates {
(1,8.905069e-01) (1.5,8.806970e-01) (2,8.583985e-01) (2.5,8.077197e-01) (3,6.389298e-01) (3.5,5.081813e-01)
(4,4.803347e-01) (4.5,4.065547e-01) (5,3.734449e-01) (5.5,3.139546e-01) (6,3.109679e-01) (6.5,2.556601e-01)
(7,2.429374e-01) (7.5,2.254108e-01) (8,2.247647e-01) (8.5,2.182167e-01) (9,2.178921e-01) (9.5,2.116410e-01)
};
\addplot coordinates {
(1,8.905148e-01) (1.5,8.807122e-01) (2,8.583650e-01) (2.5,8.075339e-01) (3,6.389949e-01) (3.5,5.078829e-01)
(4,4.763401e-01) (4.5,3.982338e-01) (5,3.617309e-01) (5.5,2.941064e-01) (6,2.605628e-01) (6.5,1.875965e-01)
(7,1.632893e-01) (7.5,1.422152e-01) (8,1.378209e-01) (8.5,1.227680e-01) (9,1.167975e-01) (9.5,1.050989e-01)
};
\addplot coordinates {
(1,8.904821e-01) (1.5,8.807025e-01) (2,8.583936e-01) (2.5,8.075539e-01) (3,6.390540e-01) (3.5,5.077668e-01)
(4,4.730648e-01) (4.5,3.974056e-01) (5,3.532489e-01) (5.5,2.870509e-01) (6,2.448139e-01) (6.5,1.651410e-01)
(7,1.401373e-01) (7.5,1.229702e-01) (8,1.169399e-01) (8.5,9.665591e-02) (9,8.911025e-02) (9.5,7.691987e-02)
};
\addplot coordinates {
(1,8.904784e-01) (1.5,8.807091e-01) (2,8.584292e-01) (2.5,8.075397e-01) (3,6.392294e-01) (3.5,5.079089e-01)
(4,4.731829e-01) (4.5,3.987930e-01) (5,3.545434e-01) (5.5,2.916522e-01) (6,2.440984e-01) (6.5,1.645096e-01)
(7,1.410138e-01) (7.5,1.243401e-01) (8,1.171708e-01) (8.5,9.765252e-02) (9,9.169340e-02) (9.5,8.249878e-02)
};
\addplot coordinates {
(1,8.904742e-01) (1.5,8.807064e-01) (2,8.584110e-01) (2.5,8.075002e-01) (3,6.393750e-01) (3.5,5.080668e-01)
(4,4.731002e-01) (4.5,3.992038e-01) (5,3.543568e-01) (5.5,2.908992e-01) (6,2.392087e-01) (6.5,1.587272e-01)
(7,1.342345e-01) (7.5,1.169429e-01) (8,1.084889e-01) (8.5,8.751819e-02) (9,8.298786e-02) (9.5,7.496193e-02)
};
\addplot coordinates {
(1,8.904778e-01) (1.5,8.806996e-01) (2,8.584170e-01) (2.5,8.075500e-01) (3,6.391740e-01) (3.5,5.079084e-01)
(4,4.729517e-01) (4.5,3.978749e-01) (5,3.533221e-01) (5.5,2.855639e-01) (6,2.343381e-01) (6.5,1.511063e-01)
(7,1.236153e-01) (7.5,1.037435e-01) (8,9.484486e-02) (8.5,6.971425e-02) (9,6.274486e-02) (9.5,4.796734e-02)
};
\addplot coordinates {
(1,8.905047e-01) (1.5,8.807221e-01) (2,8.586277e-01) (2.5,8.078779e-01) (3,6.379677e-01) (3.5,5.036224e-01)
(4,4.648982e-01) (4.5,3.940108e-01) (5,3.487298e-01) (5.5,2.791646e-01) (6,2.300323e-01) (6.5,1.472535e-01)
(7,1.172597e-01) (7.5,9.538726e-02) (8,7.377173e-02) (8.5,5.662153e-02) (9,4.267085e-02) (9.5,2.941934e-02)
};
\end{semilogyaxis}
\end{tikzpicture}
}
\subfigure[$\delta=80$]{
\begin{tikzpicture}[scale=0.8]
\begin{semilogyaxis}[xmin=1, xmax=9.5,domain=1:9.5,ymin=1e-2, ymax=2,domain=1e-2:2,
    xlabel={Wavenumber},
    ylabel={$E_\mathrm{MIMDNP}$},
    grid=major,
    legend pos=south west,
    legend entries={$N_s=10$, $N_s=50$, $N_s=100$, $N_s=500$, $N_s=1000$, $N_s=2000$, RLA},
]
\addplot coordinates {
(1,8.916926e-01) (1.5,8.816125e-01) (2,8.608002e-01) (2.5,8.102633e-01) (3,6.585897e-01) (3.5,5.550325e-01)
(4,7.078735e-01) (4.5,5.851114e-01) (5,8.883447e-01) (5.5,9.769652e-01) (6,1.095853e+00) (6.5,1.062046e+00)
(7,1.069462e+00) (7.5,1.045295e+00) (8,1.044977e+00) (8.5,1.054001e+00) (9,1.051313e+00) (9.5,1.033256e+00)
};
\addplot coordinates {
(1,8.904747e-01) (1.5,8.807968e-01) (2,8.584884e-01) (2.5,8.098878e-01) (3,6.515788e-01) (3.5,5.618583e-01)
(4,6.106395e-01) (4.5,4.666692e-01) (5,4.912064e-01) (5.5,4.617492e-01) (6,5.290007e-01) (6.5,4.827920e-01)
(7,4.818244e-01) (7.5,4.340656e-01) (8,4.647295e-01) (8.5,4.392483e-01) (9,4.423165e-01) (9.5,4.231531e-01)
};
\addplot coordinates {
(1,8.904808e-01) (1.5,8.807948e-01) (2,8.583433e-01) (2.5,8.093893e-01) (3,6.513286e-01) (3.5,5.608191e-01)
(4,5.996723e-01) (4.5,4.409074e-01) (5,4.327616e-01) (5.5,3.740899e-01) (6,3.934114e-01) (6.5,3.275300e-01)
(7,3.072626e-01) (7.5,2.559002e-01) (8,3.029875e-01) (8.5,2.471794e-01) (9,2.509153e-01) (9.5,2.066837e-01)
};
\addplot coordinates {
(1,8.904785e-01) (1.5,8.809209e-01) (2,8.586945e-01) (2.5,8.097887e-01) (3,6.510237e-01) (3.5,5.558457e-01)
(4,5.756197e-01) (4.5,4.395869e-01) (5,3.971223e-01) (5.5,3.536697e-01) (6,3.260295e-01) (6.5,2.440681e-01)
(7,2.334376e-01) (7.5,2.034622e-01) (8,2.591835e-01) (8.5,1.876888e-01) (9,1.931174e-01) (9.5,1.543220e-01)
};
\addplot coordinates {
(1,8.904864e-01) (1.5,8.809697e-01) (2,8.588766e-01) (2.5,8.099545e-01) (3,6.517252e-01) (3.5,5.564763e-01)
(4,5.743494e-01) (4.5,4.457857e-01) (5,4.011396e-01) (5.5,3.684142e-01) (6,3.177219e-01) (6.5,2.373053e-01)
(7,2.315406e-01) (7.5,2.044368e-01) (8,2.545525e-01) (8.5,1.862070e-01) (9,1.962571e-01) (9.5,1.659316e-01)
};
\addplot coordinates {
(1,8.905034e-01) (1.5,8.810168e-01) (2,8.589852e-01) (2.5,8.101491e-01) (3,6.525181e-01) (3.5,5.579237e-01)
(4,5.744586e-01) (4.5,4.484304e-01) (5,4.014538e-01) (5.5,3.661215e-01) (6,2.990064e-01) (6.5,2.164845e-01)
(7,2.143942e-01) (7.5,1.835015e-01) (8,2.368609e-01) (8.5,1.609808e-01) (9,1.791913e-01) (9.5,1.515395e-01)
};
\addplot coordinates {
(1,8.904886e-01) (1.5,8.809564e-01) (2,8.588995e-01) (2.5,8.100799e-01) (3,6.518608e-01) (3.5,5.573245e-01)
(4,5.750136e-01) (4.5,4.434723e-01) (5,3.960039e-01) (5.5,3.484917e-01) (6,2.873190e-01) (6.5,1.987219e-01)
(7,1.903492e-01) (7.5,1.517460e-01) (8,2.149875e-01) (8.5,1.240890e-01) (9,1.434617e-01) (9.5,9.710325e-02)
};
\addplot coordinates {
(1,8.905047e-01) (1.5,8.807221e-01) (2,8.586277e-01) (2.5,8.078779e-01) (3,6.379677e-01) (3.5,5.036224e-01)
(4,4.648982e-01) (4.5,3.940108e-01) (5,3.487298e-01) (5.5,2.791646e-01) (6,2.300323e-01) (6.5,1.472535e-01)
(7,1.172597e-01) (7.5,9.538726e-02) (8,7.377173e-02) (8.5,5.662153e-02) (9,4.267085e-02) (9.5,2.941934e-02)
};
\end{semilogyaxis}
\end{tikzpicture}
}
 \caption{Plot of the error for the {\bf MIMDNP} algorithm at each wavenumber for different levels of noise: (a) $\delta=10$, (b) $\delta=20$, (c) $\delta=40$ and (d) $\delta=80$. In each plot, each line represents the error $E_\mathrm{MIMDNP}$ using a different number of samples $N_s$. We also include the error for the {\bf RLA} (with $\eta=0$) in each plot as a matter of comparison and as a benchmark for the best possible approximation for our problem.}\label{error_exp1_MIMDNP}
\end{figure}

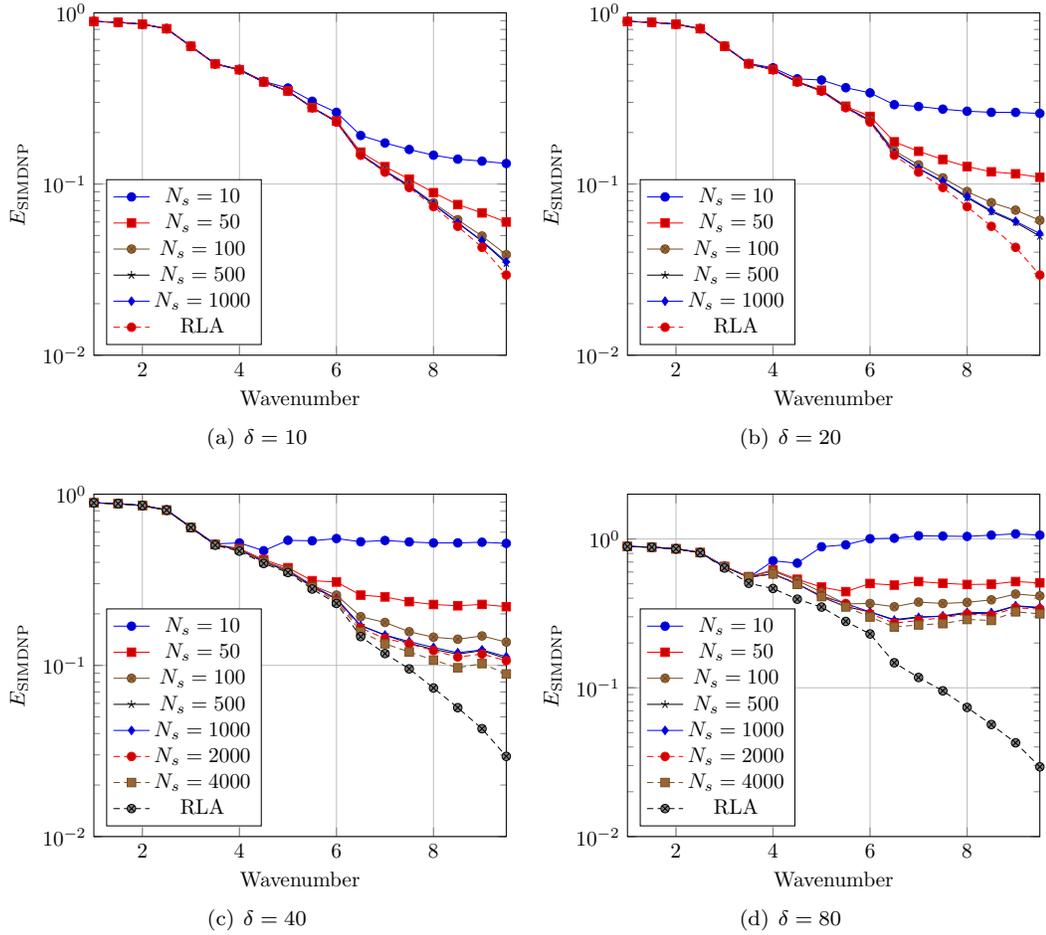
\begin{figure}[h!]
  \centering
\subfigure[$\delta=10$]{
\begin{tikzpicture}[scale=0.8]
\begin{semilogyaxis}[xmin=1, xmax=9.5,domain=1:9.5,ymin=1e-2, ymax=1,domain=1e-4:1,
    xlabel={Wavenumber},
    ylabel={$E_\mathrm{SIMDNP}$},
    grid=major,
    legend pos=south west,
    legend entries={$N_s=10$,$N_s=50$,$N_s=100$,$N_s=500$,$N_s=1000$, RLA},
]
\addplot coordinates {
(1,8.905749e-01) (1.5,8.807979e-01) (2,8.588179e-01) (2.5,8.081678e-01) (3,6.382186e-01) (3.5,5.043595e-01)
(4,4.679494e-01) (4.5,3.982710e-01) (5,3.639631e-01) (5.5,3.045385e-01) (6,2.623207e-01) (6.5,1.921674e-01)
(7,1.737544e-01) (7.5,1.591455e-01) (8,1.474193e-01) (8.5,1.397751e-01) (9,1.360705e-01) (9.5,1.317803e-01)
};
\addplot coordinates {
(1,8.905100e-01) (1.5,8.807232e-01) (2,8.586119e-01) (2.5,8.078965e-01) (3,6.379174e-01) (3.5,5.036665e-01)
(4,4.654810e-01) (4.5,3.950709e-01) (5,3.495845e-01) (5.5,2.795729e-01) (6,2.336794e-01) (6.5,1.536071e-01)
(7,1.263402e-01) (7.5,1.067356e-01) (8,8.904082e-02) (8.5,7.587476e-02) (9,6.781080e-02) (9.5,6.000195e-02)
};
\addplot coordinates {
(1,8.905121e-01) (1.5,8.807284e-01) (2,8.586091e-01) (2.5,8.078635e-01) (3,6.379878e-01) (3.5,5.036923e-01)
(4,4.648982e-01) (4.5,3.943168e-01) (5,3.495837e-01) (5.5,2.803020e-01) (6,2.312763e-01) (6.5,1.491995e-01)
(7,1.198075e-01) (7.5,9.834430e-02) (8,7.764994e-02) (8.5,6.186489e-02) (9,4.969587e-02) (9.5,3.866111e-02)
};
\addplot coordinates {
(1,8.905019e-01) (1.5,8.807195e-01) (2,8.586024e-01) (2.5,8.078380e-01) (3,6.379979e-01) (3.5,5.037123e-01)
(4,4.649252e-01) (4.5,3.940755e-01) (5,3.487850e-01) (5.5,2.794943e-01) (6,2.308812e-01) (6.5,1.484030e-01)
(7,1.186806e-01) (7.5,9.716200e-02) (8,7.584359e-02) (8.5,5.929455e-02) (9,4.631225e-02) (9.5,3.441522e-02)
};
\addplot coordinates {
(1,8.905003e-01) (1.5,8.807192e-01) (2,8.586057e-01) (2.5,8.078234e-01) (3,6.380187e-01) (3.5,5.037080e-01)
(4,4.649183e-01) (4.5,3.941038e-01) (5,3.486322e-01) (5.5,2.795419e-01) (6,2.308817e-01) (6.5,1.482568e-01)
(7,1.187532e-01) (7.5,9.734225e-02) (8,7.616870e-02) (8.5,5.956722e-02) (9,4.662430e-02) (9.5,3.511864e-02)
};
\addplot coordinates {
(1,8.905047e-01) (1.5,8.807221e-01) (2,8.586277e-01) (2.5,8.078779e-01) (3,6.379677e-01) (3.5,5.036224e-01)
(4,4.648982e-01) (4.5,3.940108e-01) (5,3.487298e-01) (5.5,2.791646e-01) (6,2.300323e-01) (6.5,1.472535e-01)
(7,1.172597e-01) (7.5,9.538726e-02) (8,7.377173e-02) (8.5,5.662153e-02) (9,4.267085e-02) (9.5,2.941934e-02)
};
\end{semilogyaxis}
\end{tikzpicture}
}
\subfigure[$\delta=20$]{
\begin{tikzpicture}[scale=0.8]
\begin{semilogyaxis}[xmin=1, xmax=9.5,domain=1:9.5,ymin=1e-2, ymax=1,domain=1e-4:1,
    xlabel={Wavenumber},
    ylabel={$E_\mathrm{SIMDNP}$},
    grid=major,
    legend pos=south west,
    legend entries={$N_s=10$,$N_s=50$,$N_s=100$,$N_s=500$,$N_s=1000$, RLA},
]
\addplot coordinates {
(1,8.906774e-01) (1.5,8.809049e-01) (2,8.590763e-01) (2.5,8.085045e-01) (3,6.391068e-01) (3.5,5.060646e-01)
(4,4.775988e-01) (4.5,4.117376e-01) (5,4.050556e-01) (5.5,3.653033e-01) (6,3.409222e-01) (6.5,2.907562e-01)
(7,2.837993e-01) (7.5,2.738676e-01) (8,2.663013e-01) (8.5,2.624333e-01) (9,2.621728e-01) (9.5,2.585277e-01)
};
\addplot coordinates {
(1,8.905123e-01) (1.5,8.807183e-01) (2,8.585640e-01) (2.5,8.078556e-01) (3,6.379805e-01) (3.5,5.040616e-01)
(4,4.673570e-01) (4.5,3.979329e-01) (5,3.533440e-01) (5.5,2.849711e-01) (6,2.483772e-01) (6.5,1.764739e-01)
(7,1.553525e-01) (7.5,1.391109e-01) (8,1.264933e-01) (8.5,1.178665e-01) (9,1.147707e-01) (9.5,1.095583e-01)
};
\addplot coordinates {
(1,8.905166e-01) (1.5,8.807283e-01) (2,8.585555e-01) (2.5,8.077817e-01) (3,6.380851e-01) (3.5,5.040496e-01)
(4,4.661322e-01) (4.5,3.956053e-01) (5,3.516653e-01) (5.5,2.826105e-01) (6,2.356249e-01) (6.5,1.564147e-01)
(7,1.294607e-01) (7.5,1.086661e-01) (8,9.039038e-02) (8.5,7.789120e-02) (9,7.037264e-02) (9.5,6.140755e-02)
};
\addplot coordinates {
(1,8.904968e-01) (1.5,8.807128e-01) (2,8.585472e-01) (2.5,8.077454e-01) (3,6.381063e-01) (3.5,5.040857e-01)
(4,4.654819e-01) (4.5,3.947385e-01) (5,3.492379e-01) (5.5,2.804119e-01) (6,2.330976e-01) (6.5,1.518517e-01)
(7,1.232117e-01) (7.5,1.028312e-01) (8,8.332247e-02) (8.5,6.892011e-02) (9,5.971604e-02) (9.5,4.955864e-02)
};
\addplot coordinates {
(1,8.904938e-01) (1.5,8.807133e-01) (2,8.585565e-01) (2.5,8.077217e-01) (3,6.381611e-01) (3.5,5.041006e-01)
(4,4.655118e-01) (4.5,3.948769e-01) (5,3.492792e-01) (5.5,2.813257e-01) (6,2.332125e-01) (6.5,1.516700e-01)
(7,1.234578e-01) (7.5,1.037662e-01) (8,8.449448e-02) (8.5,6.992418e-02) (9,6.061363e-02) (9.5,5.134449e-02)
};
\addplot coordinates {
(1,8.905047e-01) (1.5,8.807221e-01) (2,8.586277e-01) (2.5,8.078779e-01) (3,6.379677e-01) (3.5,5.036224e-01)
(4,4.648982e-01) (4.5,3.940108e-01) (5,3.487298e-01) (5.5,2.791646e-01) (6,2.300323e-01) (6.5,1.472535e-01)
(7,1.172597e-01) (7.5,9.538726e-02) (8,7.377173e-02) (8.5,5.662153e-02) (9,4.267085e-02) (9.5,2.941934e-02)
};
\end{semilogyaxis}
\end{tikzpicture}
}
\subfigure[$\delta=40$]{
\begin{tikzpicture}[scale=0.8]
\begin{semilogyaxis}[xmin=1, xmax=9.5,domain=1:9.5,ymin=1e-2, ymax=1,domain=1e-4:1,
    xlabel={Wavenumber},
    ylabel={$E_\mathrm{SIMDNP}$},
    grid=major,
    legend pos=south west,
    legend entries={$N_s=10$, $N_s=50$, $N_s=100$, $N_s=500$, $N_s=1000$, $N_s=2000$, $N_s=4000$, RLA},
]
\addplot coordinates {
(1,8.909588e-01) (1.5,8.811670e-01) (2,8.596884e-01) (2.5,8.091582e-01) (3,6.425944e-01) (3.5,5.126259e-01)
(4,5.193277e-01) (4.5,4.676624e-01) (5,5.379845e-01) (5.5,5.343362e-01) (6,5.508966e-01) (6.5,5.275316e-01)
(7,5.365419e-01) (7.5,5.262875e-01) (8,5.191113e-01) (8.5,5.195854e-01) (9,5.240684e-01) (9.5,5.160079e-01)
};
\addplot coordinates {
(1,8.905072e-01) (1.5,8.806979e-01) (2,8.584023e-01) (2.5,8.077307e-01) (3,6.388956e-01) (3.5,5.079511e-01)
(4,4.816027e-01) (4.5,4.144171e-01) (5,3.720738e-01) (5.5,3.124714e-01) (6,3.069449e-01) (6.5,2.573295e-01)
(7,2.511546e-01) (7.5,2.356525e-01) (8,2.268237e-01) (8.5,2.229515e-01) (9,2.271691e-01) (9.5,2.201972e-01)
};
\addplot coordinates {
(1,8.905151e-01) (1.5,8.807134e-01) (2,8.583689e-01) (2.5,8.075453e-01) (3,6.389578e-01) (3.5,5.076649e-01)
(4,4.774009e-01) (4.5,4.076204e-01) (5,3.625674e-01) (5.5,2.931648e-01) (6,2.574007e-01) (6.5,1.925478e-01)
(7,1.781524e-01) (7.5,1.573047e-01) (8,1.461484e-01) (8.5,1.422770e-01) (9,1.485553e-01) (9.5,1.368453e-01)
};
\addplot coordinates {
(1,8.904823e-01) (1.5,8.807025e-01) (2,8.583967e-01) (2.5,8.075636e-01) (3,6.390370e-01) (3.5,5.076309e-01)
(4,4.738399e-01) (4.5,4.030964e-01) (5,3.541500e-01) (5.5,2.874887e-01) (6,2.446907e-01) (6.5,1.710609e-01)
(7,1.504852e-01) (7.5,1.346493e-01) (8,1.240765e-01) (8.5,1.164198e-01) (9,1.220035e-01) (9.5,1.095352e-01)
};
\addplot coordinates {
(1,8.904786e-01) (1.5,8.807089e-01) (2,8.584322e-01) (2.5,8.075493e-01) (3,6.392202e-01) (3.5,5.077880e-01)
(4,4.739246e-01) (4.5,4.033941e-01) (5,3.555226e-01) (5.5,2.925248e-01) (6,2.449931e-01) (6.5,1.703725e-01)
(7,1.503820e-01) (7.5,1.378193e-01) (8,1.273044e-01) (8.5,1.186763e-01) (9,1.231745e-01) (9.5,1.124500e-01)
};
\addplot coordinates {
(1,8.904743e-01) (1.5,8.807060e-01) (2,8.584139e-01) (2.5,8.075093e-01) (3,6.393709e-01) (3.5,5.079588e-01)
(4,4.738073e-01) (4.5,4.030163e-01) (5,3.554854e-01) (5.5,2.923042e-01) (6,2.414751e-01) (6.5,1.652396e-01)
(7,1.436341e-01) (7.5,1.337246e-01) (8,1.221245e-01) (8.5,1.118690e-01) (9,1.164364e-01) (9.5,1.066063e-01)
};
\addplot coordinates {
(1,8.904779e-01) (1.5,8.806992e-01) (2,8.584210e-01) (2.5,8.075604e-01) (3,6.391679e-01) (3.5,5.077982e-01)
(4,4.736977e-01) (4.5,4.023148e-01) (5,3.543605e-01) (5.5,2.868414e-01) (6,2.360155e-01) (6.5,1.577113e-01)
(7,1.337810e-01) (7.5,1.198480e-01) (8,1.074263e-01) (8.5,9.685450e-02) (9,1.024575e-01) (9.5,8.913123e-02)
};
\addplot coordinates {
(1,8.905047e-01) (1.5,8.807221e-01) (2,8.586277e-01) (2.5,8.078779e-01) (3,6.379677e-01) (3.5,5.036224e-01)
(4,4.648982e-01) (4.5,3.940108e-01) (5,3.487298e-01) (5.5,2.791646e-01) (6,2.300323e-01) (6.5,1.472535e-01)
(7,1.172597e-01) (7.5,9.538726e-02) (8,7.377173e-02) (8.5,5.662153e-02) (9,4.267085e-02) (9.5,2.941934e-02)
};
\end{semilogyaxis}
\end{tikzpicture}
}
\subfigure[$\delta=80$]{
\begin{tikzpicture}[scale=0.8]
\begin{semilogyaxis}[xmin=1, xmax=9.5,domain=1:9.5,ymin=1e-2, ymax=2,domain=1e-2:2,
    xlabel={Wavenumber},
    ylabel={$E_\mathrm{SIMDNP}$},
    grid=major,
    legend pos=south west,
    legend entries={$N_s=10$, $N_s=50$, $N_s=100$, $N_s=500$, $N_s=1000$, $N_s=2000$, $N_s=4000$, RLA},
]
\addplot coordinates {
(1,8.916994e-01) (1.5,8.816462e-01) (2,8.608212e-01) (2.5,8.102263e-01) (3,6.576749e-01) (3.5,5.510701e-01)
(4,7.130147e-01) (4.5,6.882500e-01) (5,8.861811e-01) (5.5,9.138538e-01) (6,1.004222e+00) (6.5,1.011319e+00)
(7,1.051096e+00) (7.5,1.045514e+00) (8,1.040457e+00) (8.5,1.061268e+00) (9,1.082810e+00) (9.5,1.061291e+00)
};
\addplot coordinates {
(1,8.904751e-01) (1.5,8.807835e-01) (2,8.584640e-01) (2.5,8.096984e-01) (3,6.511993e-01) (3.5,5.582816e-01)
(4,6.192643e-01) (4.5,5.360598e-01) (5,4.752015e-01) (5.5,4.426302e-01) (6,5.021888e-01) (6.5,4.897863e-01)
(7,5.174652e-01) (7.5,5.036462e-01) (8,4.945731e-01) (8.5,4.969310e-01) (9,5.173578e-01) (9.5,5.071919e-01)
};
\addplot coordinates {
(1,8.904815e-01) (1.5,8.807848e-01) (2,8.583188e-01) (2.5,8.091974e-01) (3,6.509182e-01) (3.5,5.573313e-01)
(4,6.083074e-01) (4.5,5.264440e-01) (5,4.393380e-01) (5.5,3.661674e-01) (6,3.686904e-01) (6.5,3.504226e-01)
(7,3.766917e-01) (7.5,3.685694e-01) (8,3.751286e-01) (8.5,3.898838e-01) (9,4.265754e-01) (9.5,4.144210e-01)
};
\addplot coordinates {
(1,8.904780e-01) (1.5,8.809060e-01) (2,8.586733e-01) (2.5,8.096257e-01) (3,6.508167e-01) (3.5,5.534875e-01)
(4,5.832718e-01) (4.5,5.003564e-01) (5,4.081531e-01) (5.5,3.511362e-01) (6,3.234143e-01) (6.5,2.869333e-01)
(7,2.997024e-01) (7.5,3.018333e-01) (8,3.162527e-01) (8.5,3.169254e-01) (9,3.539131e-01) (9.5,3.437252e-01)
};
\addplot coordinates {
(1,8.904857e-01) (1.5,8.809540e-01) (2,8.588561e-01) (2.5,8.097963e-01) (3,6.515764e-01) (3.5,5.542804e-01)
(4,5.819264e-01) (4.5,4.985734e-01) (5,4.136258e-01) (5.5,3.677501e-01) (6,3.233609e-01) (6.5,2.839182e-01)
(7,2.964462e-01) (7.5,3.044848e-01) (8,3.213876e-01) (8.5,3.200056e-01) (9,3.557385e-01) (9.5,3.474439e-01)
};
\addplot coordinates {
(1,8.905024e-01) (1.5,8.809988e-01) (2,8.589635e-01) (2.5,8.099856e-01) (3,6.524059e-01) (3.5,5.558059e-01)
(4,5.815584e-01) (4.5,4.959341e-01) (5,4.151393e-01) (5.5,3.677139e-01) (6,3.144800e-01) (6.5,2.722624e-01)
(7,2.831758e-01) (7.5,2.961605e-01) (8,3.131020e-01) (8.5,3.071098e-01) (9,3.446144e-01) (9.5,3.382348e-01)
};
\addplot coordinates {
(1,8.904878e-01) (1.5,8.809397e-01) (2,8.588778e-01) (2.5,8.099167e-01) (3,6.517177e-01) (3.5,5.551635e-01)
(4,5.821847e-01) (4.5,4.954680e-01) (5,4.097290e-01) (5.5,3.494091e-01) (6,2.990897e-01) (6.5,2.561392e-01)
(7,2.645941e-01) (7.5,2.699211e-01) (8,2.878308e-01) (8.5,2.836239e-01) (9,3.240500e-01) (9.5,3.141806e-01)
};
\addplot coordinates {
(1,8.905047e-01) (1.5,8.807221e-01) (2,8.586277e-01) (2.5,8.078779e-01) (3,6.379677e-01) (3.5,5.036224e-01)
(4,4.648982e-01) (4.5,3.940108e-01) (5,3.487298e-01) (5.5,2.791646e-01) (6,2.300323e-01) (6.5,1.472535e-01)
(7,1.172597e-01) (7.5,9.538726e-02) (8,7.377173e-02) (8.5,5.662153e-02) (9,4.267085e-02) (9.5,2.941934e-02)
};
\end{semilogyaxis}
\end{tikzpicture}
}
 \caption{Error for the {\bf SIMDNP} algorithm as a function of the maximum wavenumber used in the reconstruction for different levels of noise: (a) $\delta=10$, (b) $\delta=20$, (c) $\delta=40$ and (d) $\delta=80$. In each plot, each line represents the error $E_\mathrm{SIMDNP}$ using a different number of samples $N_s$. We also include the error for the {\bf RLA} (with data obtained with $\eta=0) $in each plot as a matter of comparison and as a benchmark for the best possible approximation for our problem.}
 \label{error_exp1_SIMDNP}
\end{figure}

In Figures \ref{error_exp1_MIMDNP} and \ref{error_exp1_SIMDNP} we report  the errors 
\begin{equation*}
E_\mathrm{SIMDNP}=\frac{\|q_\mathrm{MIMDNP}-\qb\|_2}{\|\qb\|_2}
\end{equation*}
and 
\begin{equation*}
E_\mathrm{SIMDNP}=\frac{\|q_\mathrm{SIMDNP}-\qb\|_2}{\|\qb\|_2}
\end{equation*}
as a function of the maximum wavenumber $\kappa_Q$ used in the reconstruction (see~\ref{s:setup} for $\kappa_Q$), for different noise level $\delta$.  As a matter of comparison, we apply the {\bf RLA} on scattered data generated by the scatterer $q$ with no background noisy medium ($\eta=0$) to obtain the approximation $q_\mathrm{RLA}$. The error $\|q_\mathrm{RLA}-\qb\|_2/\|\qb\|_2$ of the solution using the {\bf RLA} is presented in each one of the images in Figures \ref{error_exp1_MIMDNP} and \ref{error_exp1_SIMDNP}. 

\begin{figure}[!htp]
\centering
\subfigure[$\delta=10$ and $N_s=10$]{
\includegraphics[width=0.3\textwidth]{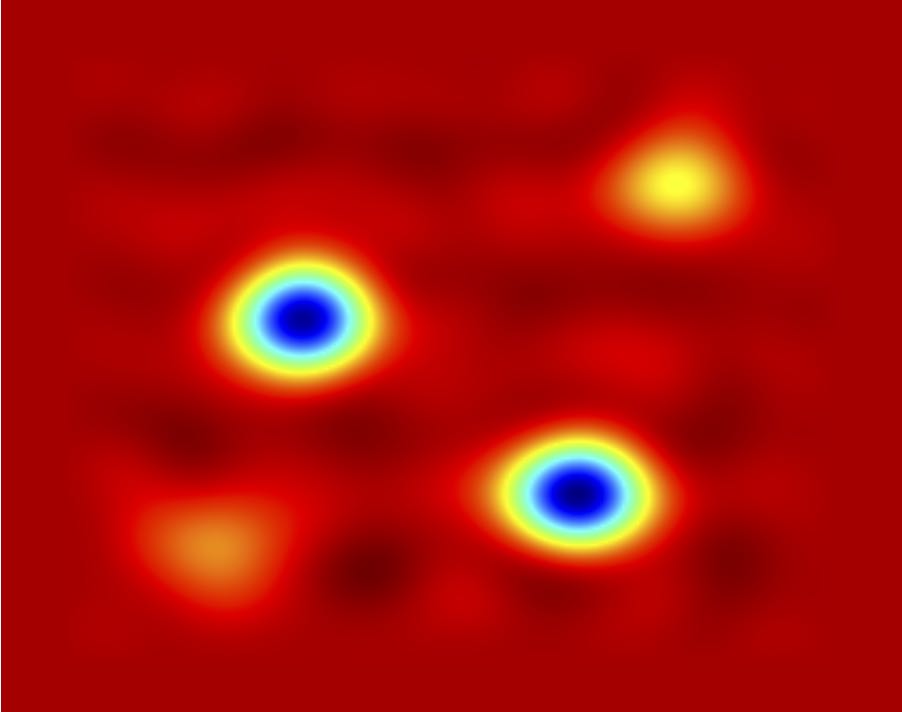}\label{sol_npid1_d10_ns10}
}
\subfigure[$\delta=10$ and $N_s=100$]{
\includegraphics[width=0.3\textwidth]{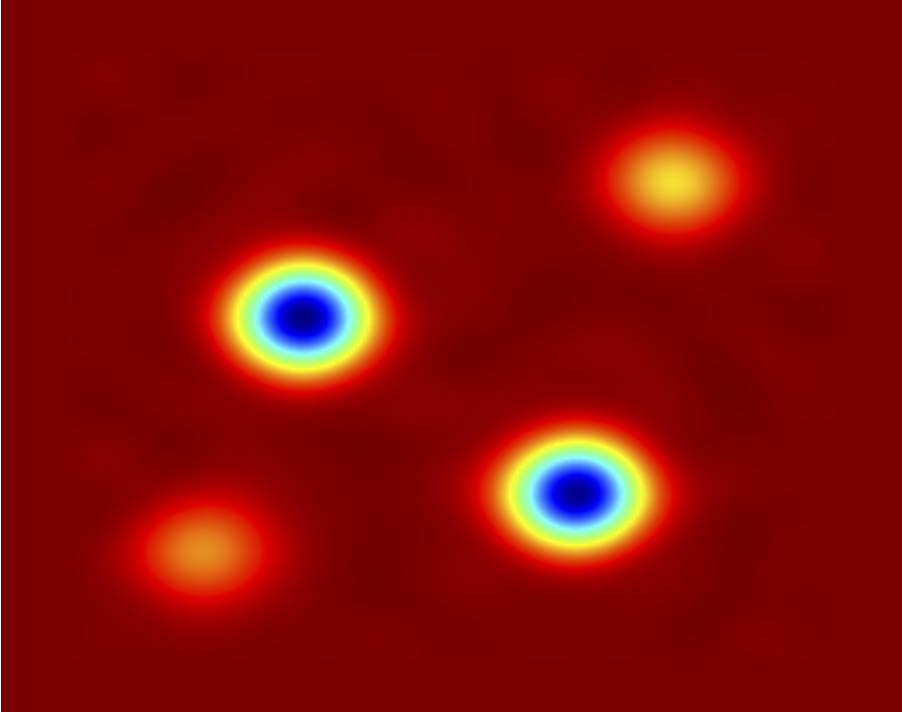}\label{sol_npid1_d_10_ns100}
}
\subfigure[$\delta=10$ and $N_s=1000$]{
\includegraphics[width=0.3\textwidth]{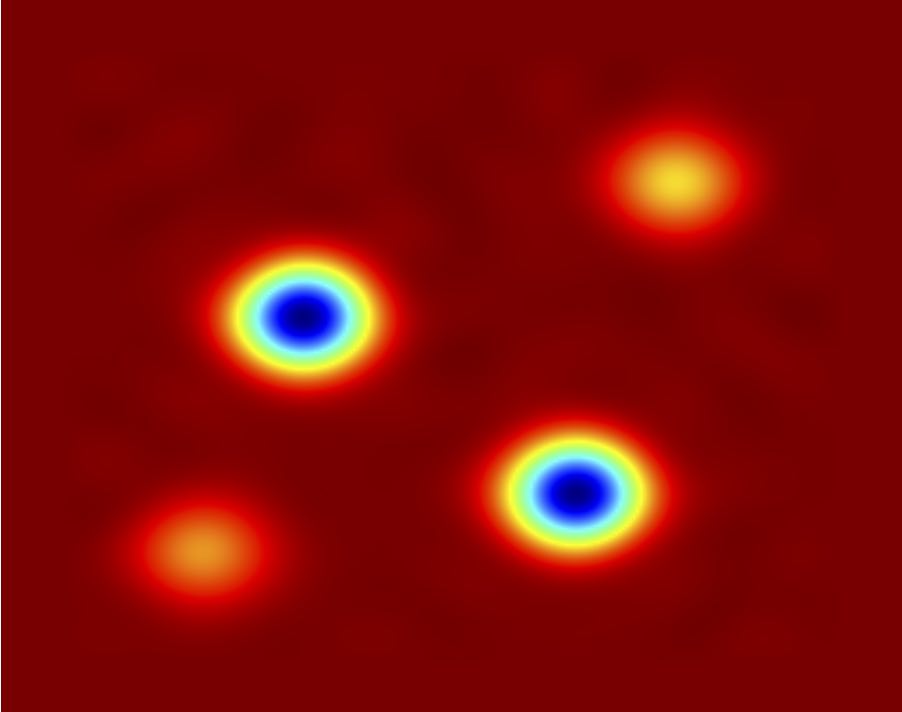}\label{sol_npid1_d10_ns1000}
}

\subfigure[$\delta=20$ and $N_s=10$]{
\includegraphics[width=0.3\textwidth]{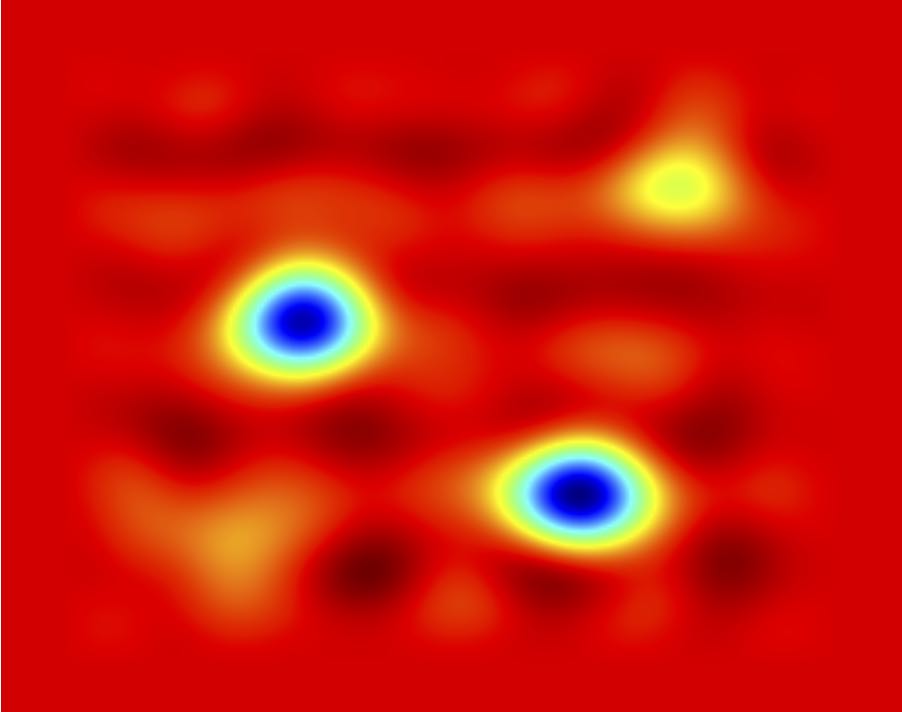}\label{sol_npid1_d20_ns10}
}
\subfigure[$\delta=20$ and $N_s=100$]{
\includegraphics[width=0.3\textwidth]{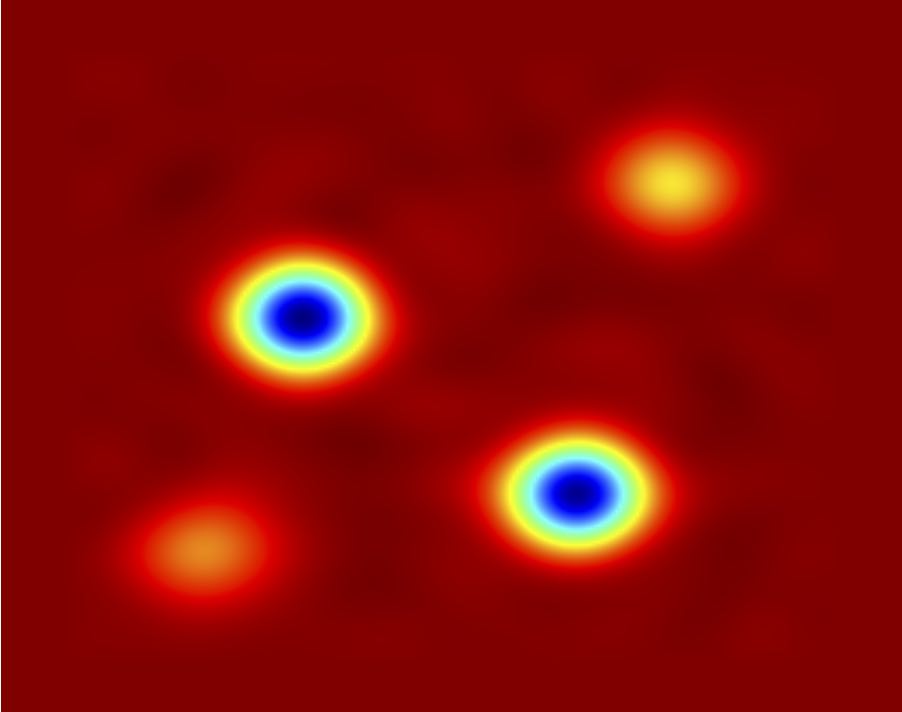}\label{sol_npid1_d20_ns100}
}
\subfigure[$\delta=20$ and $N_s=1000$]{
\includegraphics[width=0.3\textwidth]{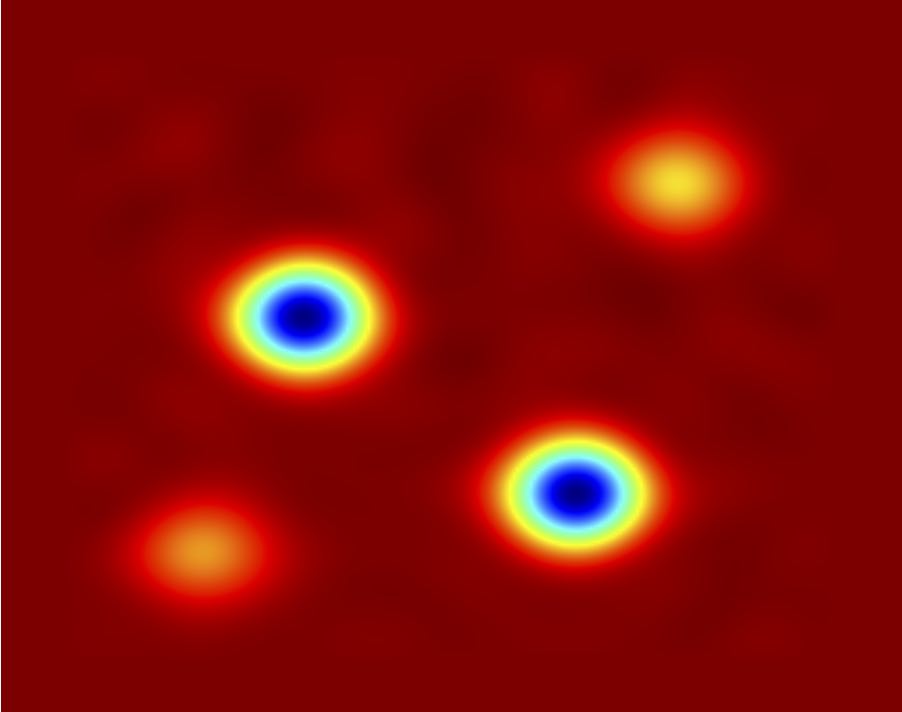}\label{sol_npid1_d20_ns1000}
}

\subfigure[$\delta=40$ and $N_s=10$]{
\includegraphics[width=0.3\textwidth]{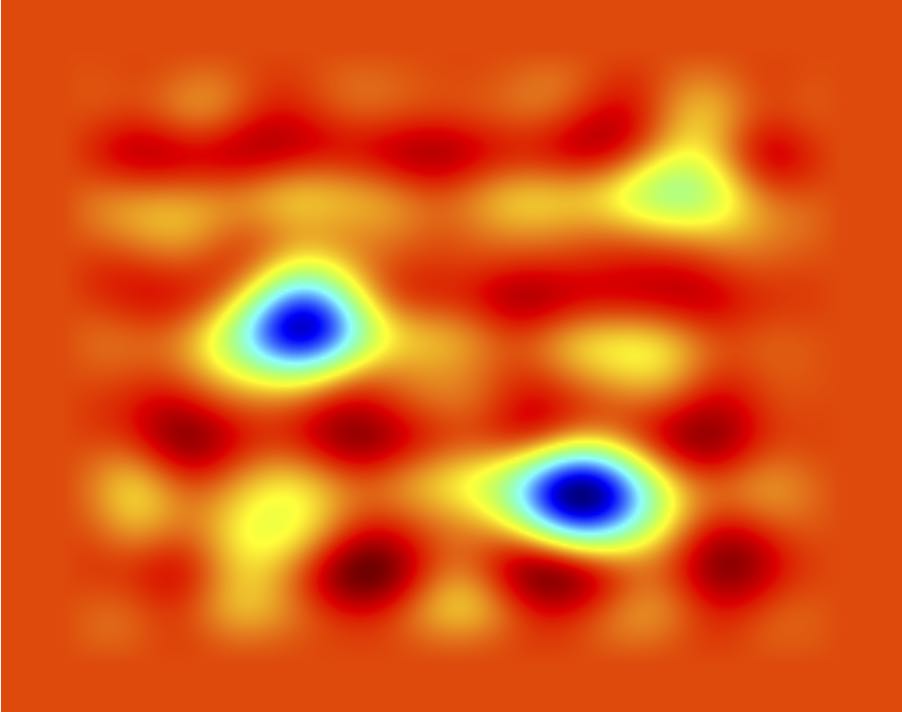}\label{sol_npid1_d40_ns10}
}
\subfigure[$\delta=40$ and $N_s=500$]{
\includegraphics[width=0.3\textwidth]{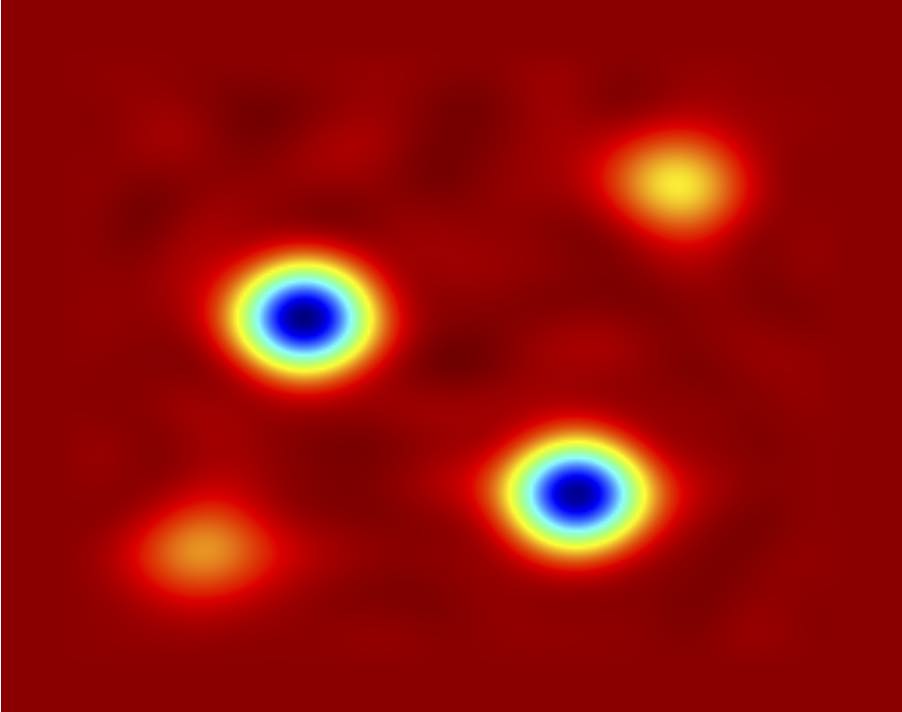}\label{sol_npid1_d40_ns500}
}
\subfigure[$\delta=40$ and $N_s=4000$]{
\includegraphics[width=0.3\textwidth]{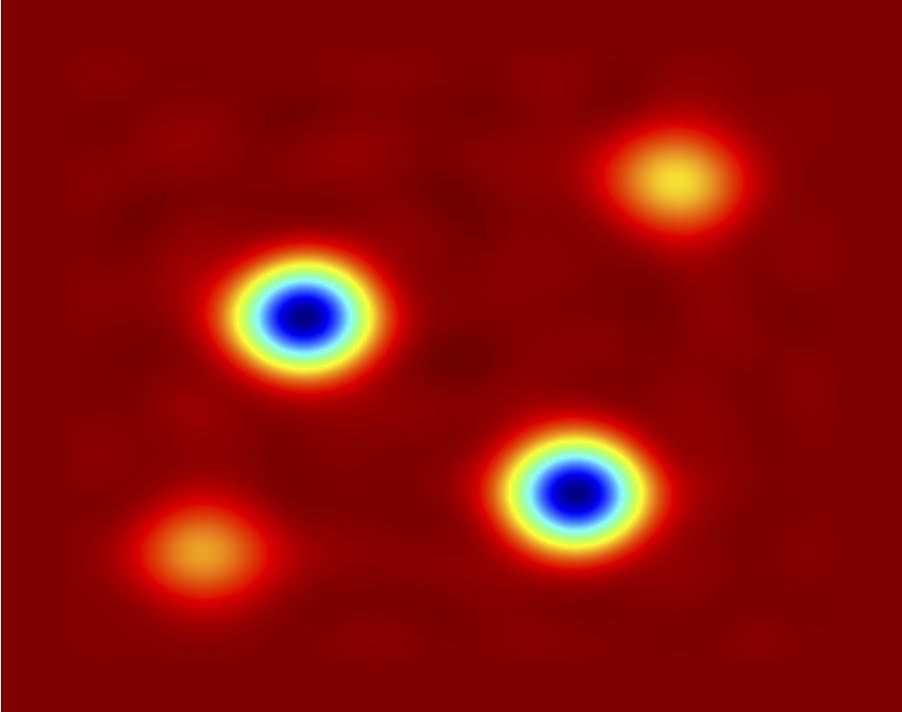}\label{sol_npid1_d40_ns4000}
}

\subfigure[$\delta=80$ and $N_s=10$]{
\includegraphics[width=0.3\textwidth]{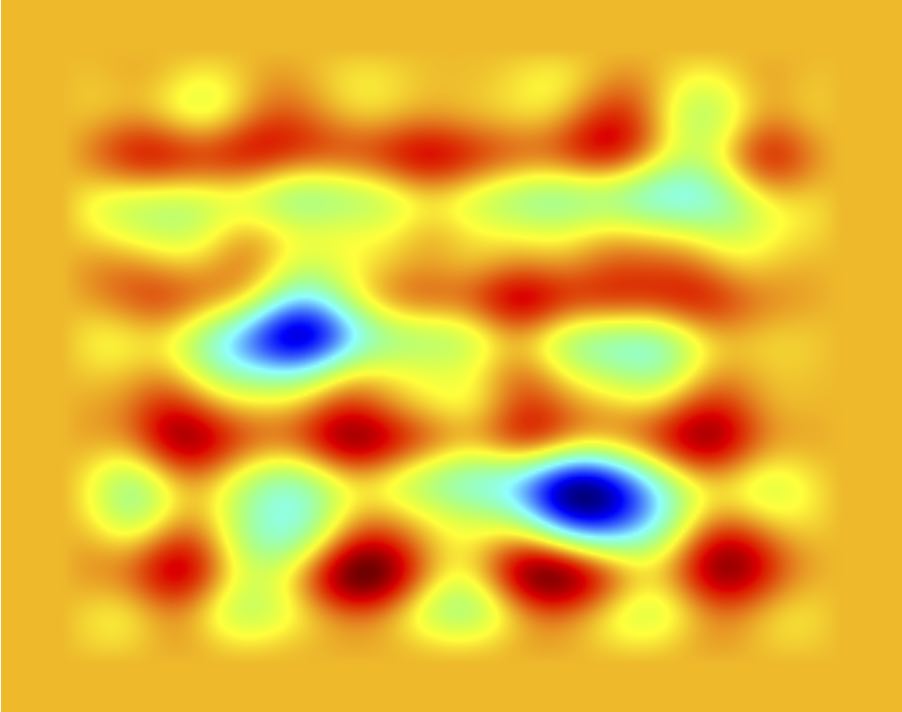}\label{sol_npid1_d80_ns10}
}
\subfigure[$\delta=80$ and $N_s=500$]{
\includegraphics[width=0.3\textwidth]{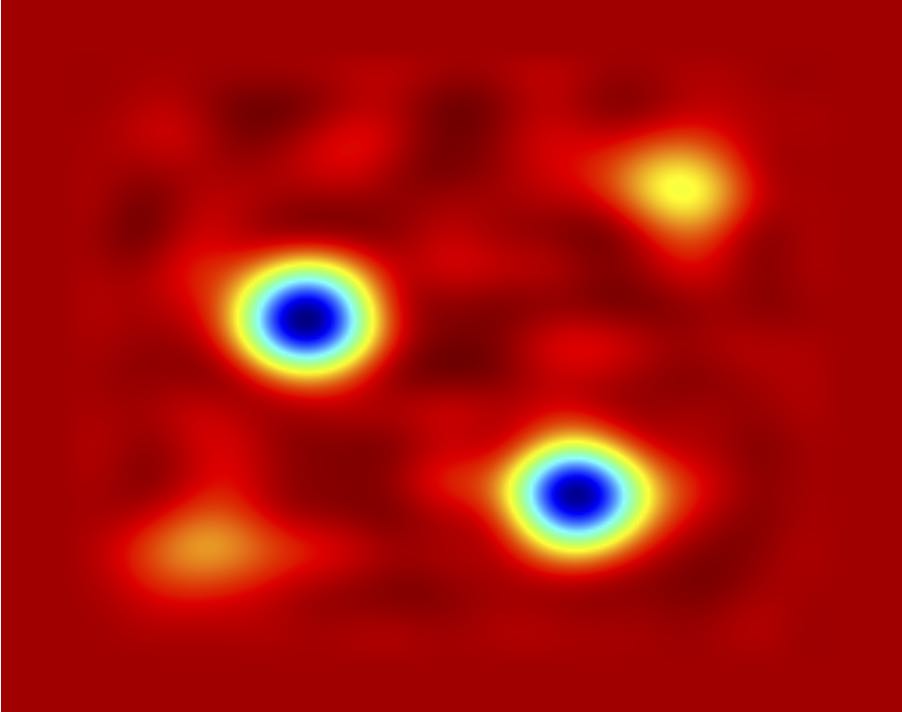}\label{sol_npid1_d80_ns500}
}
\subfigure[$\delta=80$ and $N_s=4000$]{
\includegraphics[width=0.3\textwidth]{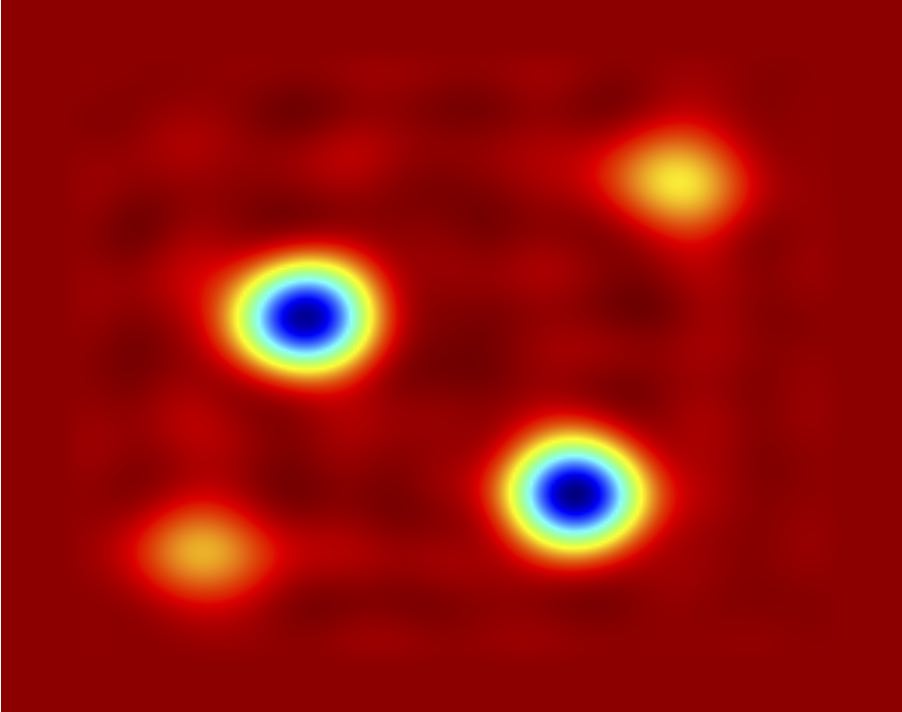}\label{sol_npid1_d80_ns4000}
}
\caption{Reconstruction of $\qb$ for Example \ref{example:MDNP}. The {\bf MIMDNP} algorithm is used to reconstruct the domain $\qb$ in the presence of a background medium. From top to bottom, we present the solution $q_\mathrm{MIMDNP}$ when the background medium is generated using the parameter $\delta=10$, $20$, $40$ and $80$, with different number of samples $N_s$ of the background medium.}\label{ex1_MIMDNP_results}
\end{figure}

\begin{figure}[!htp]
\centering
\subfigure[$\delta=10$ and $N_s=10$]{
\includegraphics[width=0.3\textwidth]{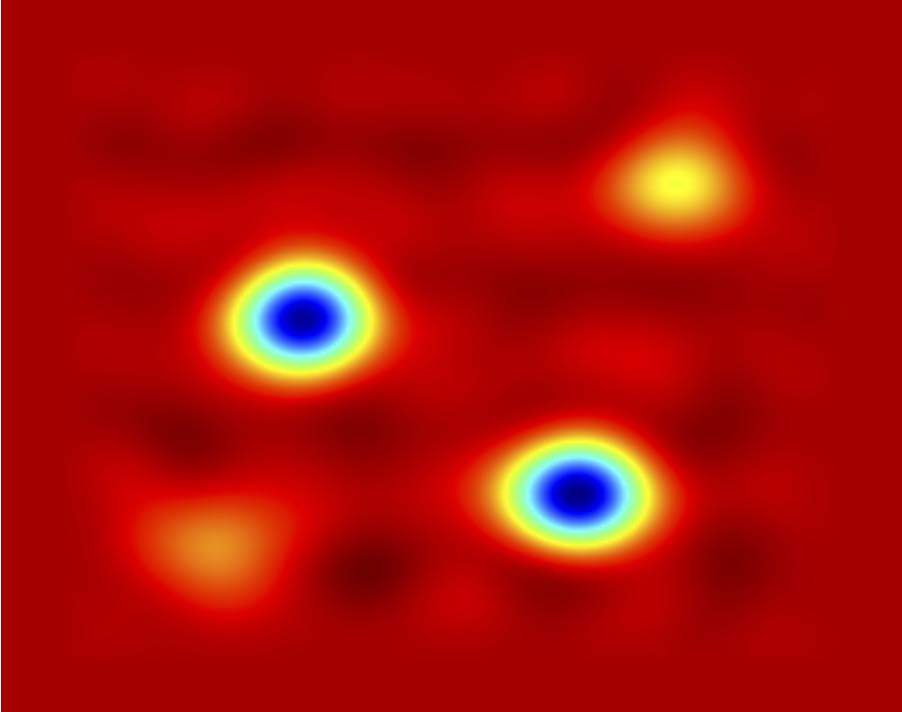}\label{sol_npid2_d10_ns10}
}
\subfigure[$\delta=10$ and $N_s=100$]{
\includegraphics[width=0.3\textwidth]{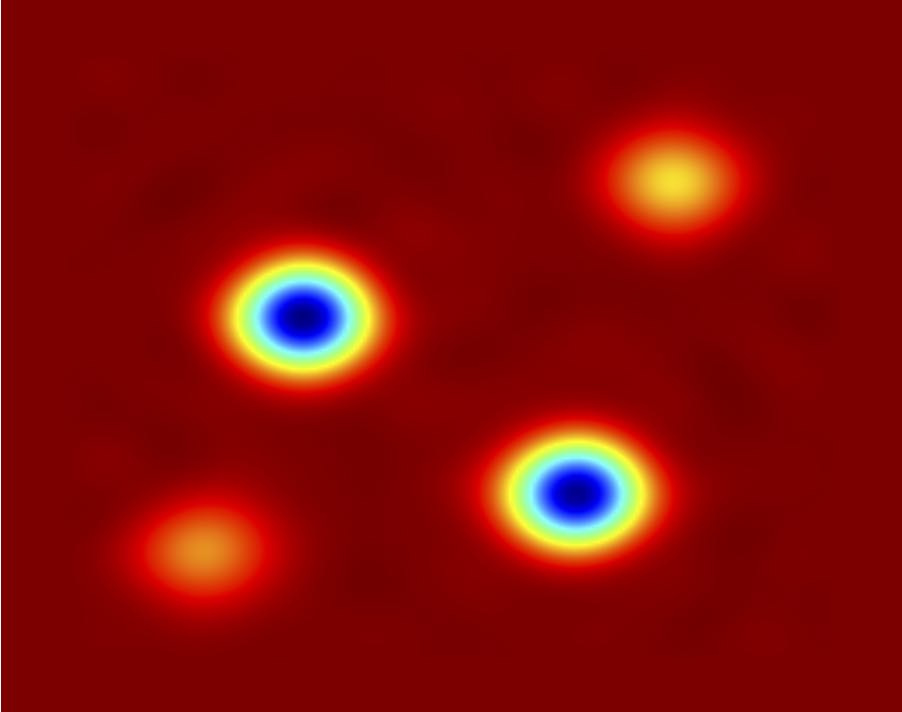}\label{sol_npid2_d_10_ns100}
}
\subfigure[$\delta=10$ and $N_s=1000$]{
\includegraphics[width=0.3\textwidth]{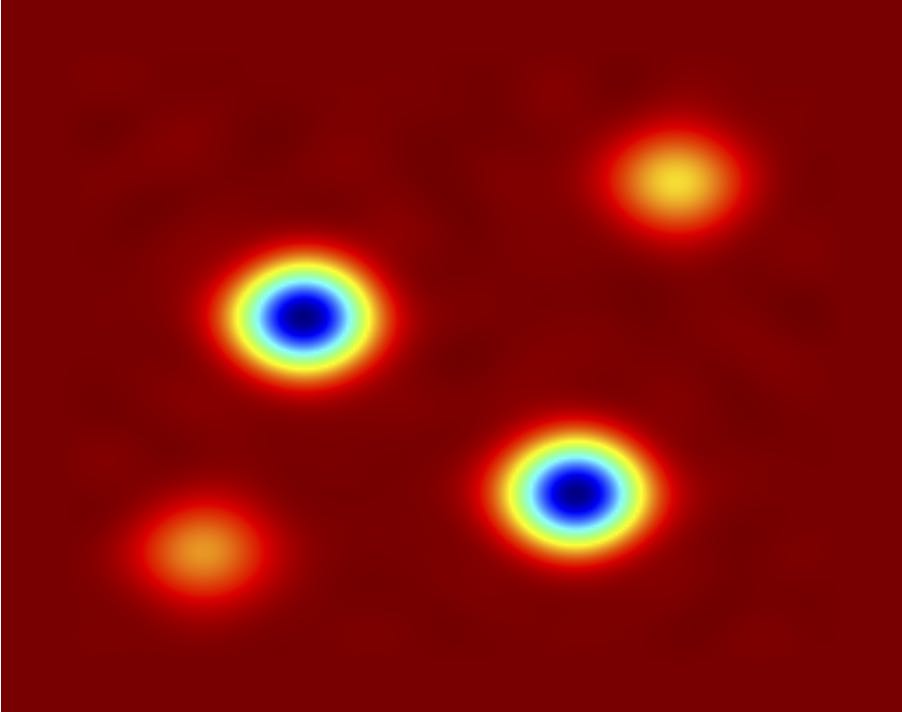}\label{sol_npid2_d10_ns1000}
}

\subfigure[$\delta=20$ and $N_s=10$]{
\includegraphics[width=0.3\textwidth]{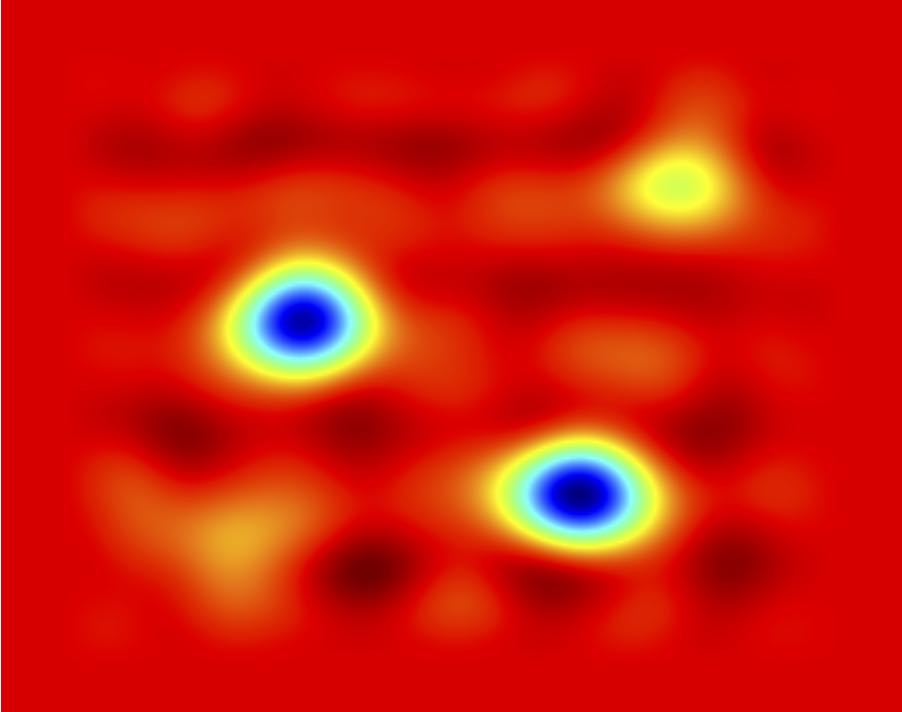}\label{sol_npid2_d20_ns10}
}
\subfigure[$\delta=20$ and $N_s=100$]{
\includegraphics[width=0.3\textwidth]{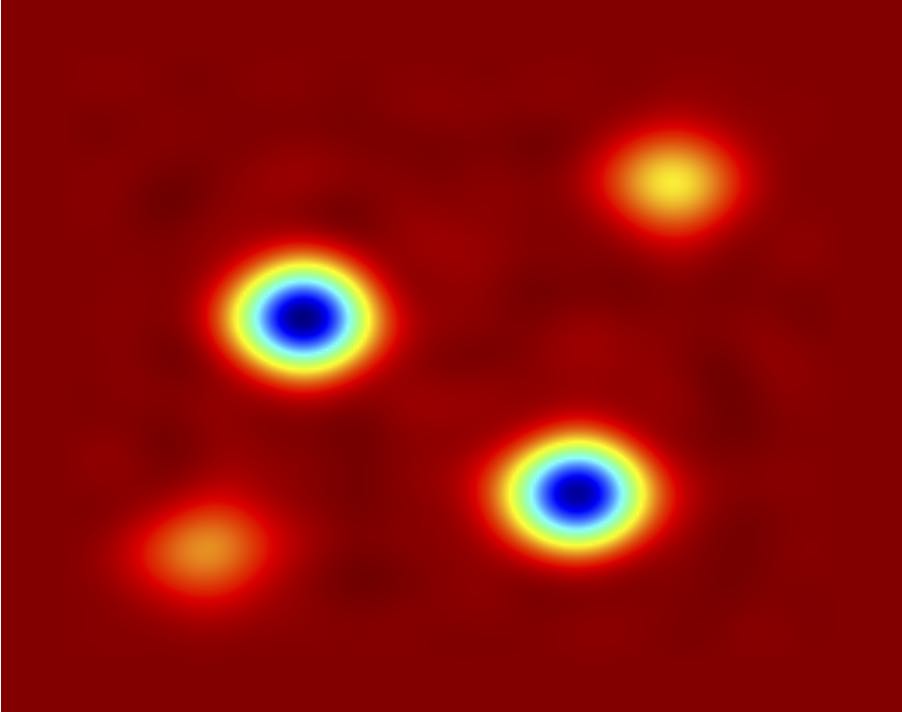}\label{sol_npid2_d20_ns100}
}
\subfigure[$\delta=20$ and $N_s=1000$]{
\includegraphics[width=0.3\textwidth]{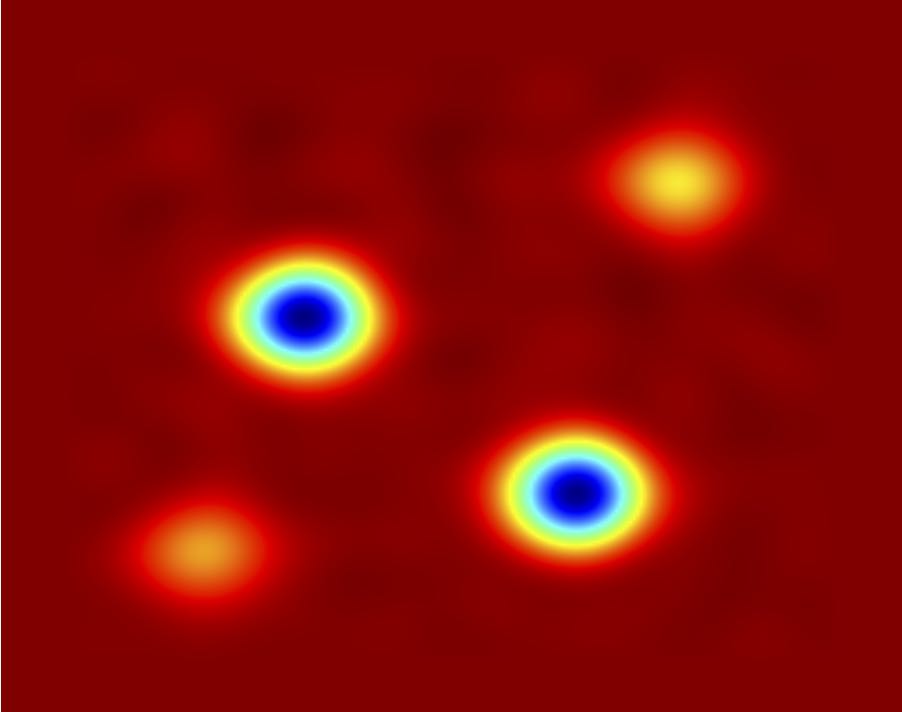}\label{sol_npid2_d20_ns1000}
}

\subfigure[$\delta=40$ and $N_s=10$]{
\includegraphics[width=0.3\textwidth]{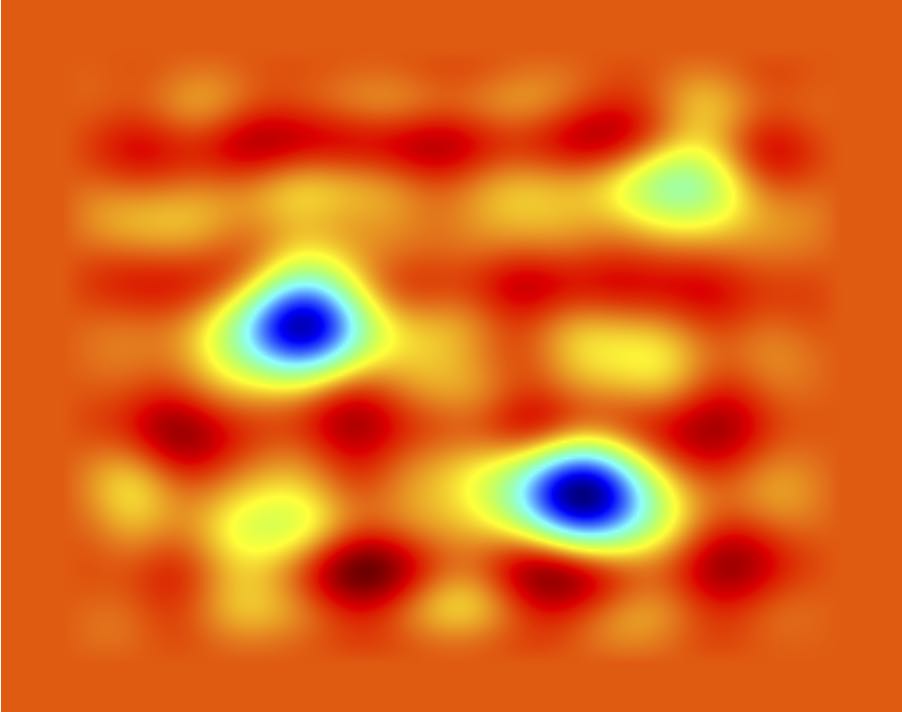}\label{sol_npid2_d40_ns10}
}
\subfigure[$\delta=40$ and $N_s=500$]{
\includegraphics[width=0.3\textwidth]{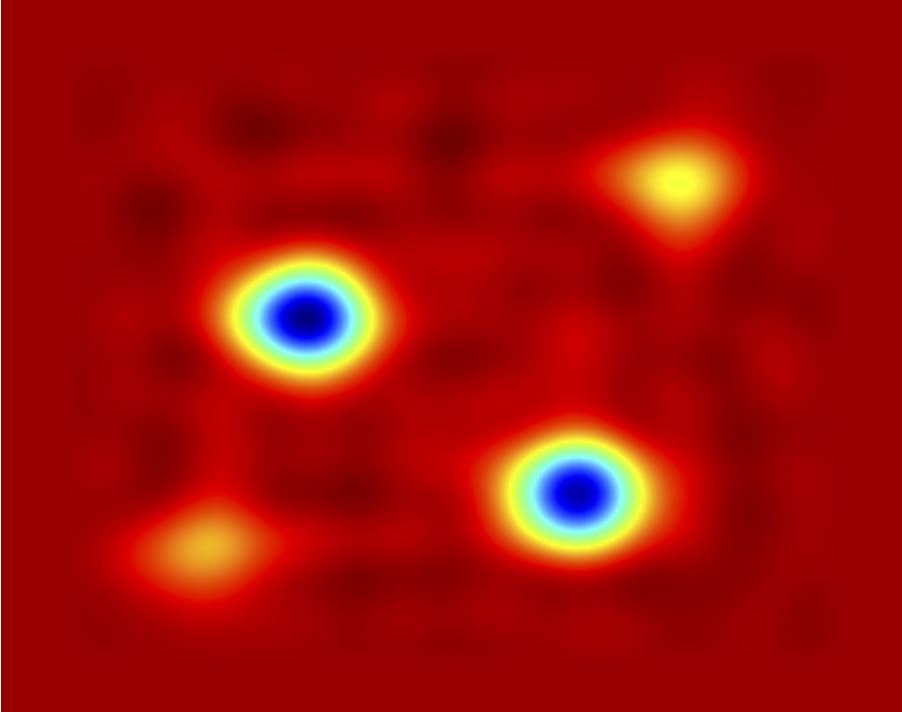}\label{sol_npid2_d40_ns500}
}
\subfigure[$\delta=40$ and $N_s=4000$]{
\includegraphics[width=0.3\textwidth]{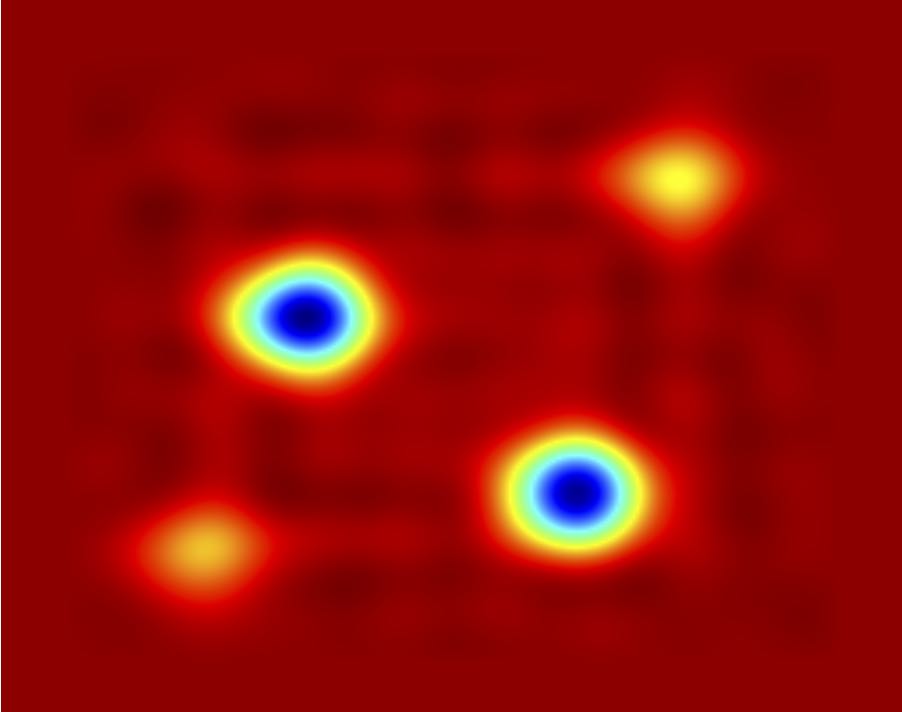}\label{sol_npid2_d40_ns4000}
}

\subfigure[$\delta=80$ and $N_s=10$]{
\includegraphics[width=0.3\textwidth]{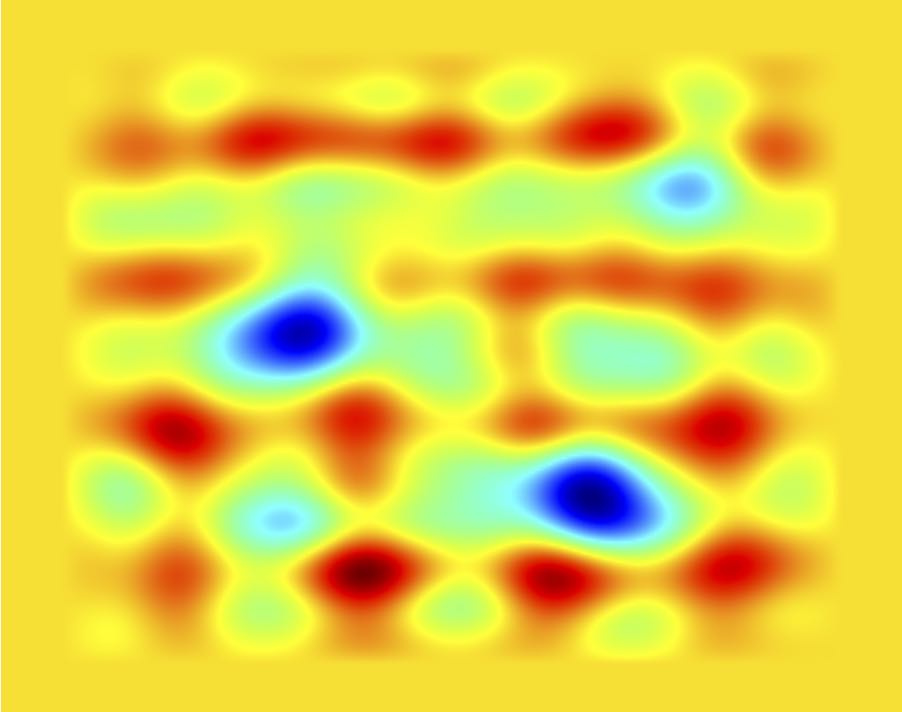}\label{sol_npid2_d80_ns10}
}
\subfigure[$\delta=80$ and $N_s=500$]{
\includegraphics[width=0.3\textwidth]{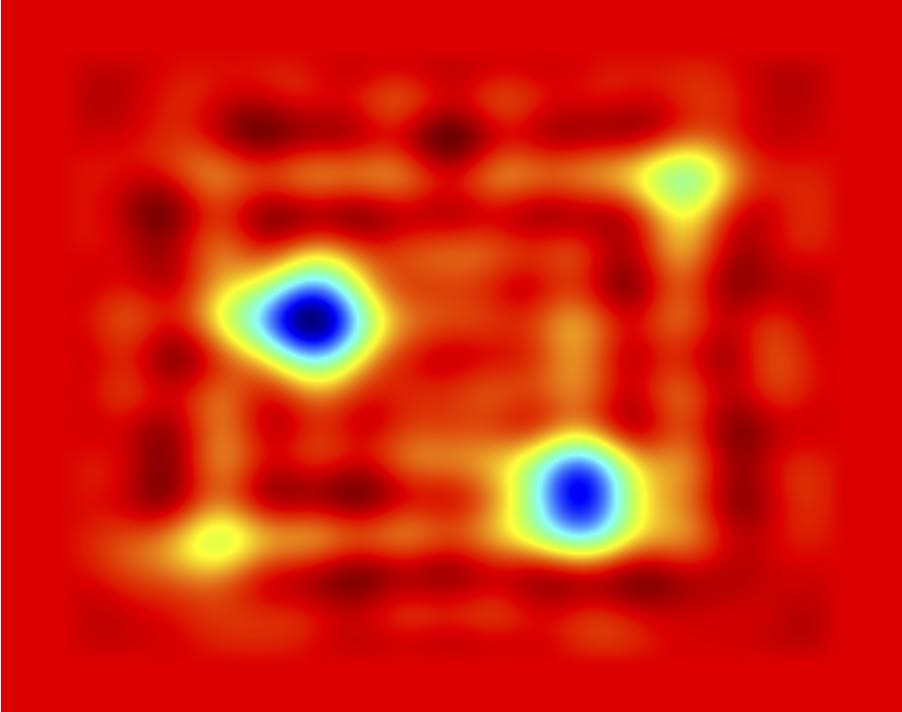}\label{sol_npid2_d80_ns500}
}
\subfigure[$\delta=80$ and $N_s=4000$]{
\includegraphics[width=0.3\textwidth]{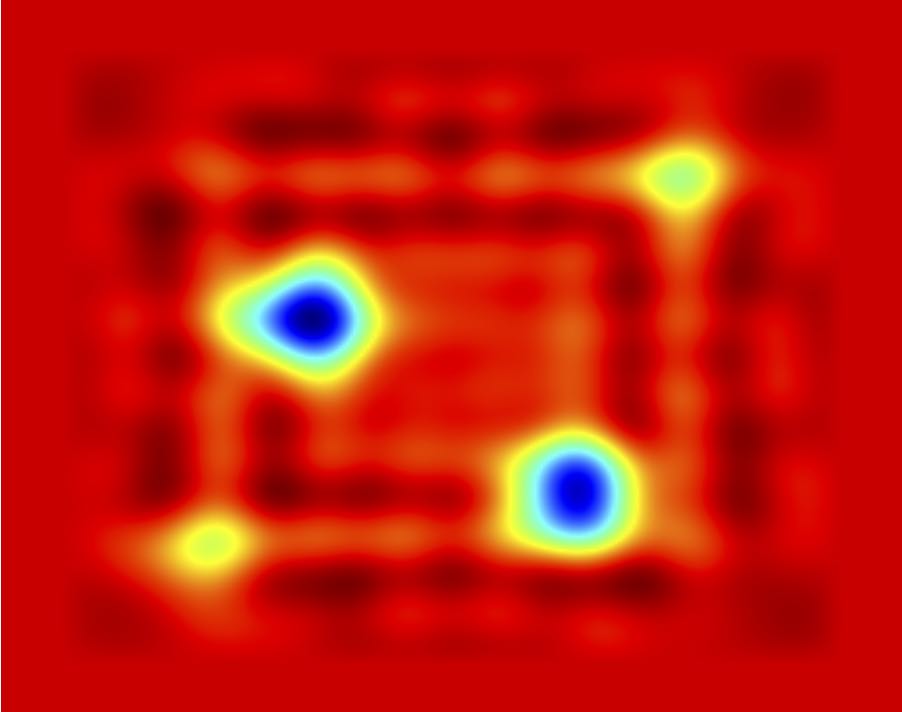}\label{sol_npid2_d80_ns4000}
}
\caption{Reconstruction of $\qb$ for Example \ref{example:MDNP}. The {\bf SIMDNP} algorithm is used to reconstruct the domain $\qb$ in the presence of a background medium. From top to bottom, we present the solution $q_\mathrm{SIMDNP}$ when the background medium is generated using the parameter $\delta=10$, $20$, $40$ and $80$, with different number of samples $N_s$ of the background medium.}\label{ex1_SIMDNP_results}
\end{figure}

In Figures \ref{ex1_MIMDNP_results} and \ref{ex1_SIMDNP_results}, respectively, we present the reconstructions obtained using the algorithms {\bf MIMDNP} and {\bf SIMDNP} for the different levels of noise $\delta=10$, $20$, $40$ and $80$. For $\delta=10$ and $20$, we present the reconstructions using $N_s=10$, $100$ and $1000$, and for $\delta=40$ and $80$, we present the results using $N_s=10$, $500$ and $4000$.

This example confirms that {\bf MIMDNP} has better accuracy than {\bf SISDNP} and the extra cost is justified. From the experimental results in Figure \ref{error_exp1_MIMDNP}, it is clear that if we use a sufficiently large number of samples, the solution of $q_\mathrm{MIMDNP}$ should converge asymptotically to $\qb$. We can also see, from the results in Figure \ref{error_exp1_SIMDNP}, that when the hypotheses of Theorem \ref{thm:SIMDNP} are satisfied, meaning that the product of the domain with the square of the wavenumber is small enough, the forward operator becomes approximately linear, and the approximation $q_\mathrm{SIMDNP}$ is very close to the approximation obtained by the {\bf RLA} (with data obtained with $\eta=0$). 

\subsection{Data from single $\eta$ and inversion for both $q$ and $\eta$}
\label{example:SDP}

In this example, we compare the solutions by {\bf SISDP} and {\bf RLA} in three test cases: 
\begin{enumerate}[label=(\alph*)]
\item $q$ and $\eta$ have the same prior;
\item $q$ and $\eta$ have completely separated priors;
\item $q$ and $\eta$ have partially overlapping priors.
\end{enumerate}

Similar to the discretization for $q$~\eqref{eq:sineseries}, the random field $\eta$ is given by
\begin{equation}\label{e:spectral-eta}
\eta(x_1,x_2)=\sum_{m_1,m_2=1}^{M_\eta}\eta_{m_1 m_2}\sin(m_1 x_1)\sin(m_2 x_2), \quad \mbox{with}\quad \eta_{m_1 m_2} = \frac{\zeta}{a_{11}m_1^2+a_{22}m_2^2}.
\end{equation}
Here $a_{11}$ and $a_{22}$ are the diagonal elements of $\mathcal{T}_\eta$, and
$\zeta$ is drawn from a standard normal distribution with mean $0$ and variance $\barsigma$. In this example and the next one, we used $M_\eta=30$.

Our data measurements are obtained at the same receivers as in our previous examples. We have data measurements at wavenumbers $k_j=1+j \delta k$, $\delta k=0.5$, $j=0,\ldots,29$, so that $k_{min}=k_0=1$ and $k_{max}=k_{29}=15$. 

a) In this case, we use the {\bf SISDP} and {\bf RLA} (with the same data)  algorithms for the reconstruction of $\qb$ and $\qp$ in the presence of an isotropic background medium composed of all the frequencies in the chosen spectrum. The background medium is isotropic, with the $\mathcal{T}_\eta$ being the $2\times2$ identity matrix. For the reconstruction of $\qb$, the noise level of $\eta$ is $\barsigma=5$, while for the reconstruction of $\qp$, the noise level is $\barsigma=10$. 

For the regularization of $q$  we use an operator that filters out the higher spatial frequency coefficients from the sine series representation of $q$, so that $m_1+m_2> 2k$. Since we have no specific information regarding the probability distribution of $q$, we set the regularization parameter $\beta=0$. The regularization parameter $\alpha$ for the background medium is obtained using the heuristic described in Appendix B (Algorithm~\ref{alg:find_alpha}). We reconstruct the scatterer using as regularization parameters $\alpha$, $10\alpha$, $50\alpha$, and $100\alpha$.

As we can see in Figures \ref{ex:SISDP:qb:af:results} and \ref{ex:SISDP:qp:af:results} no matter the value of regularization parameter $\alpha$ we are not able to obtain accurate reconstruction of the scatterers $\qb$ and $\qp$. Unfortunately, for this particular case, it is not possible to separate the information obtained in our reconstruction of the scatterer and of the background medium.
 
 \begin{figure}[!htp]
\centering
\subfigure[$\qb+\eta^\ast$ for bumps case]{
\includegraphics[width=0.22\textwidth]{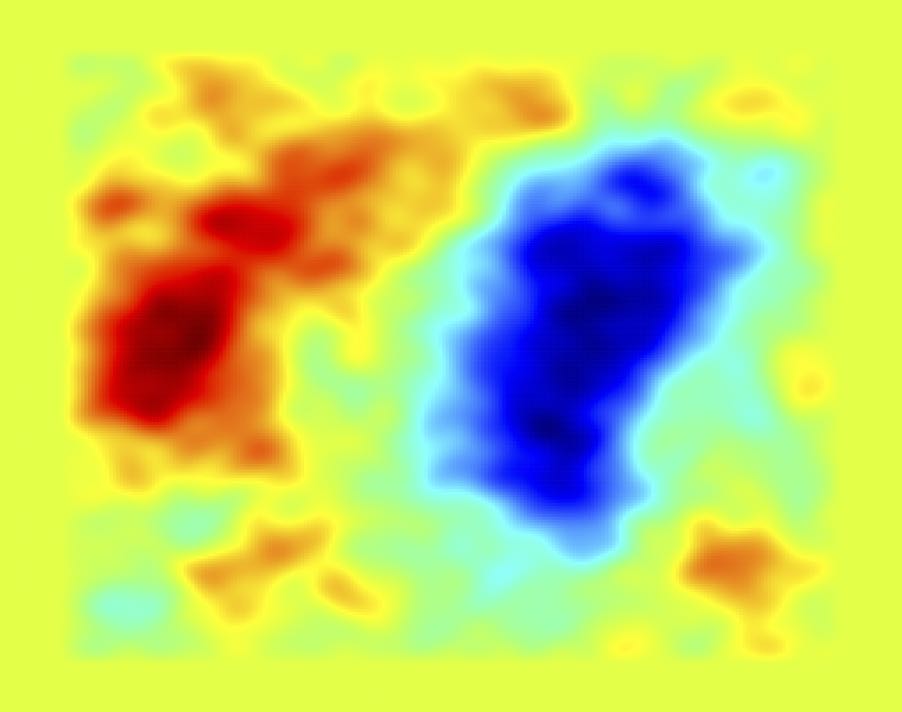}\label{bumps_af_noise_and_q}
}
\subfigure[$q_\mathrm{RLA}$ for $\qb$]{
\includegraphics[width=0.22\textwidth]{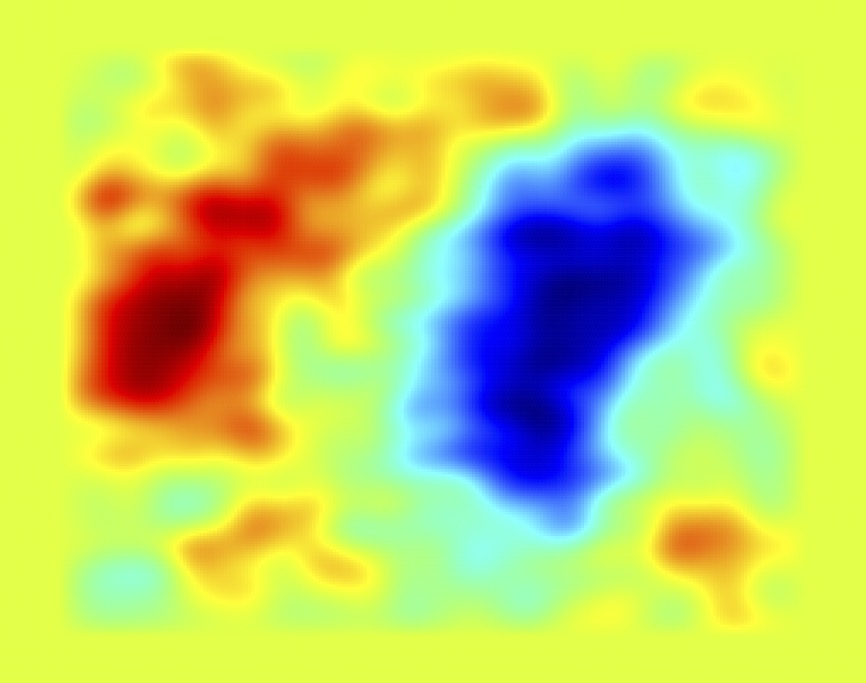}\label{bumps_af_q_newton}
}
\subfigure[$\qp+\eta^\ast$ for plane case]{
\includegraphics[width=0.22\textwidth]{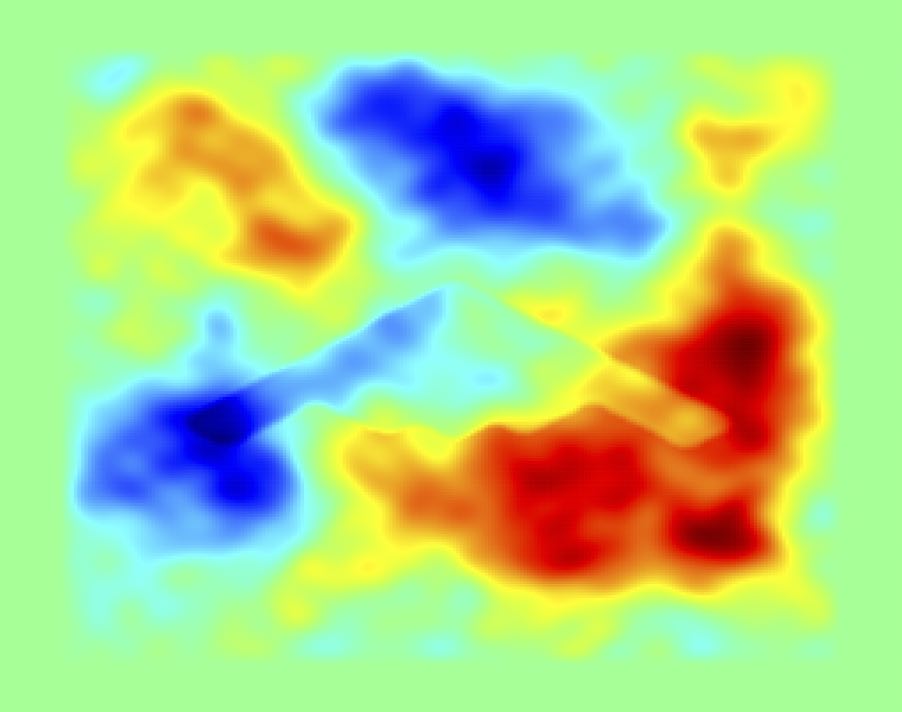}\label{plane_af_noise_and_q}
}
\subfigure[$q_\mathrm{RLA}$ for $\qp$]{
\includegraphics[width=0.22\textwidth]{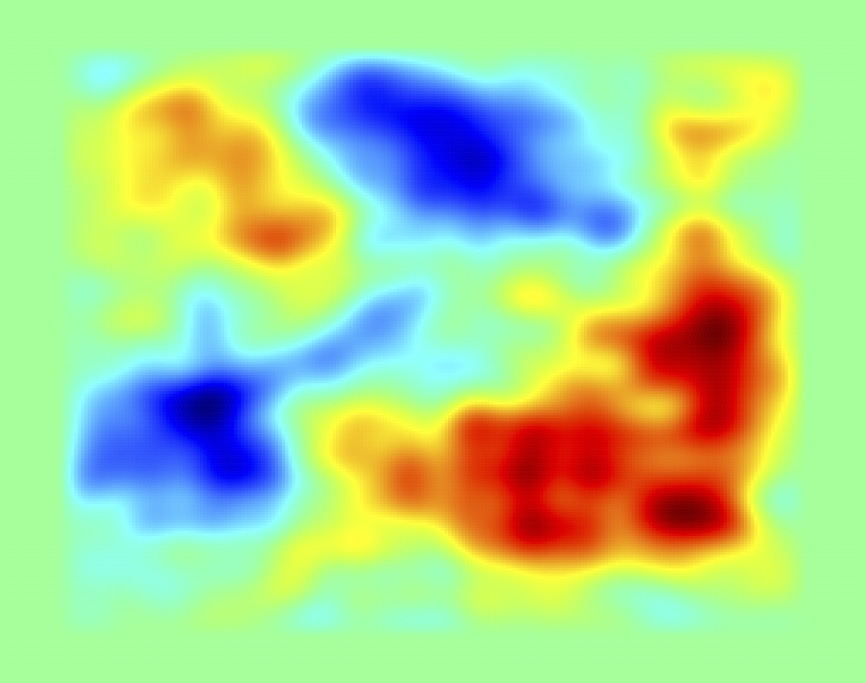}\label{plane_af_q_newton}
}
\caption{The original scatterer in the presence of the background medium and the solution obtained with {\bf RLA} are presented, respectively, in: (a) and (b) for $\qb$, and (c) and (d) for $\qp$.}\label{ex:SISDP:fig1}
\end{figure}

\begin{figure}[!htp]
\centering
\subfigure[$q_\mathrm{SISDP}+\eta_\mathrm{SISDP}$ for $\alpha$]{
\includegraphics[width=0.22\textwidth]{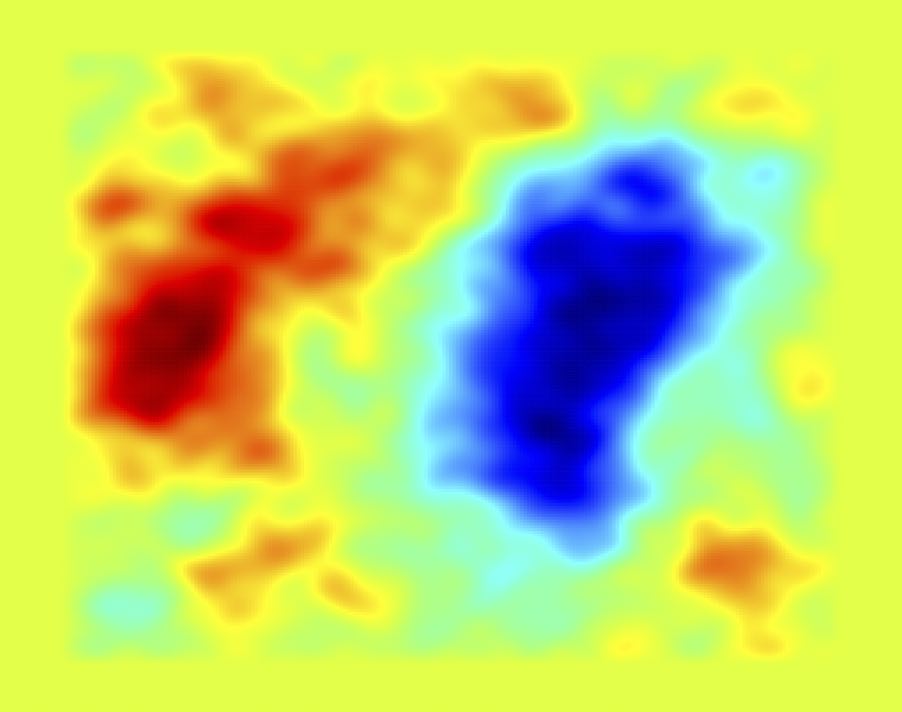}\label{bumps_af_noise_and_q_a1}
}
\subfigure[$q_\mathrm{SISDP}+\eta_\mathrm{SISDP}$ for $10\alpha$]{
\includegraphics[width=0.22\textwidth]{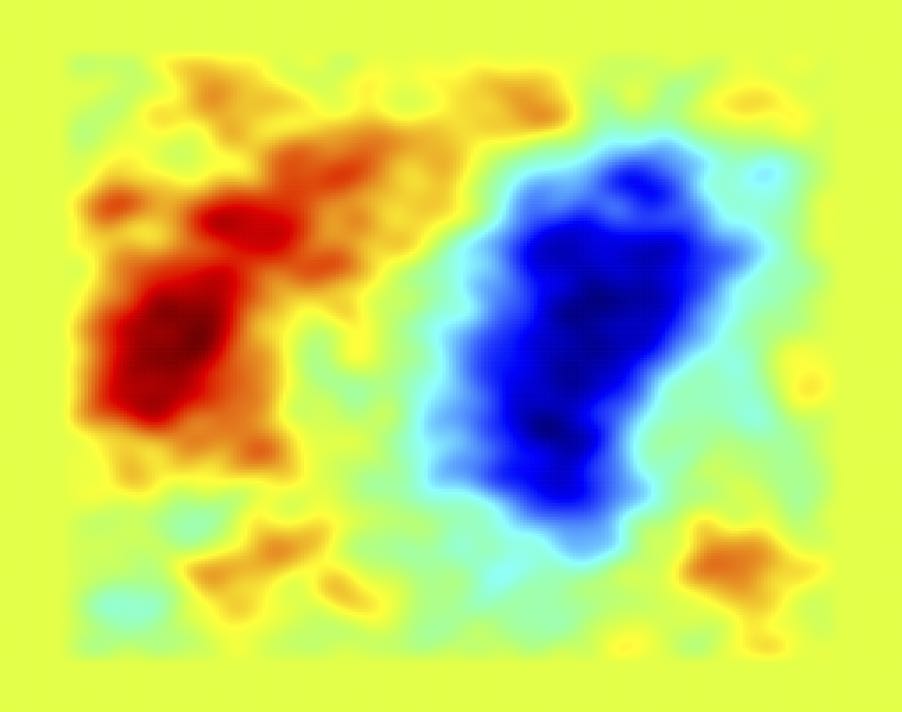}\label{bumps_af_noise_and_q_a10}
}
\subfigure[$q_\mathrm{SISDP}+\eta_\mathrm{SISDP}$ for $50\alpha$]{
\includegraphics[width=0.22\textwidth]{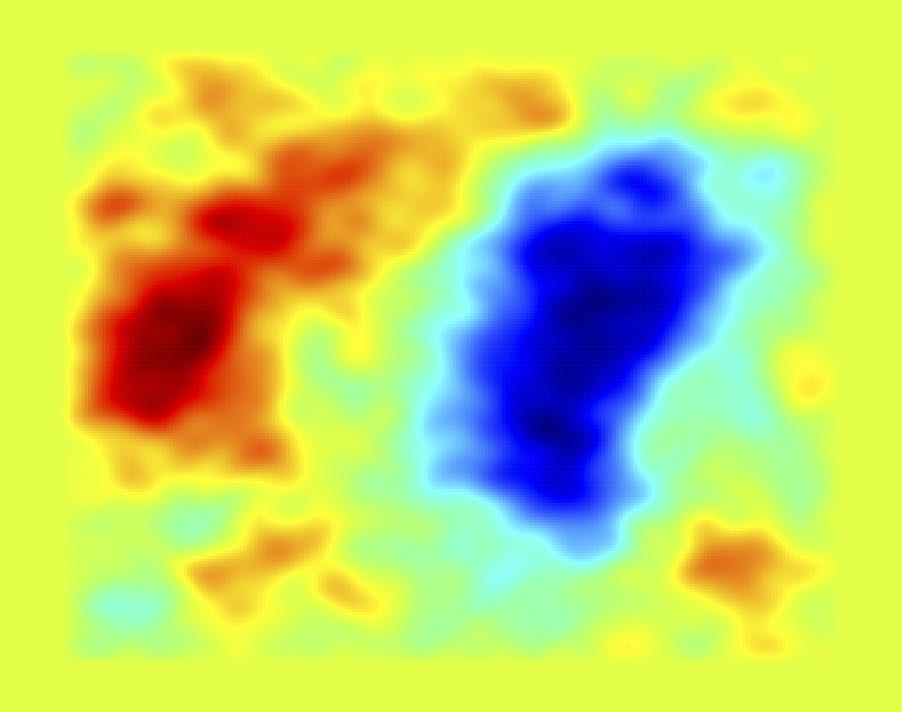}\label{bumps_af_noise_and_q_a50}
}
\subfigure[$q_\mathrm{SISDP}+\eta_\mathrm{SISDP}$ for $100\alpha$]{
\includegraphics[width=0.22\textwidth]{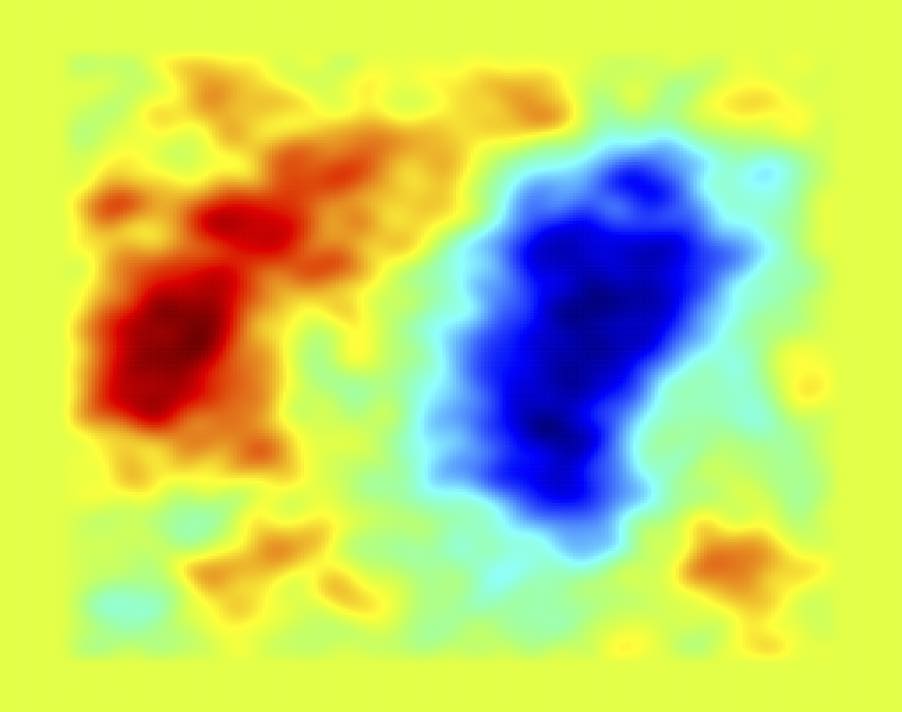}\label{bumps_af_noise_and_q_a100}
}

\subfigure[$q_\mathrm{SISDP}$ for $\alpha$]{
\includegraphics[width=0.22\textwidth]{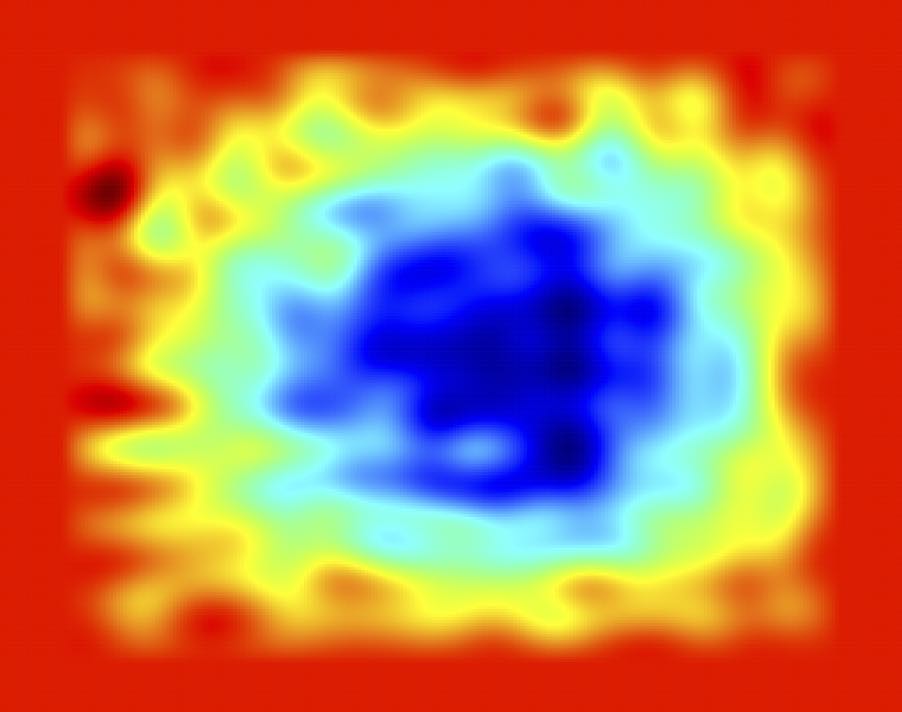}\label{bumps_af_q_a1}
}
\subfigure[$q_\mathrm{SISDP}$ for $10\alpha$]{
\includegraphics[width=0.22\textwidth]{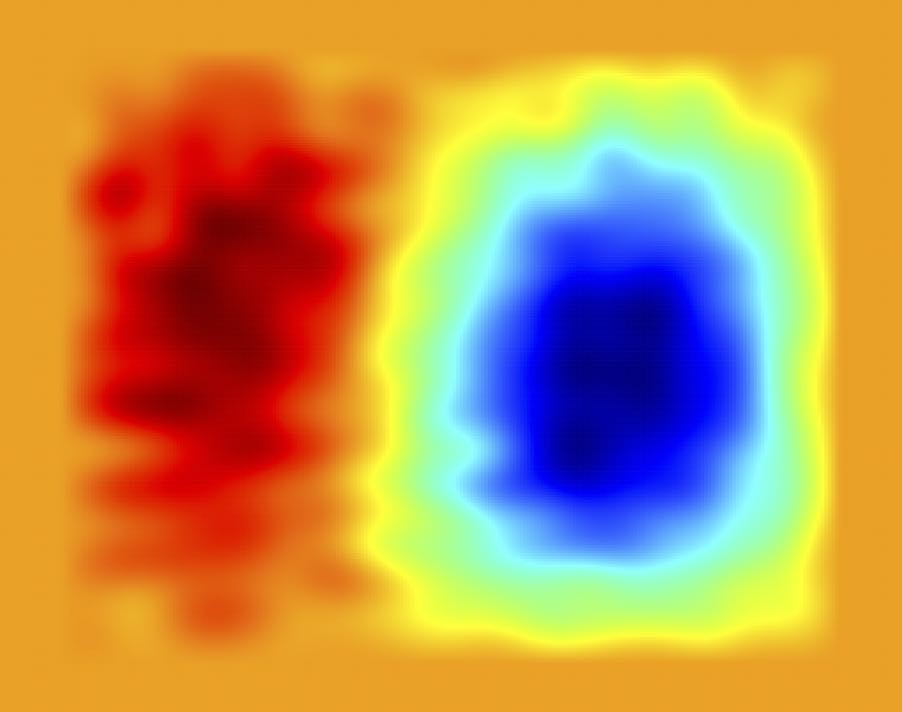}\label{bumps_af_q_a10}
}
\subfigure[$q_\mathrm{SISDP}$ for $50\alpha$]{
\includegraphics[width=0.22\textwidth]{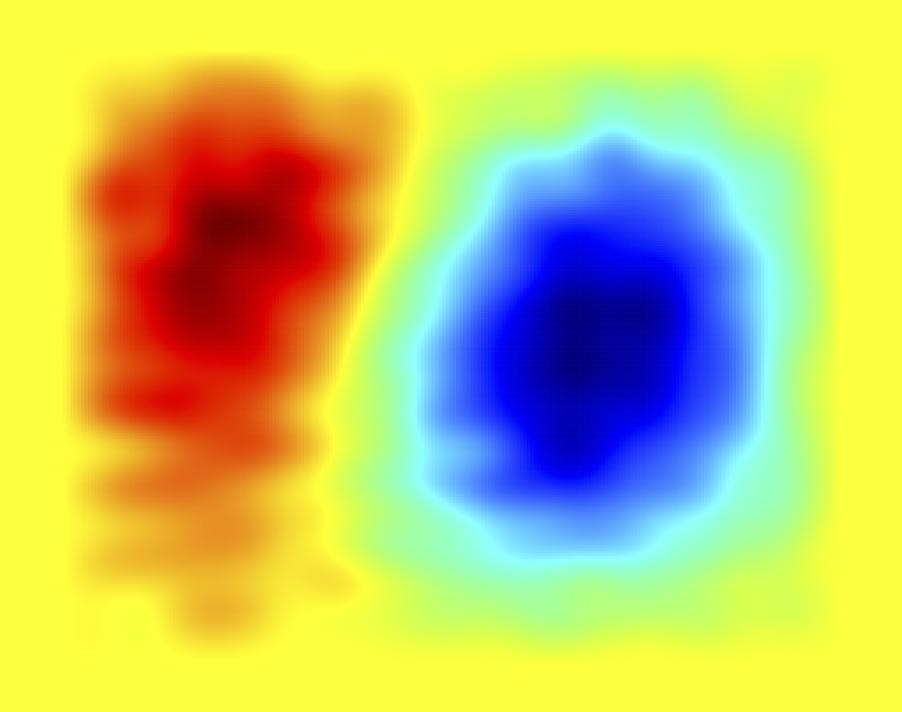}\label{bumps_af_q_a50}
}
\subfigure[$q_\mathrm{SISDP}$ for $100\alpha$]{
\includegraphics[width=0.22\textwidth]{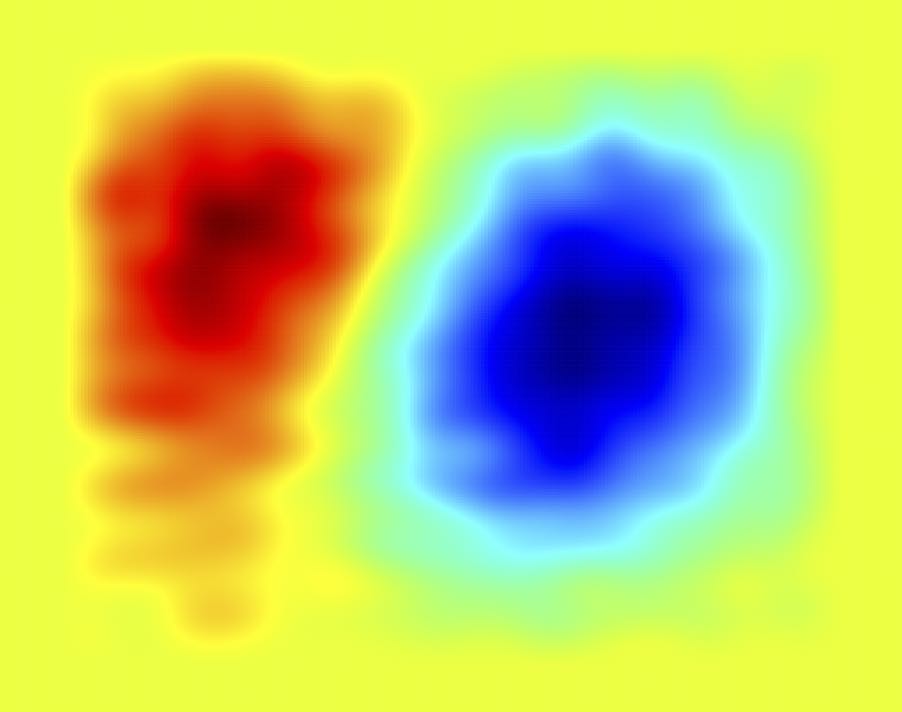}\label{bumps_af_q_a100}
}
\caption{Reconstruction of $\qb$ for Example \ref{example:SDP} part (a). The {\bf SISDP} algorithm is used to reconstruct the domain $\qp$ in the presence of an isotropic background medium generated by the parameter $\barsigma=10$. On the top row, we present the solution $q_\mathrm{SISDP}+\eta_\mathrm{SISDP}$ when the regularization parameter $\alpha$ is multiplied by the constants $1$, $10$, $50$, and $100$. On the bottom row, we present the solution $q_\mathrm{SISDP}$ using the regularization parameter $\alpha$ (initially determined by Algorithm ~\ref{alg:find_alpha} in the Appendix) multiplied by the constants $1$, $10$, $50$, and $100$.}
\label{ex:SISDP:qb:af:results}
\end{figure}

\begin{figure}[!htp]
\centering
\subfigure[$q_\mathrm{SISDP}+\eta_\mathrm{SISDP}$ for $\alpha$]{
\includegraphics[width=0.22\textwidth]{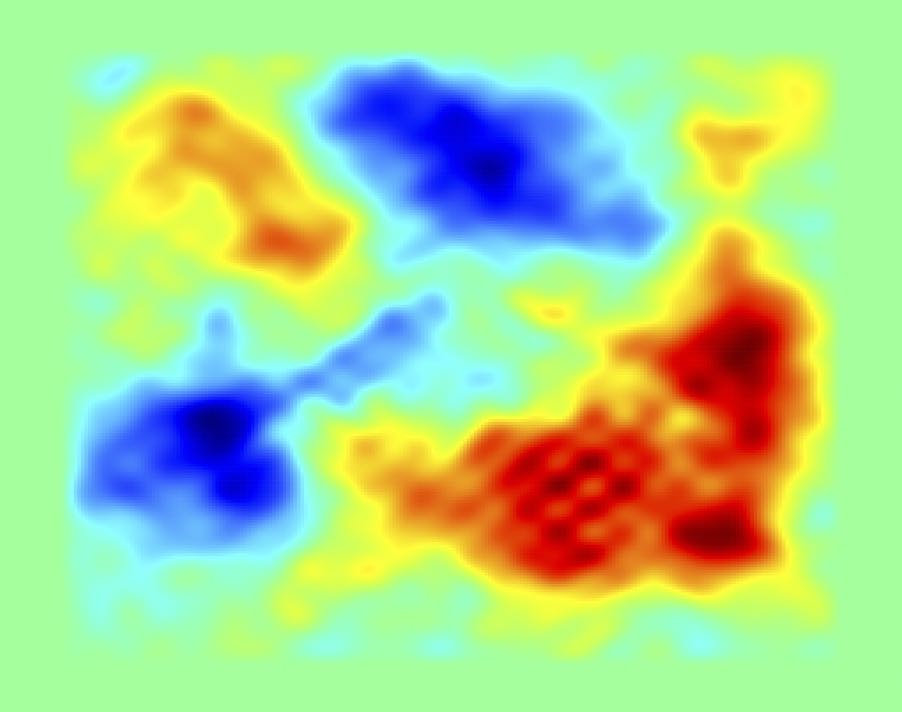}\label{bomber_af_noise_and_q_a1}
}
\subfigure[$q_\mathrm{SISDP}+\eta_\mathrm{SISDP}$ for $10\alpha$]{
\includegraphics[width=0.22\textwidth]{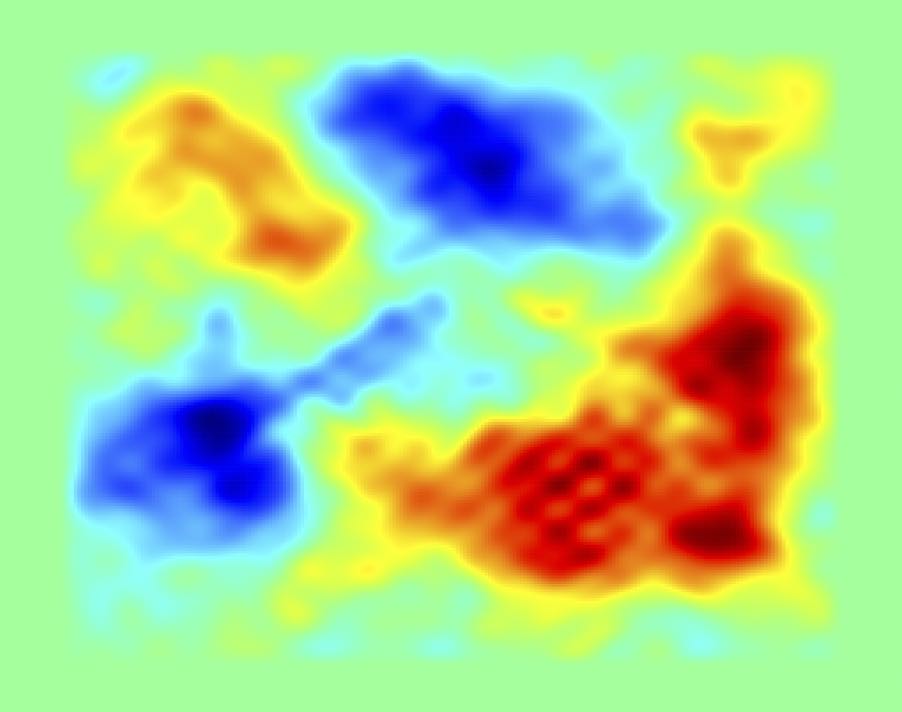}\label{bomber_af_noise_and_q_a10}
}
\subfigure[$q_\mathrm{SISDP}+\eta_\mathrm{SISDP}$ for $50\alpha$]{
\includegraphics[width=0.22\textwidth]{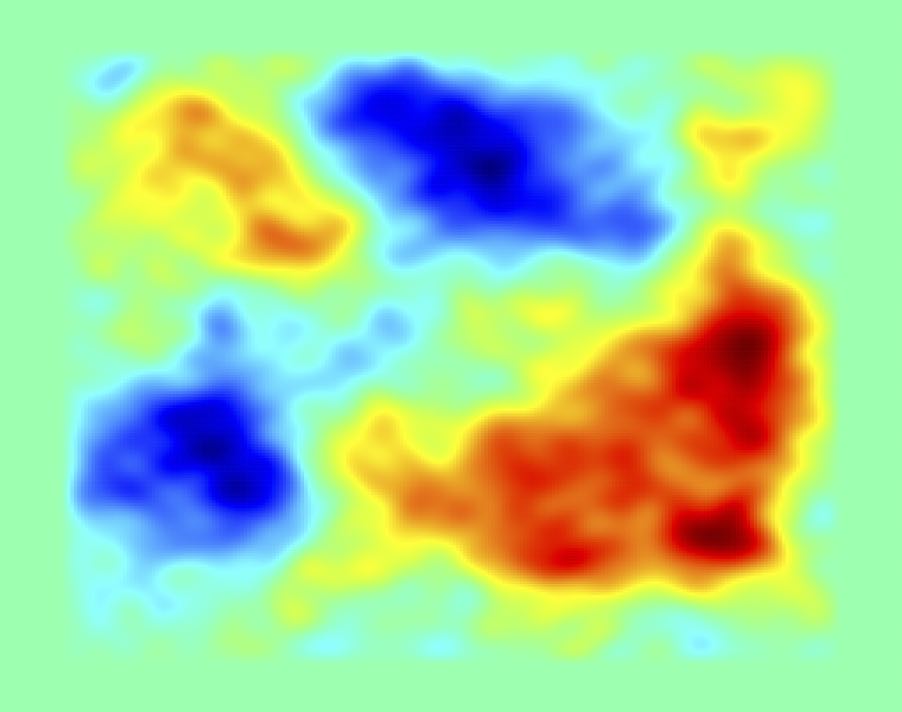}\label{bomber_af_noise_and_q_a50}
}
\subfigure[$q_\mathrm{SISDP}+\eta_\mathrm{SISDP}$ for $100\alpha$]{
\includegraphics[width=0.22\textwidth]{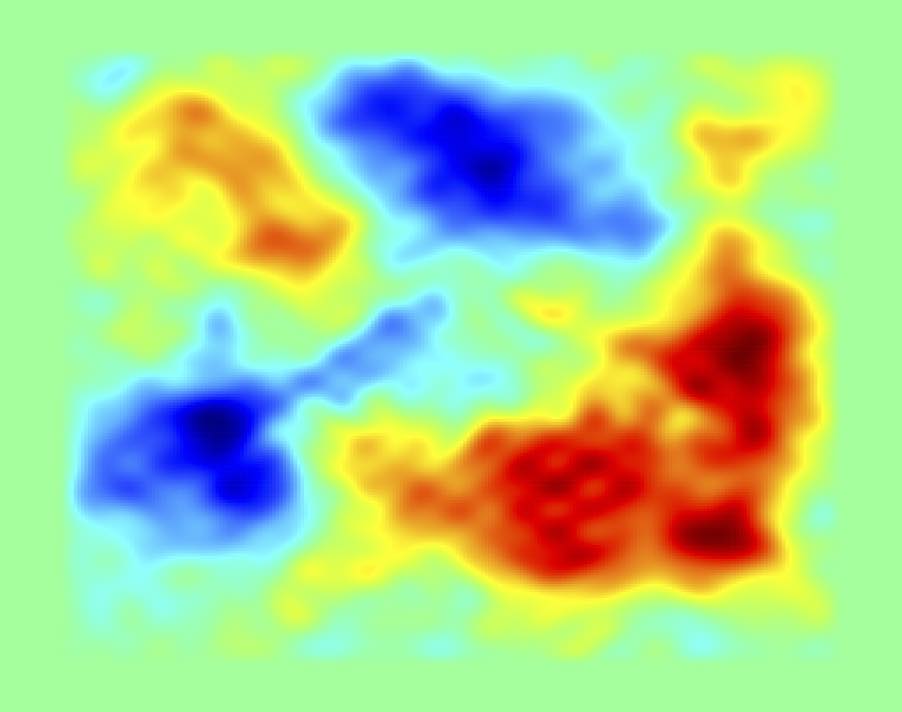}\label{bomber_af_noise_and_q_a100}
}

\subfigure[$q_\mathrm{SISDP}$ for $\alpha$]{
\includegraphics[width=0.22\textwidth]{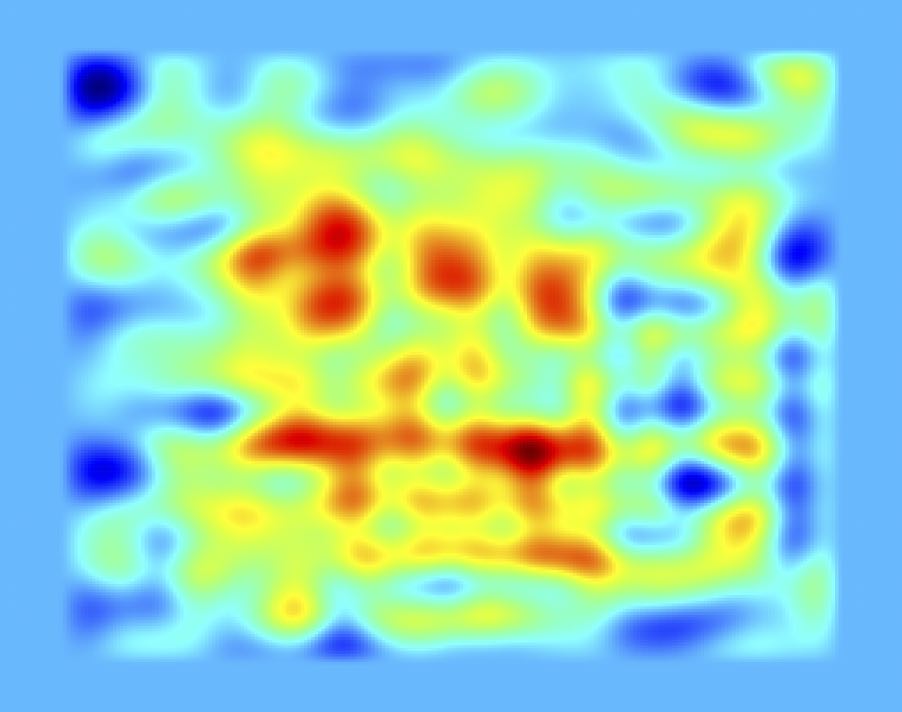}\label{bomber_af_q_a1}
}
\subfigure[$q_\mathrm{SISDP}$ for $10\alpha$]{
\includegraphics[width=0.22\textwidth]{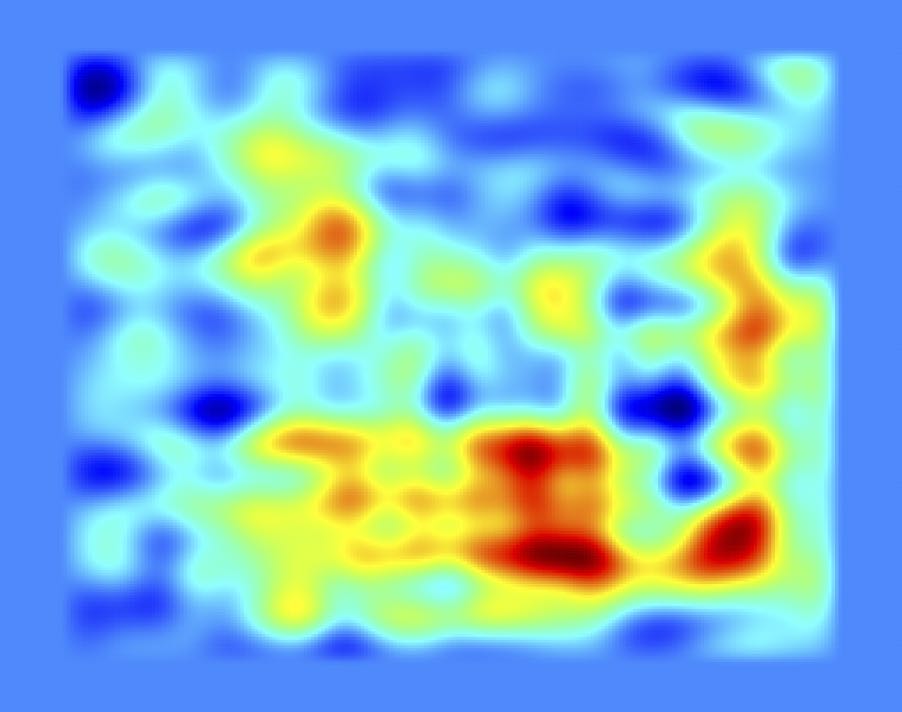}\label{bomber_af_q_a10}
}
\subfigure[$q_\mathrm{SISDP}$ for $50\alpha$]{
\includegraphics[width=0.22\textwidth]{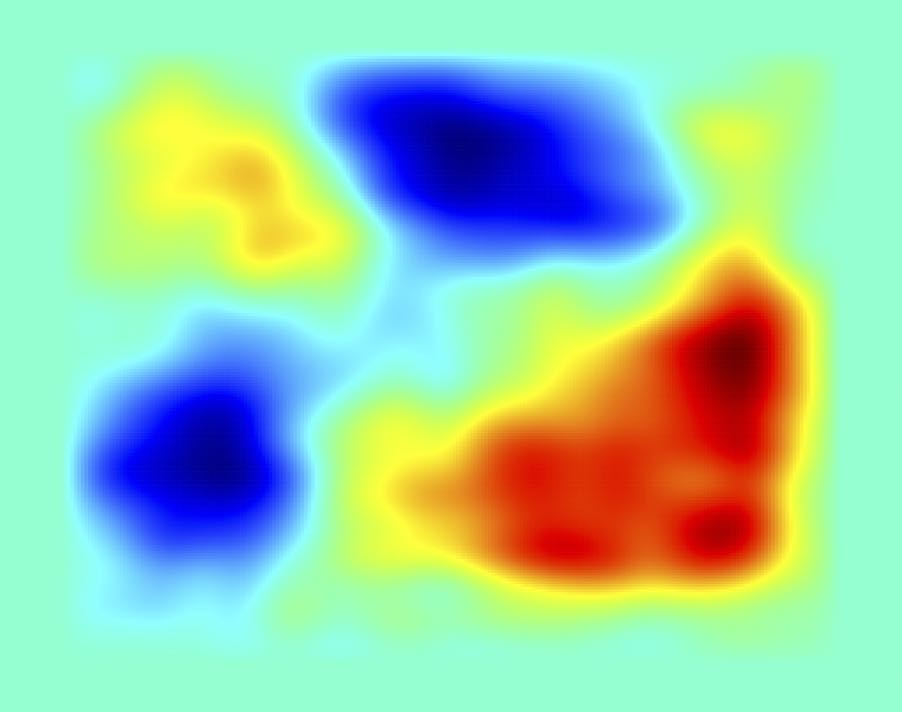}\label{bomber_af_q_a50}
}
\subfigure[$q_\mathrm{SISDP}$ for $100\alpha$]{
\includegraphics[width=0.22\textwidth]{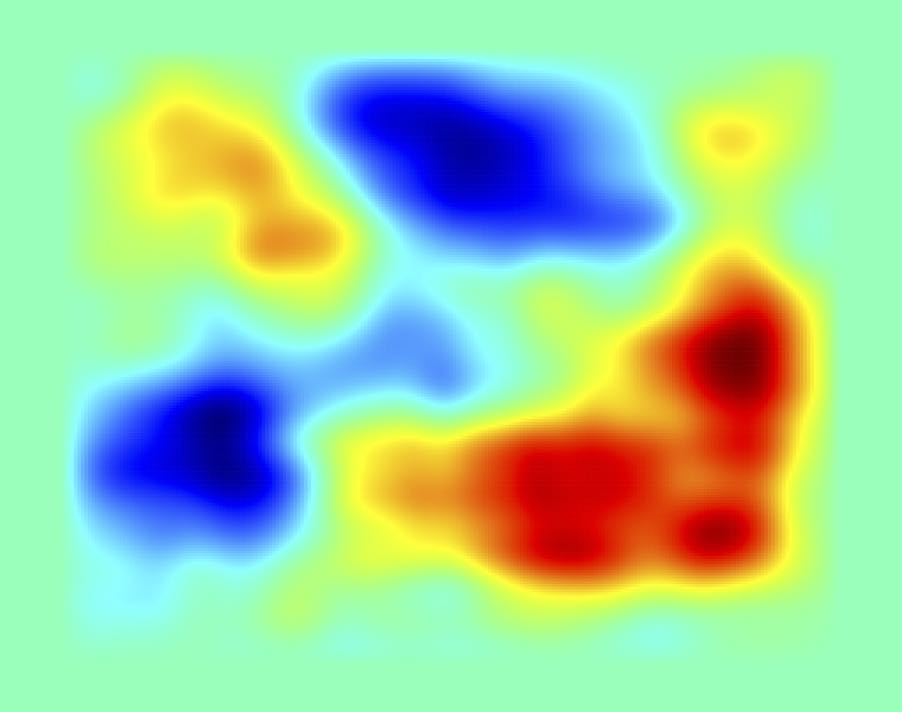}\label{bomber_af_q_a100}
}
\caption{Reconstruction of $\qp$ for Example \ref{example:SDP} part (a). The {\bf SISDP} algorithm is used to reconstruct the domain $\qp$ in the presence of an isotropic background medium generated by the parameter $\barsigma=10$. On the top row, we present the solution $q_\mathrm{SISDP}+\eta_\mathrm{SISDP}$ when the regularization parameter $\alpha$ is multiplied by the constants $1$, $10$, $50$, and $100$. On the bottom row, we present the solution $q_\mathrm{SISDP}$ using the regularization parameter $\alpha$ multiplied by the constants $1$, $10$, $50$, and $100$.}
\label{ex:SISDP:qp:af:results}
\end{figure}

b) In this case, we used the {\bf SISDP} and {\bf RLA} algorithms  for the reconstruction of $\qb$ and $\qp$ in the presence of an isotropic background medium composed only of the higher frequencies in the chosen spectrum. 
 
Here $\eta$ is isotropic, with  $\mathcal{T}_\eta$ being the $2\times2$ identity matrix. We set $\eta_{m_1 m_2}=0$ for $m_1+m_2\leq30$, so that $\eta$ has non-zero components only in  high frequencies, which inevitably will pollute the the high-frequency components of the target scatterer $q$. We expect that, at higher frequencies, the information of both the scatterer and the background medium would be mixed, as seen in part (a).

We apply the algorithms {\bf RLA} and {\bf SISDP} in the same fashion as in the previous case (a). For this example, we also reconstruct the domain applying the {\bf RLA} but assuming we know $\eta^\ast$. We denote this solution by $\hat{q}$. That is,
\begin{equation*}
\hat{q}=\min_q\frac{1}{2}\|\db_{\eta^\ast}-\Fb(q+\eta^\ast)\|^2.
\end{equation*}

We ran {\bf SISDP} for $N_s=30$ samples of the background medium with $\barsigma=10$ and $100$. Figures \ref{fig_boxplots_bump_isotropic_hf_pld1} and \ref{fig_boxplots_bomber_isotropic_hf_pld1} show the box plots with the error 
\begin{equation*}
E_\mathrm{SISDP}(k_j)=\frac{\|q_\mathrm{SISDP}-\hat{q}\|}{\|\hat{q}\|}
\end{equation*}
at each wavenumber for the reconstructions of $\qb$ and $\qp$, respectively. As we can see in those plots, the error decreases as we start to recover data for the background medium $\eta$.

\begin{figure}[h!]
\centering
\subfigure[$\barsigma=10$ and $N_s=30$]{
\input{boxplot_bumps_ihf_10.tex}
}
\subfigure[$\barsigma=100$ and $N_s=30$]{
\input{boxplot_bumps_ihf_100.tex}
}
\caption{Box plots of the error $\|q_\mathrm{SISDP}-\hat{q}\|/\|\hat{q}\|$in the reconstruction of $\qb$ for Example \ref{example:SDP} part (b). We use 30 different samples for the isotropic background medium with noise levels: a) $\barsigma=10$, and b) $\barsigma=100$. The background noise is composed only of high frequencies.}
\label{fig_boxplots_bump_isotropic_hf_pld1}
\end{figure}

\begin{figure}[h!]
\centering
\subfigure[$\barsigma=10$ and $N_s=30$]{
\input{boxplot_plane_ihf_10.tex}
}
\subfigure[$\barsigma=100$ and $N_s=30$]{
\input{boxplot_plane_ihf_100.tex}
}
\caption{Box plots of the error $\|q_\mathrm{SISDP}-\hat{q}\|/\|\hat{q}\|$ in the reconstruction of $\qp$ for Example \ref{example:SDP} part (b). We use 30 different samples for the isotropic background medium with noise levels: a) $\barsigma=10$, and b) $\barsigma=100$. The background noise is composed only of high frequencies.}
\label{fig_boxplots_bomber_isotropic_hf_pld1}
\end{figure}

\begin{figure}[h!]
\centering
\subfigure[$q^\ast+\eta^\ast$ for $\barsigma=10$]{
\includegraphics[width=0.21\textwidth]{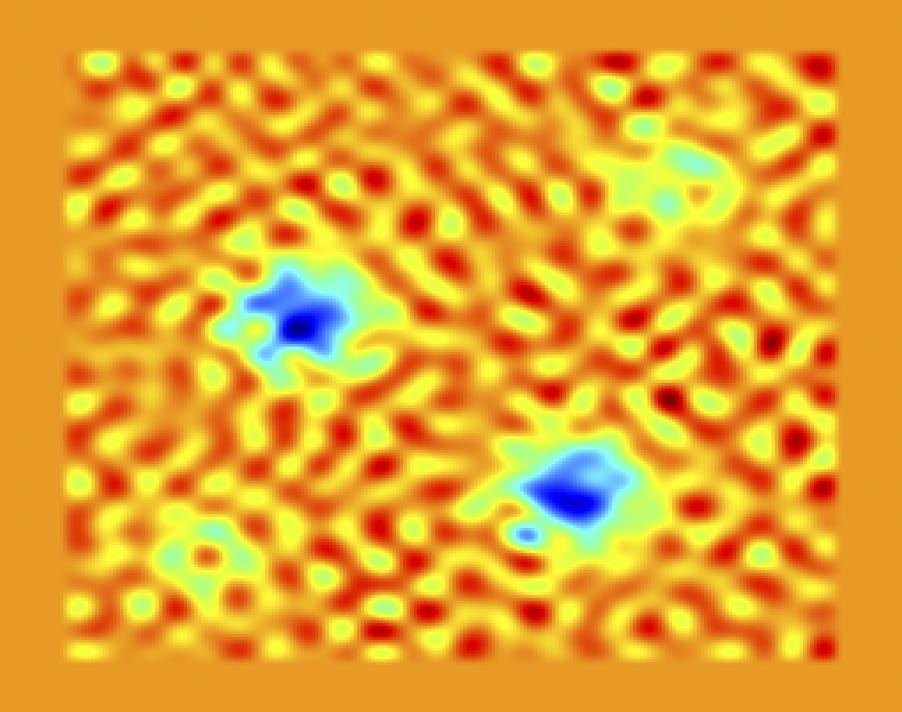}
}
\subfigure[$\hat{q}$ for $\barsigma=10$]{
\includegraphics[width=0.21\textwidth]{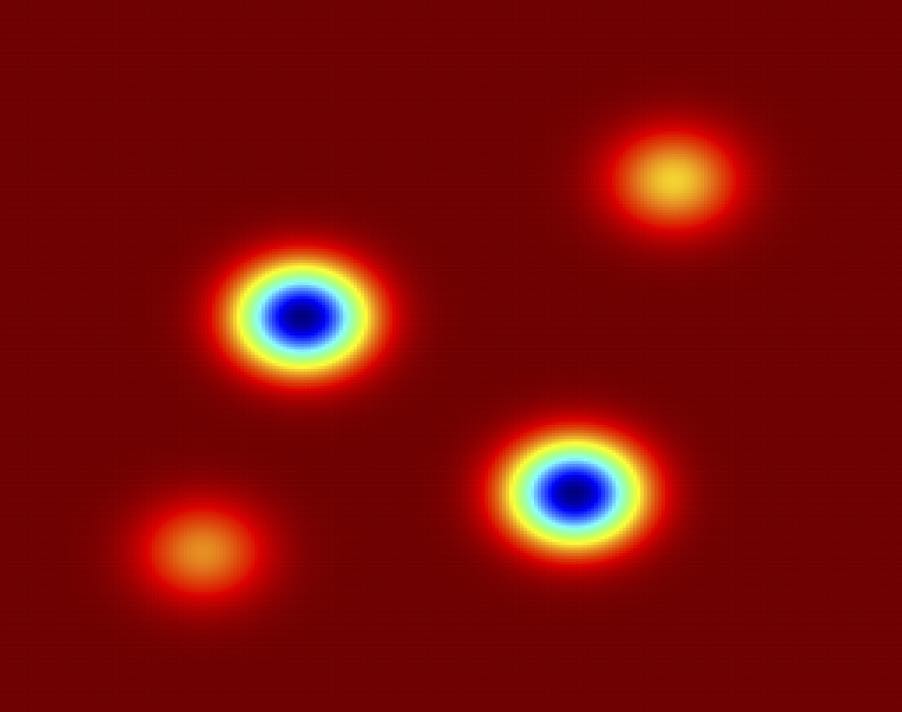}
}
\subfigure[$q_\mathrm{SISDP}$ for $\barsigma=10$]{
\includegraphics[width=0.21\textwidth]{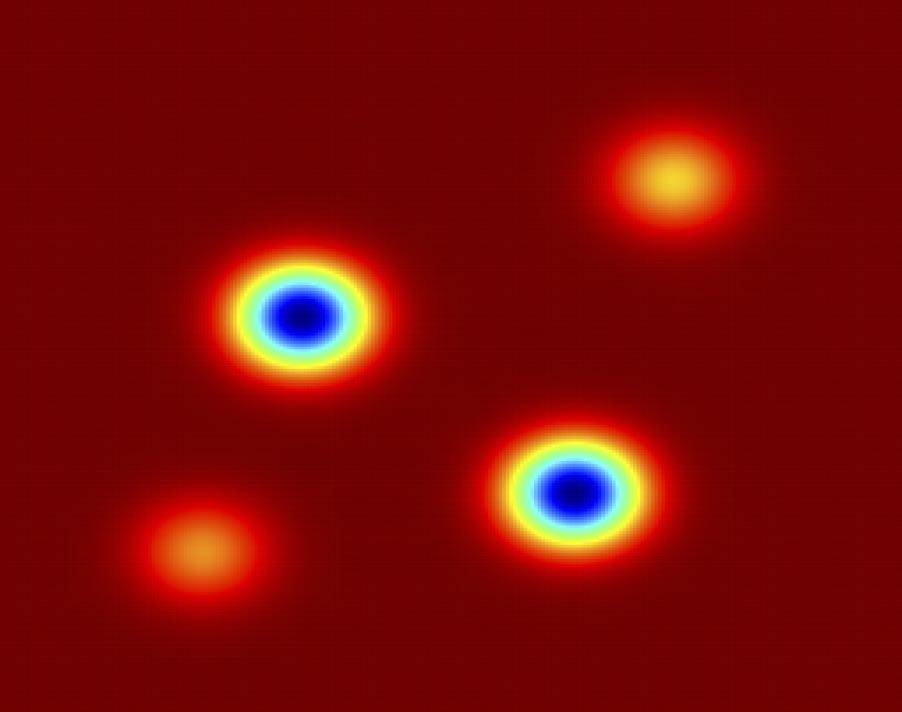}
}
\subfigure[$q_\mathrm{RLA}$ for $\barsigma=10$]{
\includegraphics[width=0.21\textwidth]{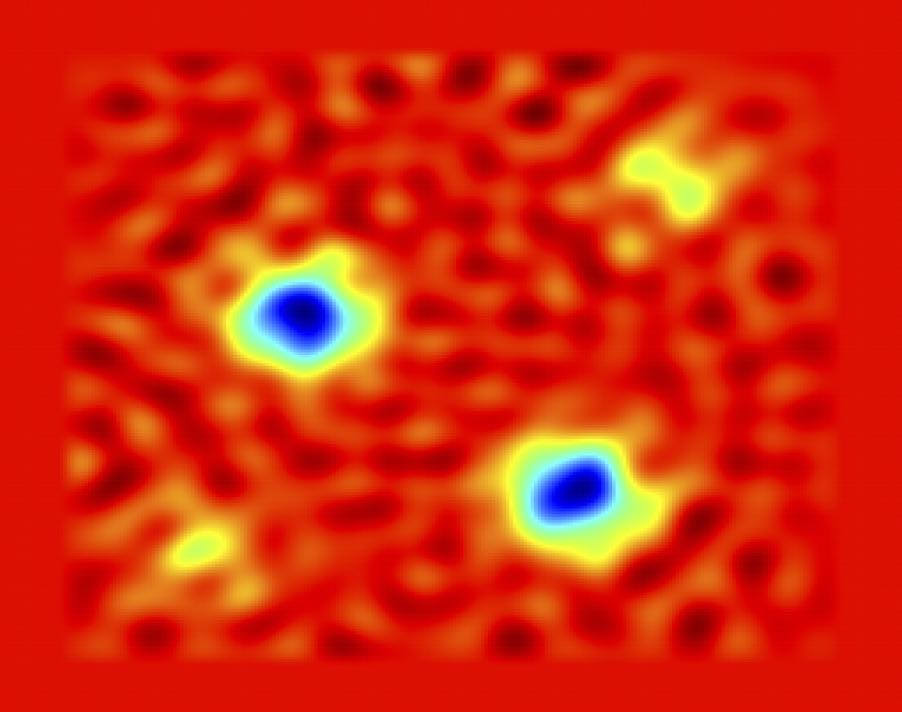}
}

\subfigure[$q^\ast+\eta^\ast$ for $\barsigma=100$]{
\includegraphics[width=0.21\textwidth]{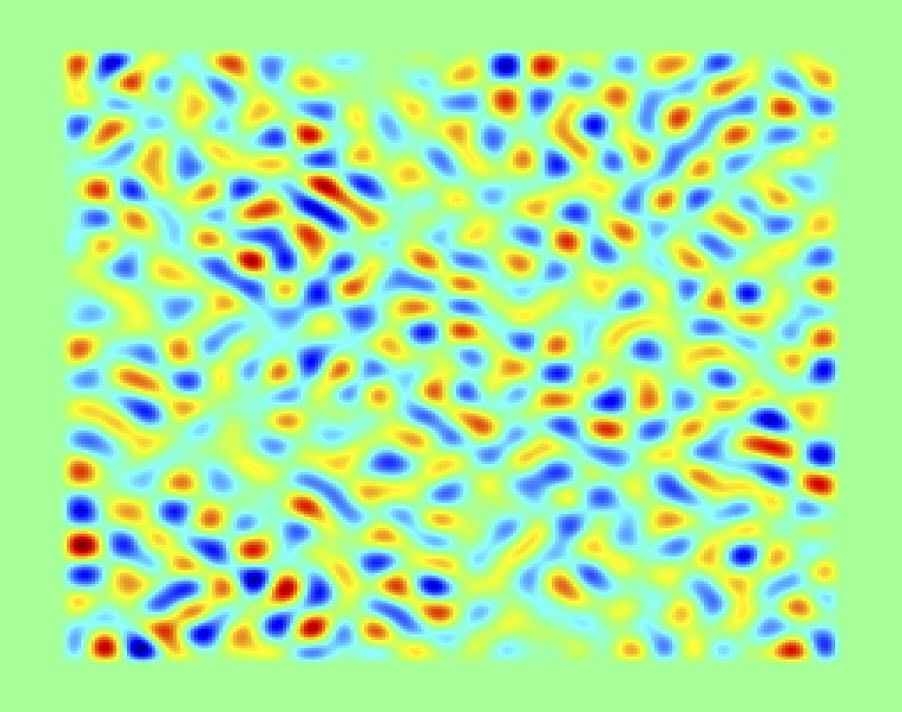}
}
\subfigure[$\hat{q}$ for $\barsigma=100$]{
\includegraphics[width=0.21\textwidth]{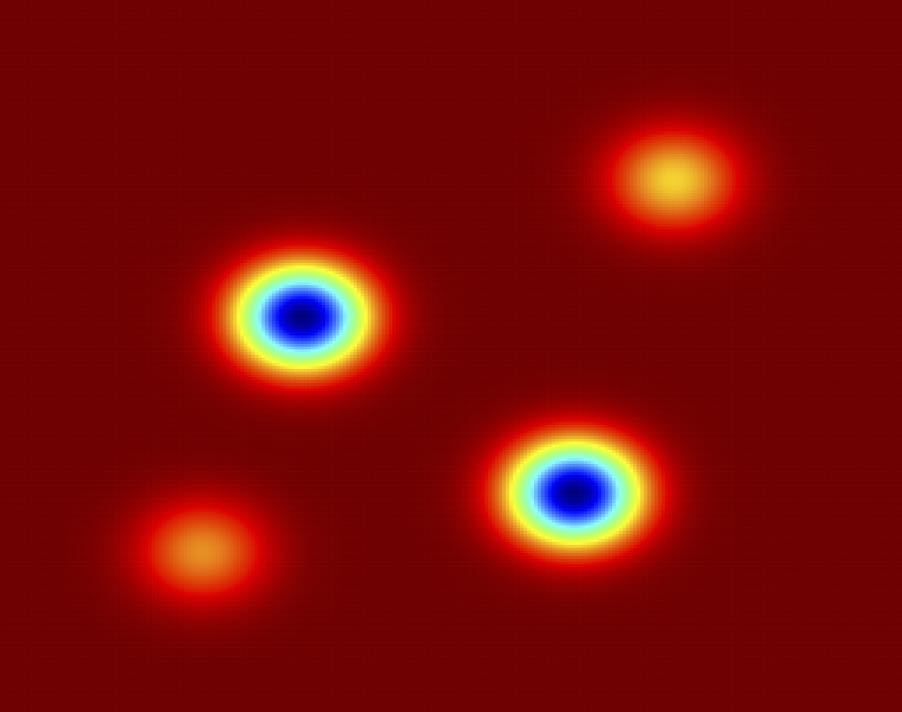}
}
\subfigure[$q_\mathrm{SISDP}$ for $\barsigma=100$]{
\includegraphics[width=0.21\textwidth]{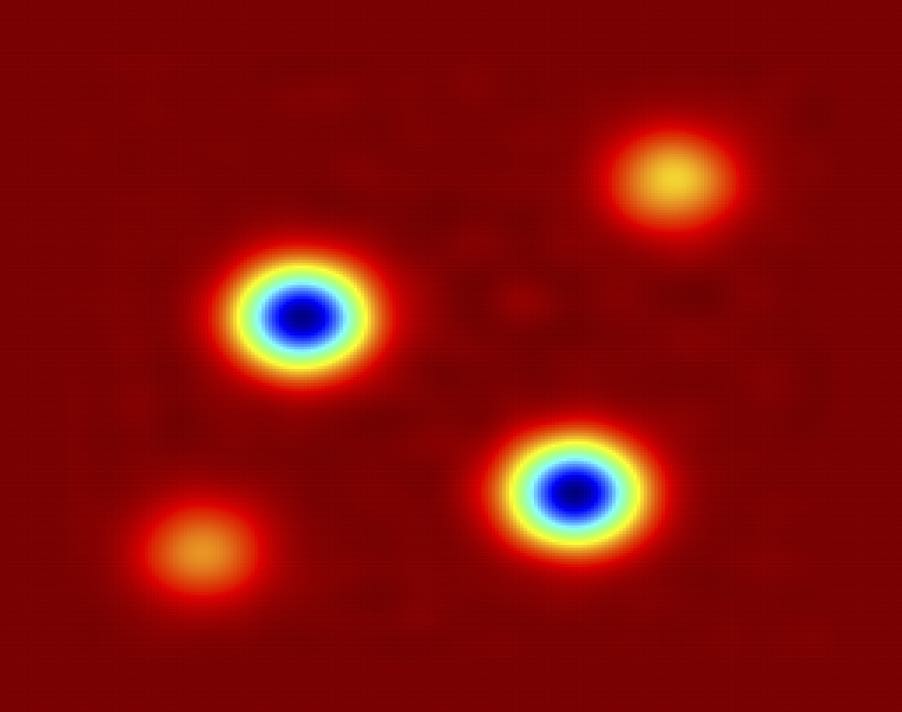}
}
\subfigure[$q_\mathrm{RLA}$ for $\barsigma=100$]{
\includegraphics[width=0.21\textwidth]{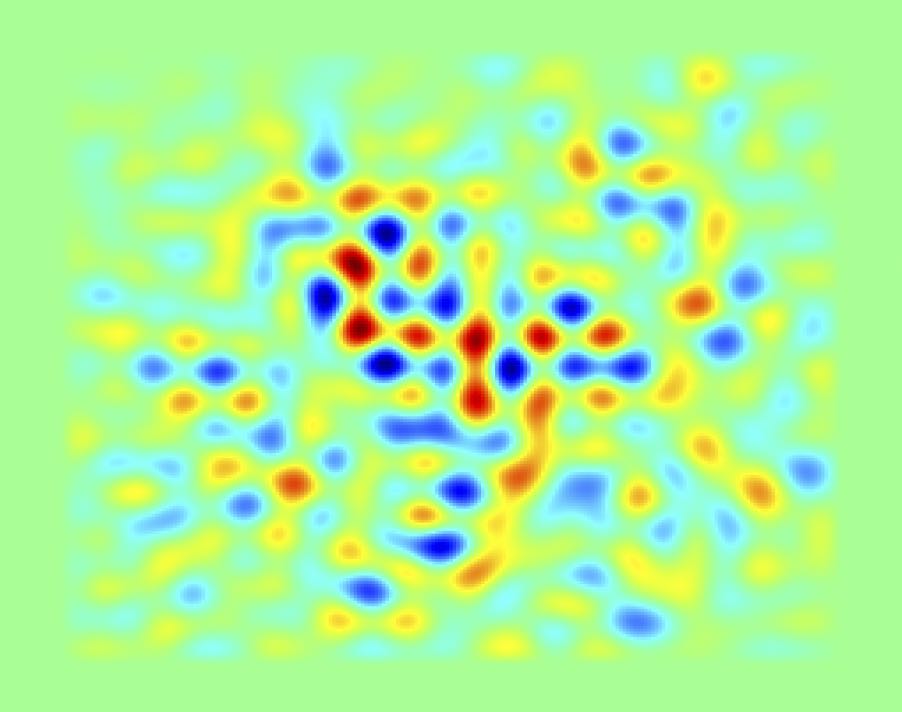}
}

\caption{Reconstruction of $\qb$ for Example \ref{example:SDP} part (b). The {\bf SISDP} algorithm is used to reconstruct the domain $\qb$ in the presence of an isotropic background medium for different noise levels. The background medium has only high frequency components. We present, from top to bottom, in each row, the domain $\qb+\eta^\ast$, the solution $q_\mathrm{SISDP}$ using {\bf SISDP}, and the solution $q_\mathrm{RLA}$ using the standard {\bf RLA} when the background medium is generated using $\barsigma=10$, and $100$.}
\label{fig_bumps_isotropic_hf_pld1}
\end{figure}

\begin{figure}[h!]
\centering
\subfigure[$q^\ast+\eta^\ast$ for $\barsigma=10$]{
\includegraphics[width=0.21\textwidth]{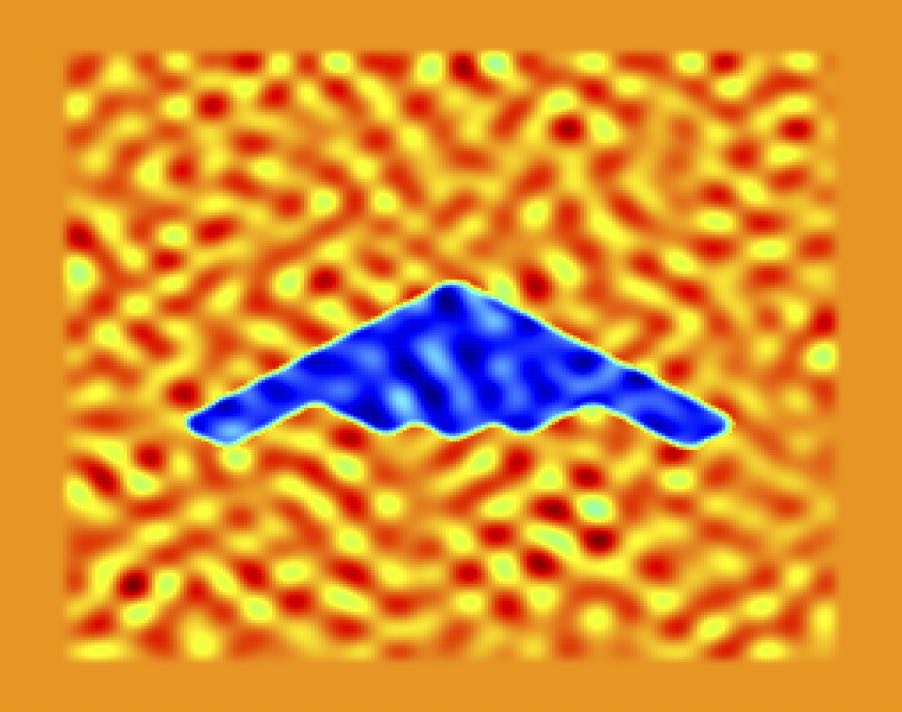}
}
\subfigure[$\hat{q}$ for $\barsigma=10$]{
\includegraphics[width=0.21\textwidth]{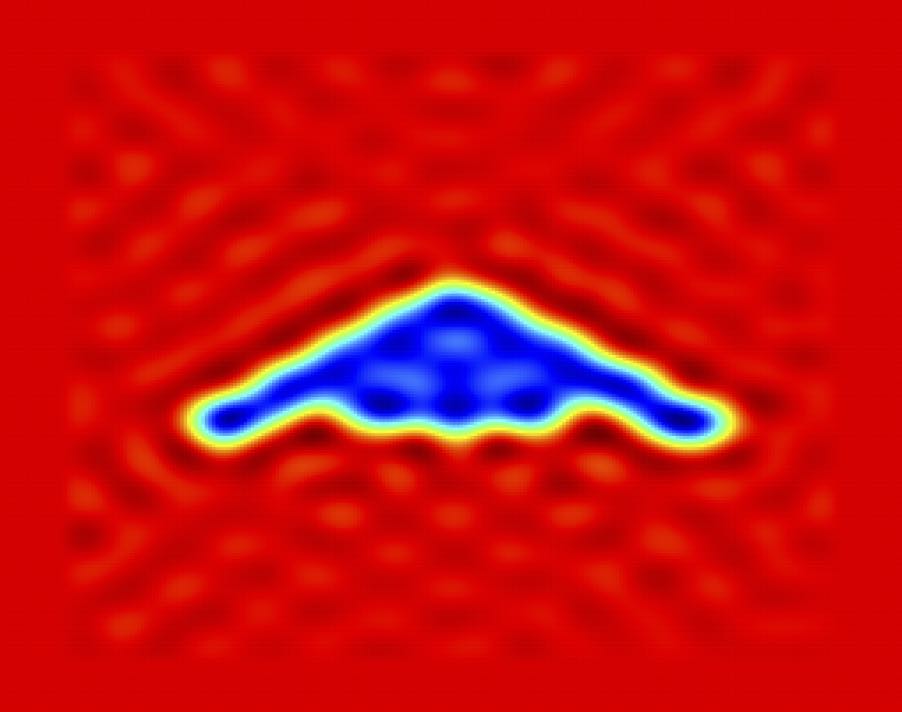}
}
\subfigure[$q_\mathrm{SISDP}$ for $\barsigma=10$]{
\includegraphics[width=0.21\textwidth]{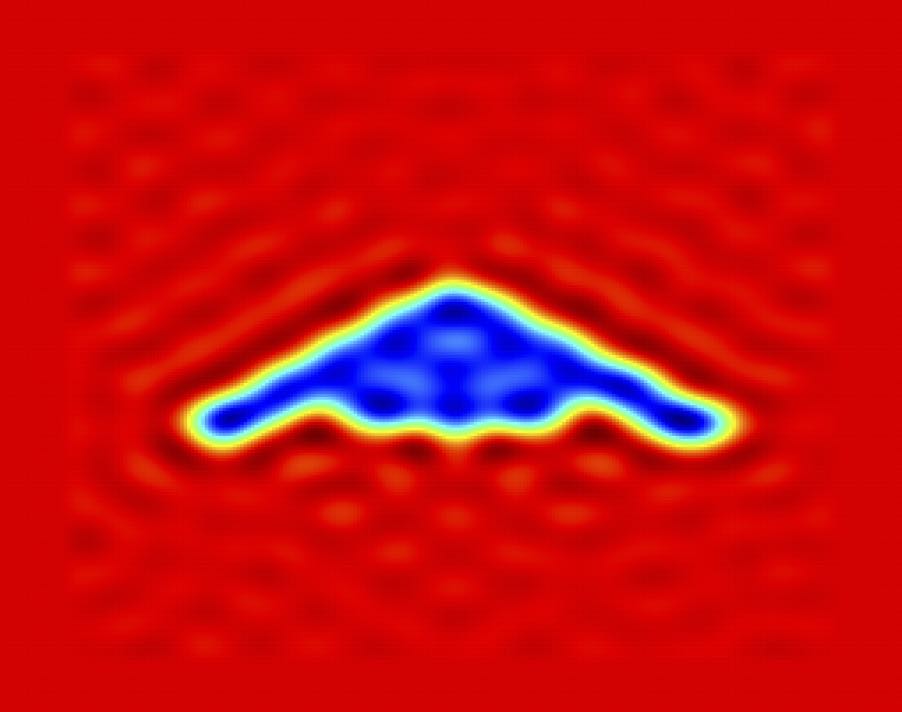}
}
\subfigure[$q_\mathrm{RLA}$ for $\barsigma=10$]{
\includegraphics[width=0.21\textwidth]{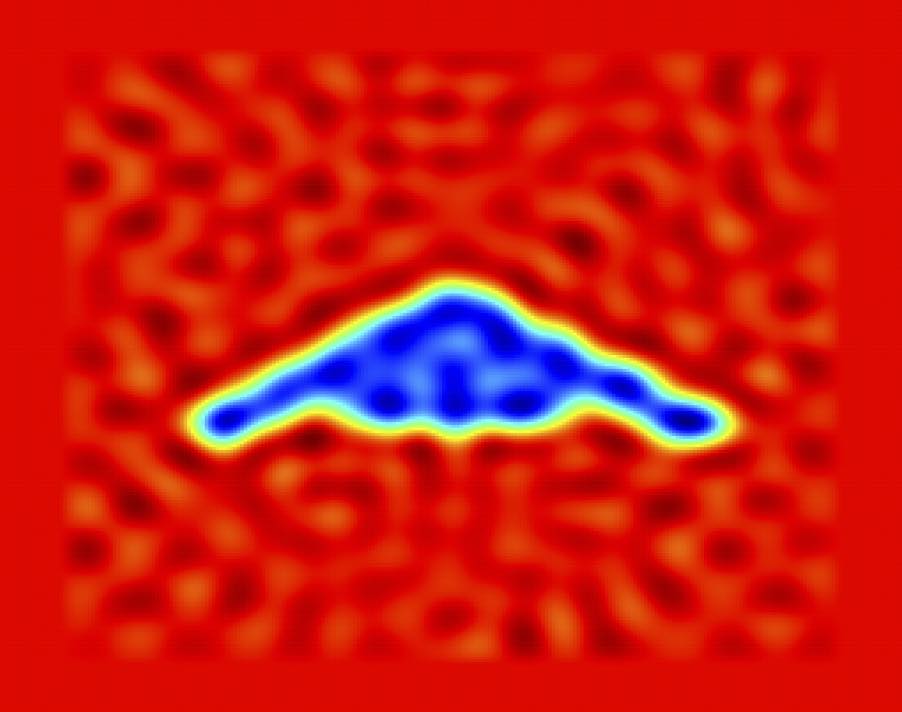}
}

\subfigure[$q^\ast+\eta^\ast$ for $\barsigma=100$]{
\includegraphics[width=0.21\textwidth]{./q_and_noise_plane_iso_hf_d100_s18.jpg}
}
\subfigure[$\hat{q}$ for $\barsigma=100$]{
\includegraphics[width=0.21\textwidth]{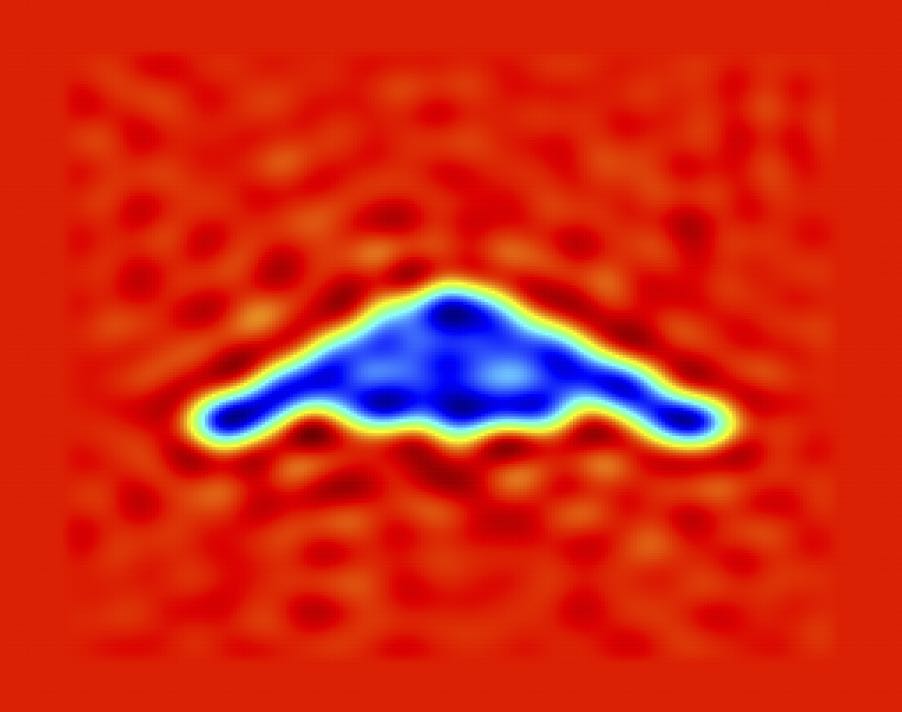}
}
\subfigure[$q_\mathrm{SISDP}$ for $\barsigma=100$]{
\includegraphics[width=0.21\textwidth]{./sol_pld1_plane_iso_hf_d100_s18.jpg}
}
\subfigure[$q_\mathrm{RLA}$ for $\barsigma=100$]{
\includegraphics[width=0.21\textwidth]{./newton_pld1_plane_iso_hf_d100_s18.jpg}
}
\caption{Reconstruction of $\qp$ for Example \ref{example:SDP} part (b). The {\bf SISDP} algorithm is used to reconstruct the domain $\qp$ in the presence of an isotropic background medium for different noise levels. The background medium has only high frequency components. We present, from top to bottom, in each row, the domain $\qp+\eta^\ast$, the solution $q_\mathrm{SISDP}$ using {\bf SISDP}, and the solution $q_\mathrm{RLA}$ using the standard {\bf RLA} when the background medium is generated using $\barsigma=10$, and $100$.}
\label{fig_plane_isotropic_hf_pld1}
\end{figure}

In Figures \ref{fig_bumps_isotropic_hf_pld1} and \ref{fig_plane_isotropic_hf_pld1}, we present the reconstruction of $\qb$ and $\qp$ respectively. We present $q^\ast+\eta^\ast$ for $q^\ast=\qb$ and $\qp$, the domain $\hat{q}$ recovered using the {\bf RLA} with  knowledge of the exact value of $\eta^\ast$, the reconstruction $q_\mathrm{SISDP}$ by the {\bf SISDP} algorithm and the reconstruction using only the {\bf RLA} for the noise levels $\barsigma=10$ and $100$.

The reconstruction using {\bf SISDP} has a quality similar to that of the reconstruction using {\bf RLA} with the knowledge of the background medium; however, it is much better than the standard {\bf RLA} without the knowledge of the background medium.

c) In this case, we used the {\bf SISDP} and {\bf RLA} algorithms for the reconstruction of $q_\mathrm{sub}$ in the presence of an anisotropic background medium.

Assuming that we have prior knowledge of the probability distribution of both $q_\mathrm{sub}$ and the background medium $\eta$, we generate four different functions for the background medium in this example using the anisotropic prior
\begin{equation}
\mathcal{T}_{\eta}=\begin{bmatrix} 10^{-4} & 0 \\ 0 & 10 \end{bmatrix}, \label{eq:tensor_sub}
\end{equation}
with noise levels of $\barsigma= 2$, $4$, $8$ and $16$. We assume that the prior for the probability distribution of the domain $q_\mathrm{sub}$ is
\begin{equation}
  \mathcal{T}_q=\begin{bmatrix} 10 & 0 \\ 0 & 10^{-4} \end{bmatrix}.
\label{eq:tensor_sub2}
\end{equation}

\begin{figure}[h!]
\centering
\subfigure[$q_\mathrm{sub}+\eta^\ast$ for $\barsigma=2$]{
\includegraphics[width=0.28\textwidth]{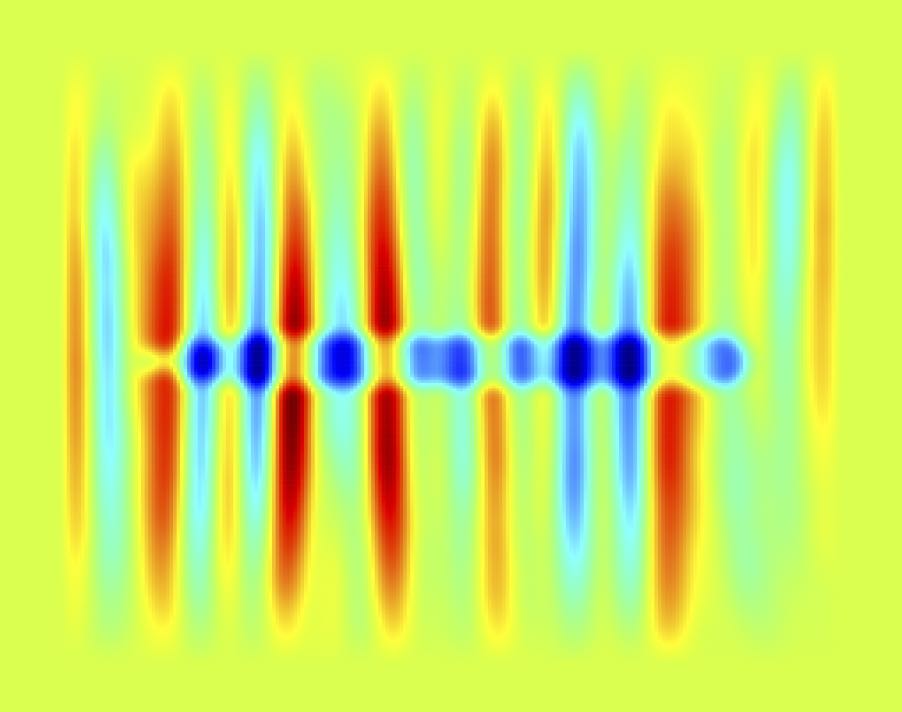}
}
\subfigure[$q_\mathrm{SISDP}$ for $\barsigma=2$]{
\includegraphics[width=0.28\textwidth]{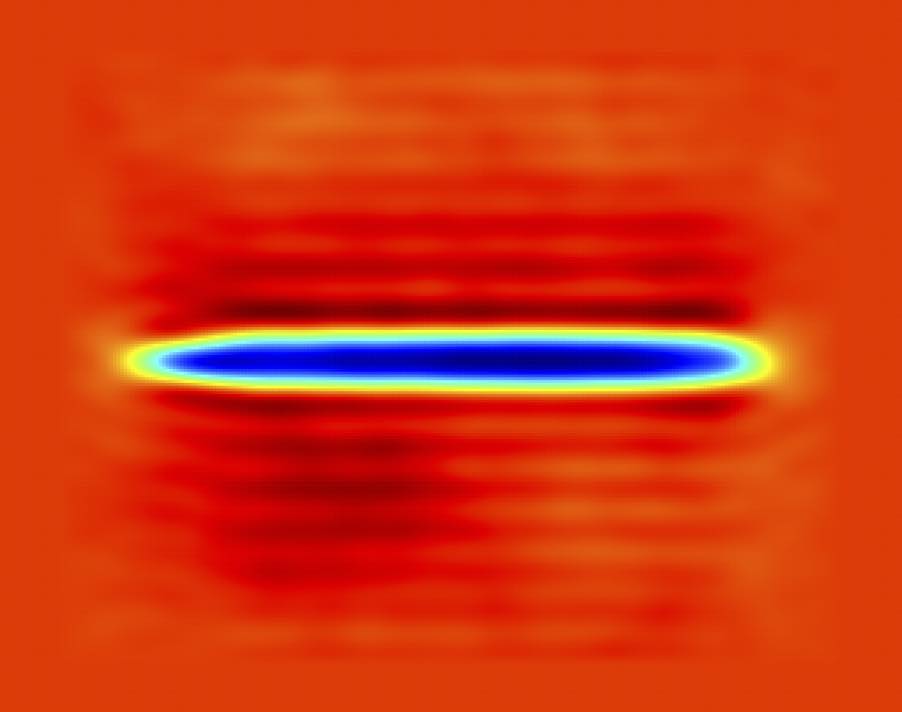}
}
\subfigure[$q_\mathrm{RLA}$ for $\barsigma=2$]{
\includegraphics[width=0.28\textwidth]{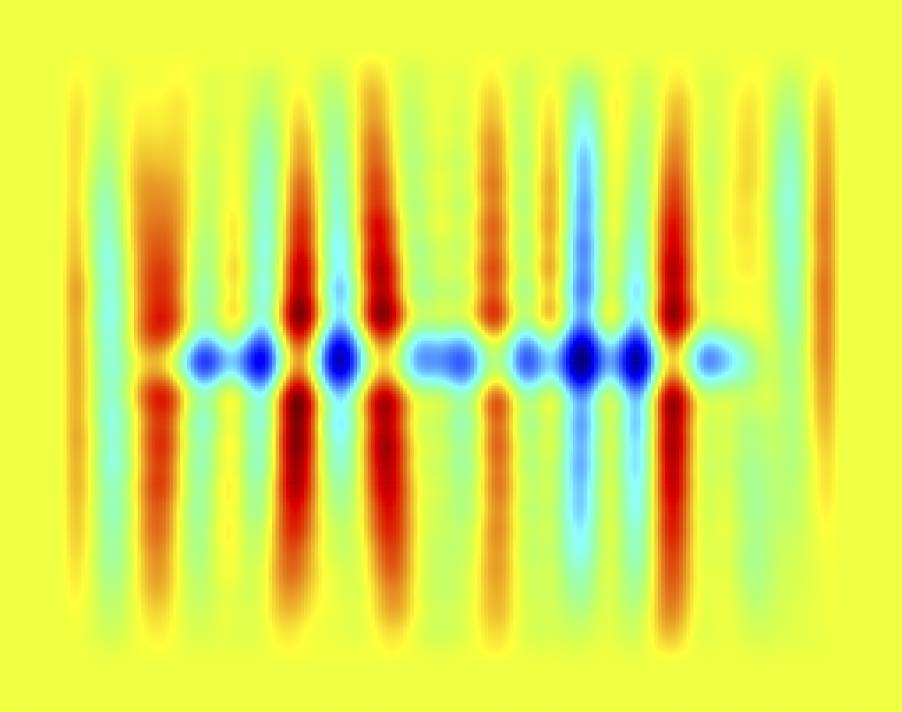}
}

\subfigure[$q_\mathrm{sub}+\eta^\ast$ for $\barsigma=4$]{
\includegraphics[width=0.28\textwidth]{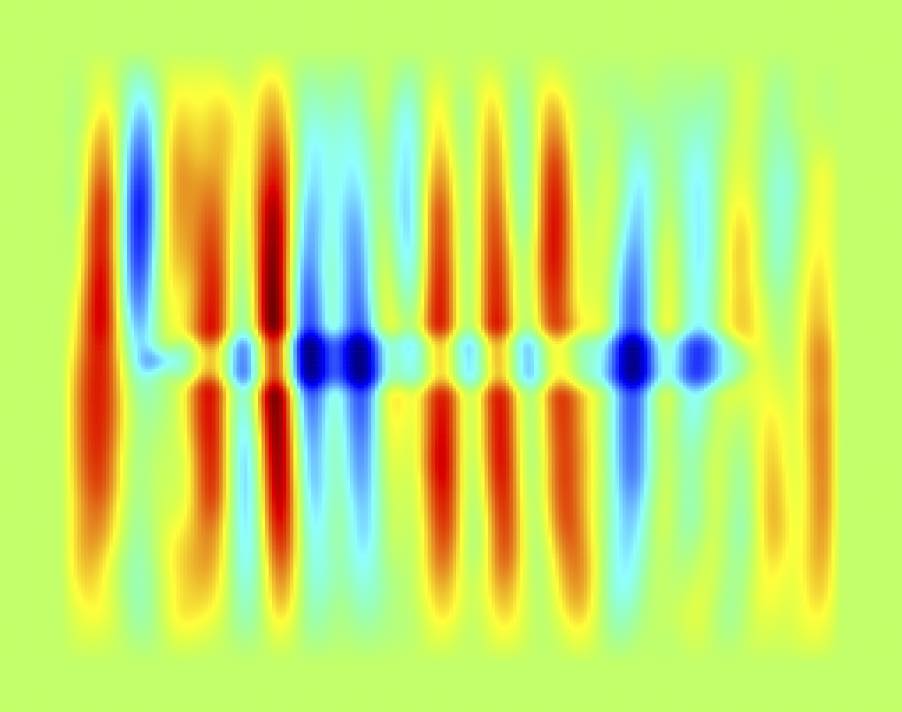}
}
\subfigure[$q_\mathrm{SISDP}$ for $\barsigma=4$]{
\includegraphics[width=0.28\textwidth]{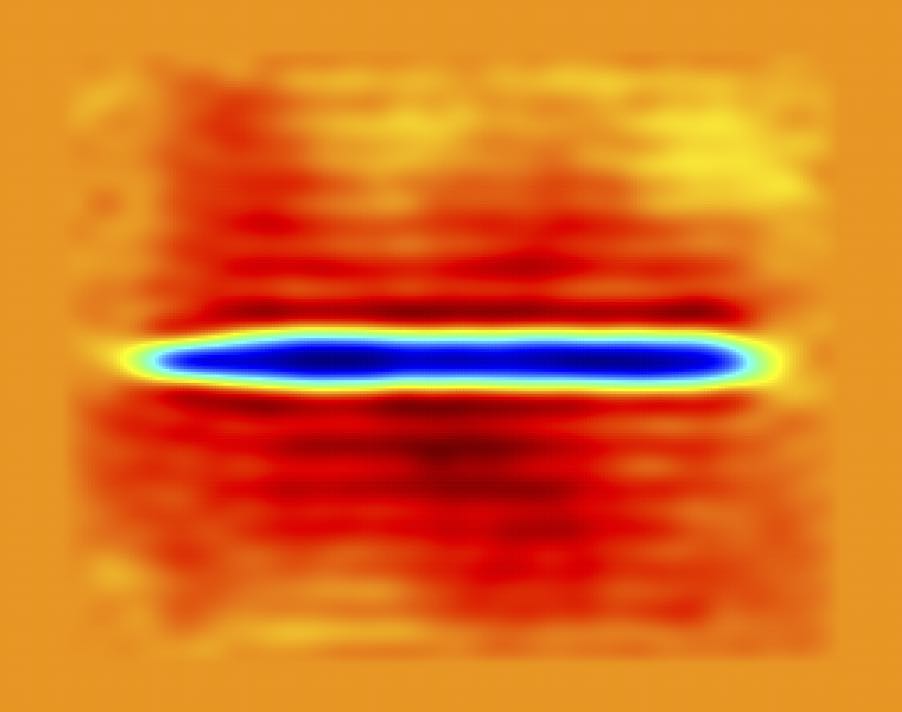}
}
\subfigure[$q_\mathrm{RLA}$ for $\barsigma=4$]{
\includegraphics[width=0.28\textwidth]{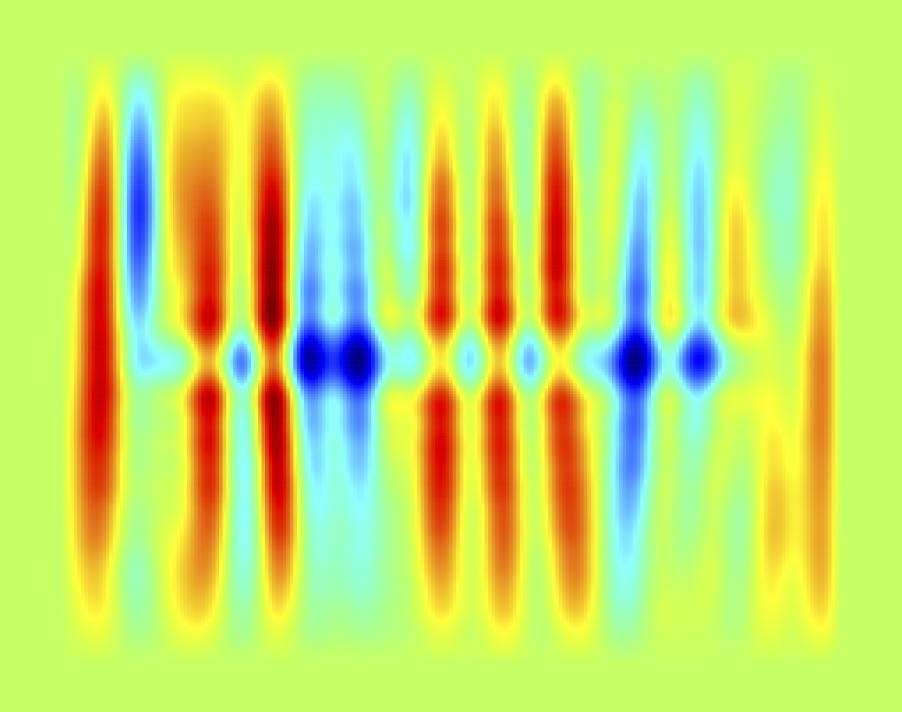}
}

\subfigure[$q_\mathrm{sub}+\eta^\ast$ for $\barsigma=8$]{
\includegraphics[width=0.28\textwidth]{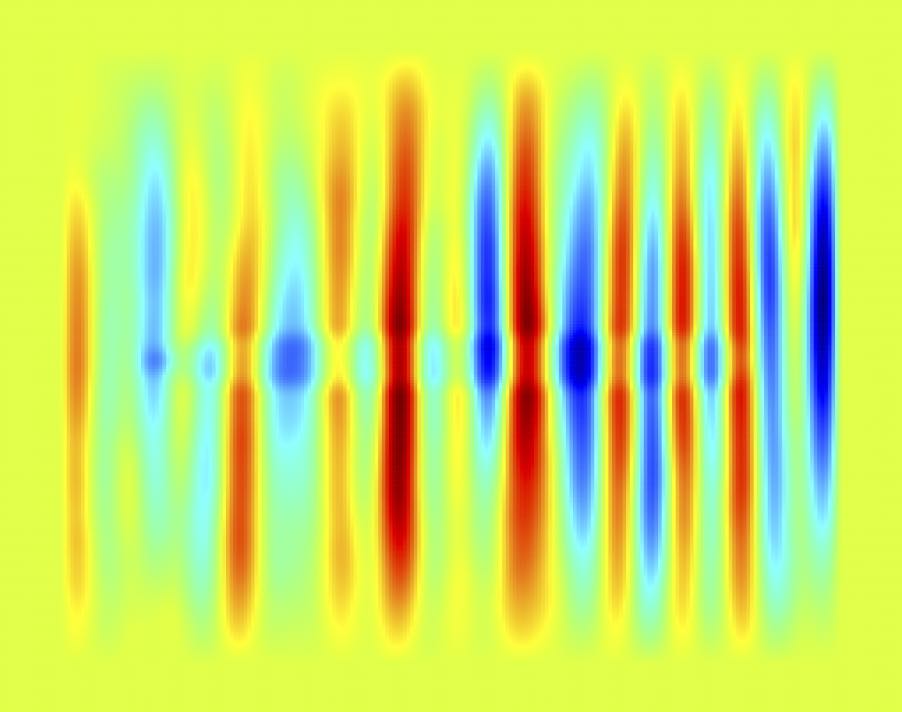}
}
\subfigure[$q_\mathrm{SISDP}$ for $\barsigma=8$]{
\includegraphics[width=0.28\textwidth]{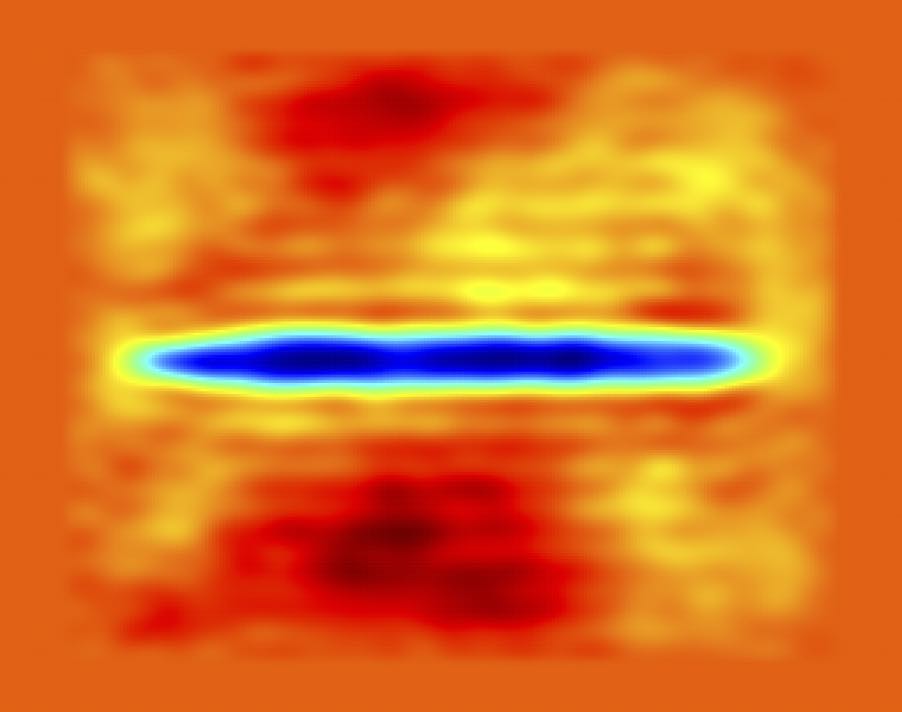}
}
\subfigure[$q_\mathrm{RLA}$ for $\barsigma=8$]{
\includegraphics[width=0.28\textwidth]{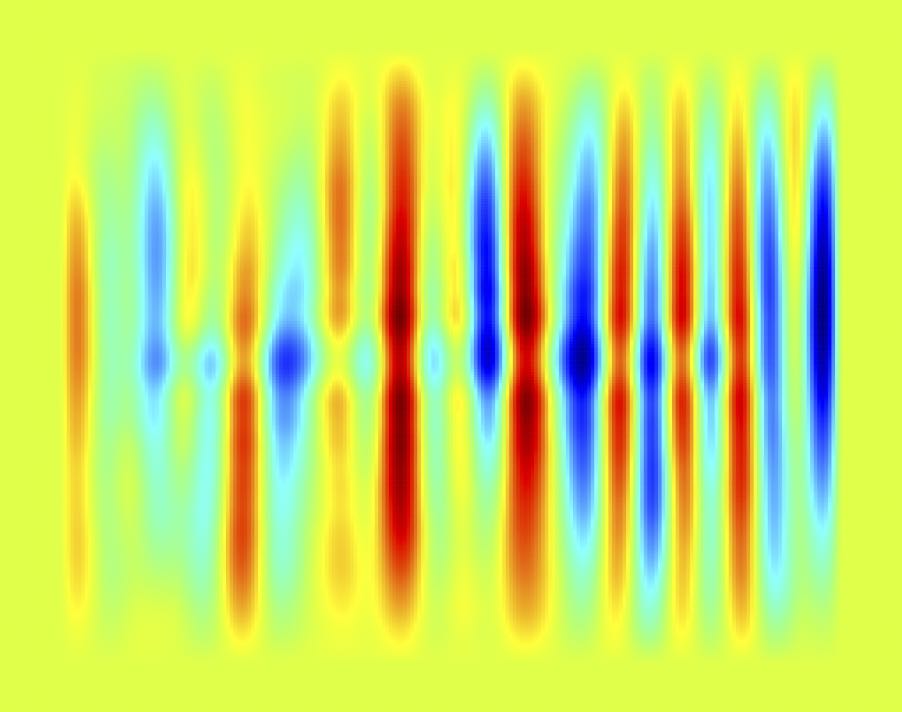}
}

\subfigure[$q_\mathrm{sub}+\eta^\ast$ for $\barsigma=16$]{
\includegraphics[width=0.28\textwidth]{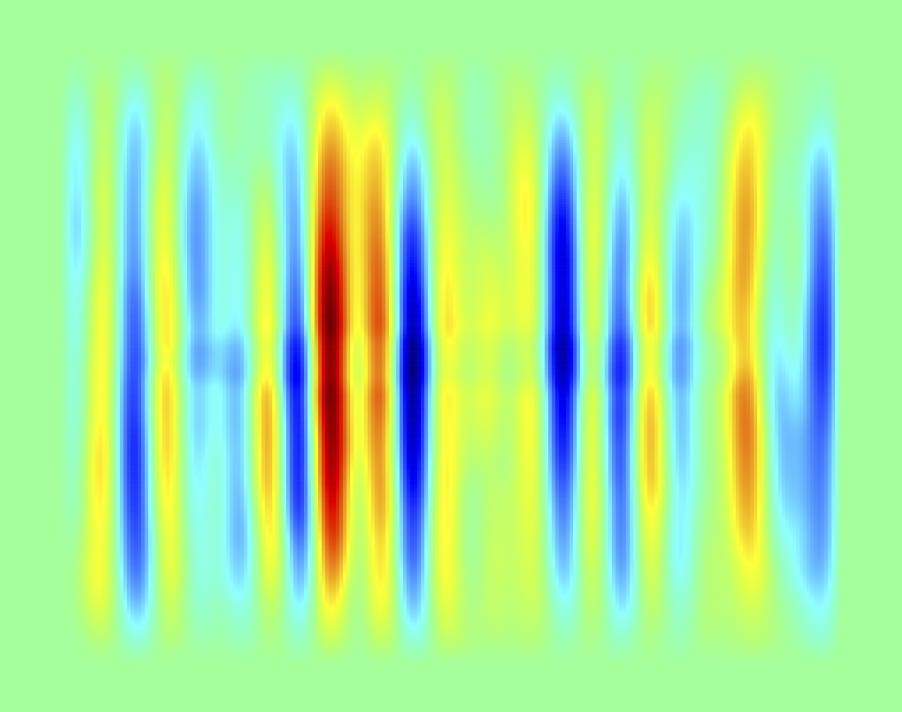}
}
\subfigure[$\eta_\mathrm{SISDP}$ for $\barsigma=16$]{
\includegraphics[width=0.28\textwidth]{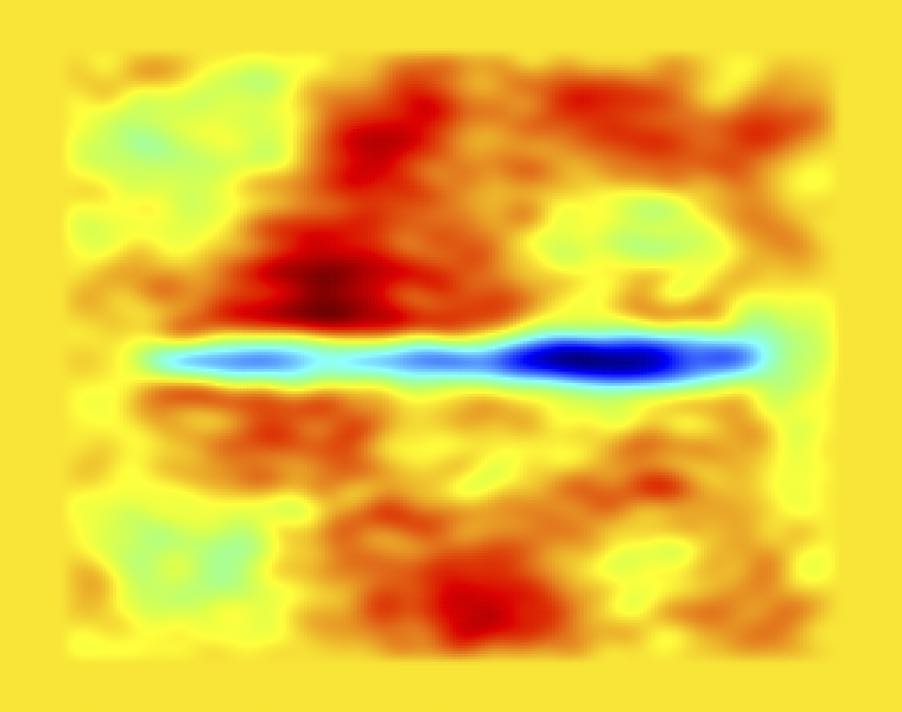}
}
\subfigure[$q_\mathrm{RLA}$ for $\barsigma=16$]{
\includegraphics[width=0.28\textwidth]{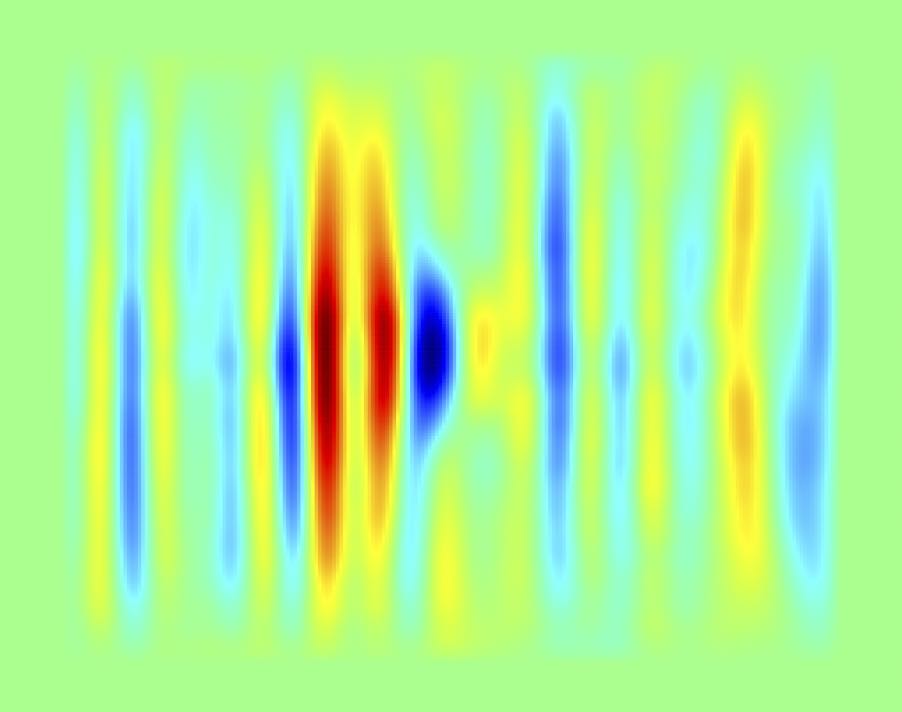}
}
\caption{Reconstruction of $q_\mathrm{sub}$ for Example \ref{example:SDP} part (c). The {\bf SISDP} algorithm is used to reconstruct the domain $q_\mathrm{sub}$ in the presence of an anisotropic background medium generated by the prior \eqref{eq:tensor_sub} with different noise levels. We present, from top to bottom, in each row, the domain $q_\mathrm{sub}+\eta^\ast$, the solution $q_\mathrm{SISDP}$ using {\bf SISDP}, and the solution $q_\mathrm{RLA}$ using the standard {\bf RLA} when the background medium is generated using $\barsigma=2$, $4$, $8$ and $16$.}
\label{fig_anisotropic_af}
\end{figure}

We apply the {\bf RLA} and {\bf SISDP} algorithms using both $p(q)$ and $p(\eta)$. We use as regularization parameters $\alpha=\beta=1$ at all frequencies. The results are presented in Figure \ref{fig_anisotropic_af}. As expected, for each noise level the results of the {\bf SISDP} algorithm are more accurate than the results obtained by the {\bf RLA} with no prior, and we are able to reconstruct the shape of $q_\mathrm{sub}$. As the noise level increases, even though the reconstruction $q_\mathrm{SISDP}$ becomes less accurate, it is still possible to see the shape of the submarine.

\begin{remark} We also used the algorithm SISDP to try to reconstruct $\qb$ in the presence of the background medium $\eta$, when $\eta$ is generated by the distribution \eqref{e:ex1-eta-reg}. As expected, from the results in Example \ref{example:SDP} part a) since $\eta$ has a significant number of components on the energy spectrum of $\qb$ it is not possible to separate $\qb$ from $\eta$.
\end{remark}

\subsection{Data from multiple realizations, inversion for both $q$ and $\eta$}
\label{example:MDP}

\begin{figure}[h!]
\centering
\subfigure[$q_{MIMDP}$ for $N_s=10$]{
\includegraphics[width=0.3\textwidth]{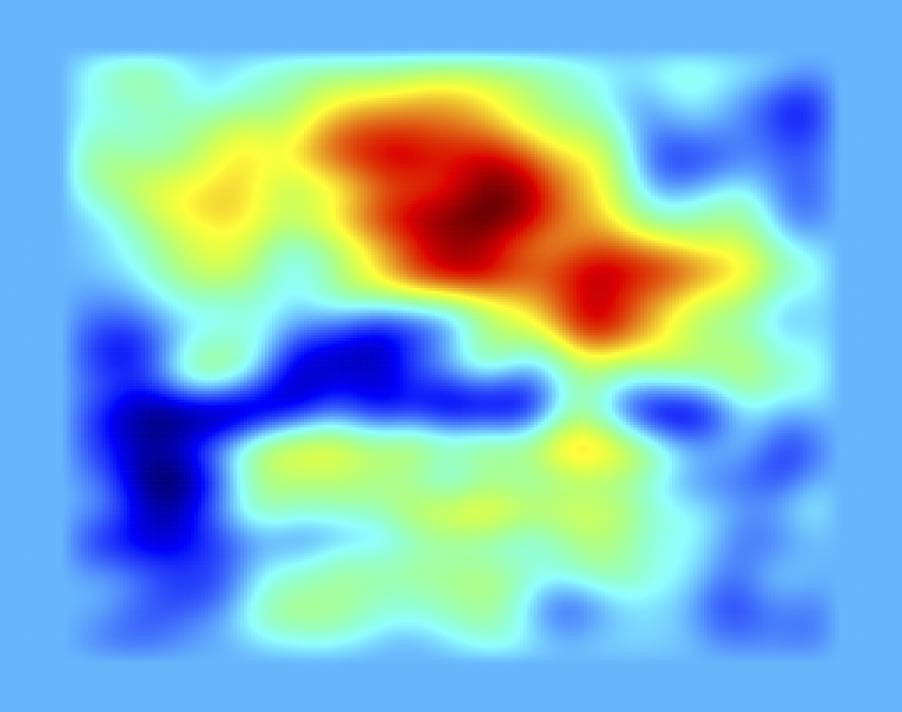}
}
\subfigure[$q_{MIMDP}$ for $N_s=50$]{
\includegraphics[width=0.3\textwidth]{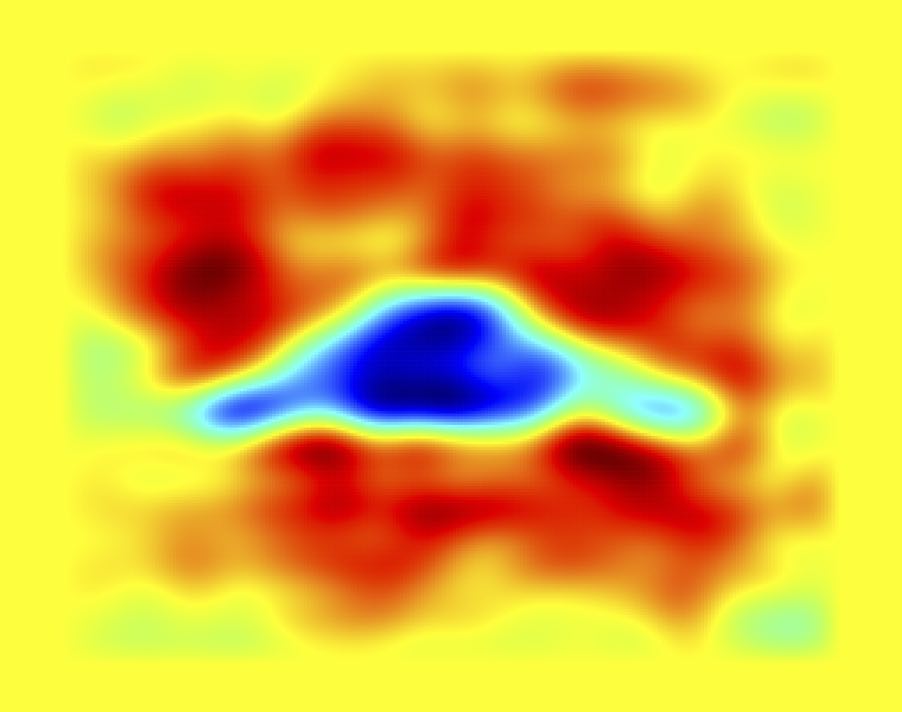}
}
\subfigure[$q_{MIMDP}$ for $N_s=100$]{
\includegraphics[width=0.3\textwidth]{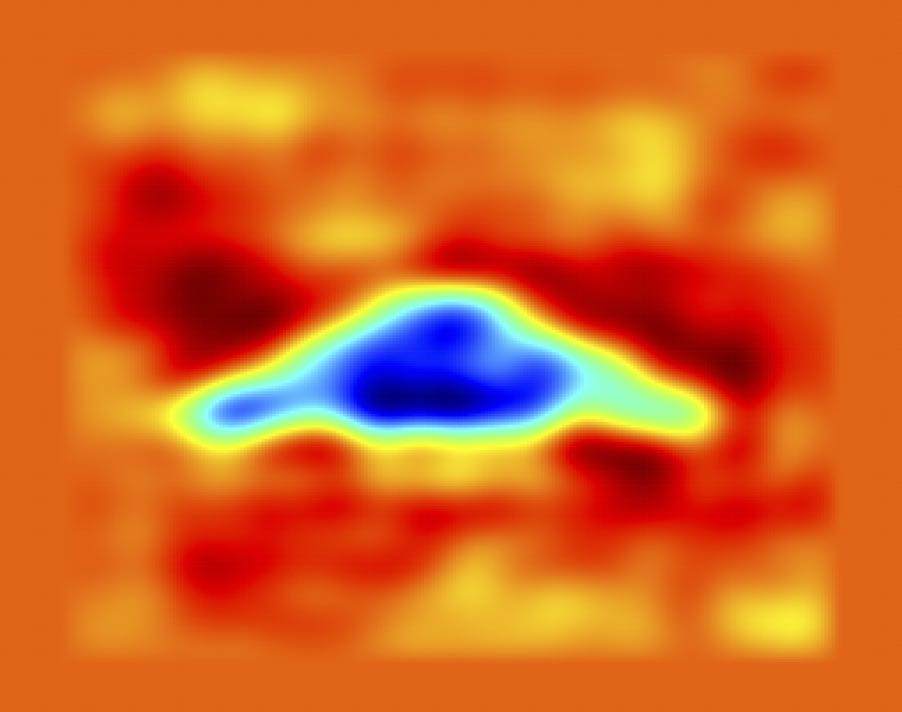}
}
\caption{Reconstruction of $\qp$ for Example \ref{example:MDP}. The {\bf MIMDP} algorithm is used to reconstruct the domain $\qp$ in the presence of an anisotropic background medium that is generated using a noise parameter $\barsigma=10$. The results are presented for the case where we have measurements of the field scattered off of $\qp$ in the presence of: (a)$N_s=10$, (b) $50$ and (c) $100$ samples of the isotropic background domain.}
\label{fig_pid_plane_samples}
\end{figure}

\begin{figure}
\centering
\subfigure[$q_{MIMDP}$ for $N_s=10$]{
\includegraphics[width=0.3\textwidth]{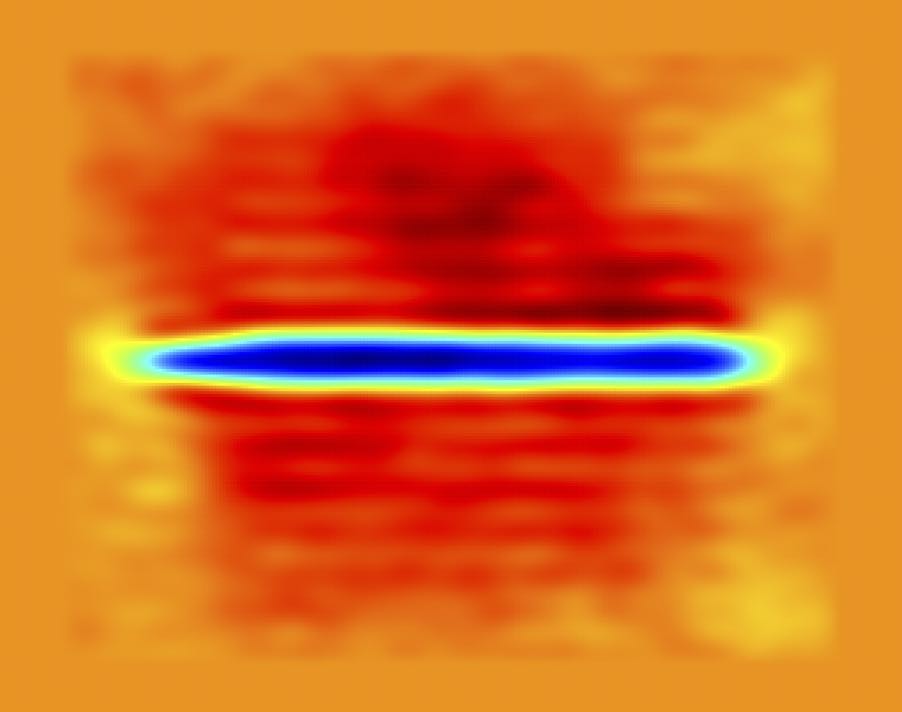}
}
\subfigure[$q_{MIMDP}$ for $N_s=50$]{
\includegraphics[width=0.3\textwidth]{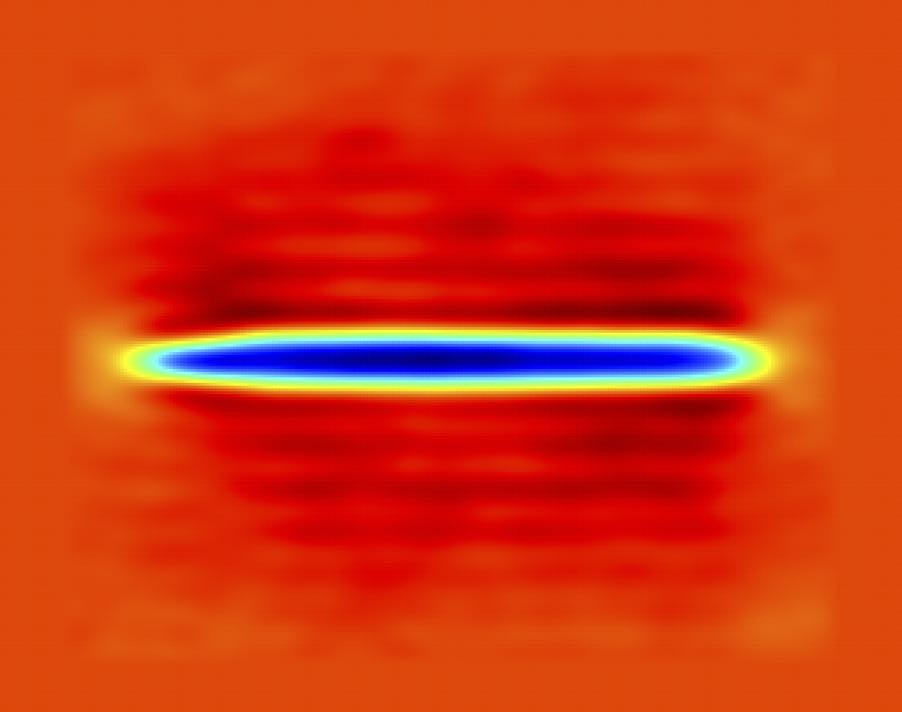}
}
\subfigure[$q_{MIMDP}$ for $N_s=100$]{
\includegraphics[width=0.3\textwidth]{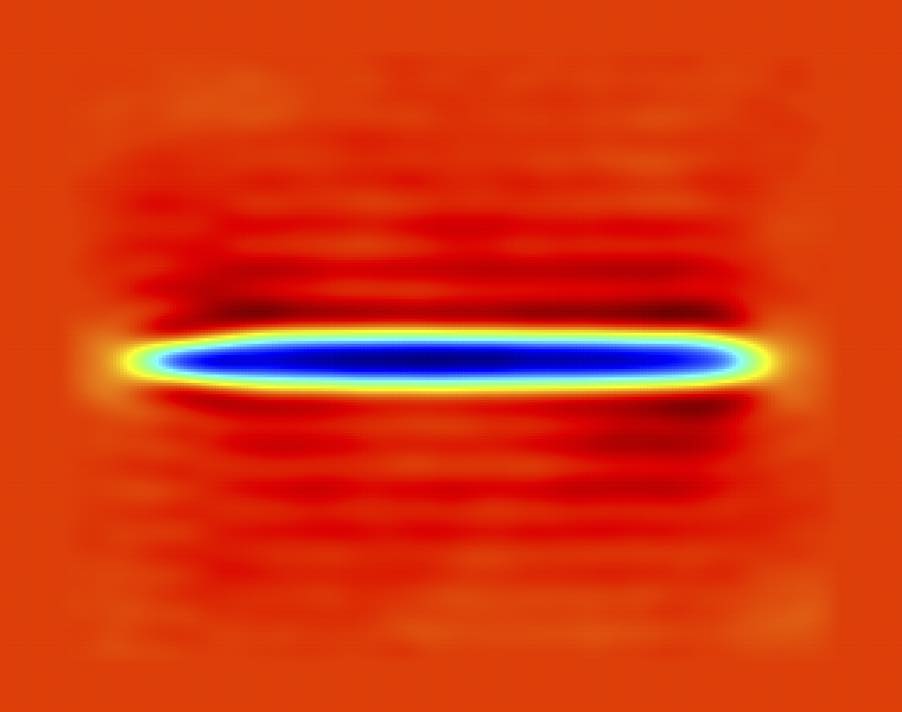}
}
\caption{Reconstruction of $q_\mathrm{sub}$ for Example \ref{example:MDP}. The {\bf MIMDP} algorithm is used to reconstruct the domain $q_\mathrm{sub}$ in the presence of an anisotropic background medium that is generated using a noise parameter $\barsigma=16$. The results are presented for the case where we have measurements of the field scattered off of $q_\mathrm{sub}$ in the presence of: (a)$N_s=10$, (b) $50$ and (c) $100$ samples of the anisotropic background domain generated using the prior \eqref{eq:tensor_sub}.}
\label{fig_pid_submarine_samples}
\end{figure}

In this example, we recover the scatterers for $\qp$ and $q_\mathrm{sub}$ using the {\bf MIMDP} algorithm. For the reconstruction of $\qp$, we consider that the background medium was generated as in Example \ref{example:SDP}, part (a) with $\barsigma=10$; meanwhile, for the reconstruction for the submarine-like object, we consider the domain to have been generated as in Example \ref{example:SDP}, part (c) with $\barsigma=16$. The {\bf SISDP} algorithm was not entirely successful at separating the data from the background medium and the scatterer at those noise levels. This time, we use the {\bf MIMDP} algorithm to recover the scatterer using data measurements generated using $N_s=10$, $50$, and $100$ samples of the background medium. The average of the reconstructions is presented in Figure \ref{fig_pid_plane_samples} and \ref{fig_pid_submarine_samples}, respectively, for $\qp$ and $q_\mathrm{sub}$. As we can see, in both examples, we were able to improve the quality of the reconstructed scatterer.

\section{Conclusion}\label{s:conclusion}
We have presented a comprehensive study of the problem of reconstructing a scatterer in the presence of a random background  medium. Six different algorithms are presented to solve this problem with different amounts of scattered field data and information about the scatterer and the background noisy medium. We consider four cases:
\begin{enumerate}[label=(\alph*)]
\item in the first case, we have data measurements of the scattered field off of $q$ in the presence of one realization of $\eta$ and $\mathbb{E}(\eta)$; 
\item in the second case, we have data measurements from several realizations of $\eta$ and $\mathbb{E}(\eta)$; 
\item in the third case, we have data measurements of the field scattered off of $q$ in the presence of one realization of $\eta$ and prior knowledge of the probability distributions of $q$ and $\eta$; and s
\item in the fourth case, we have data measurements from several realizations of $\eta$, and prior knowledge of the probability distributions of $q$ and $\eta$, and $\mathbb{E}(\eta)$. 
\end{enumerate}

The main conclusion is that, perhaps counter-intuitively, it is preferable to try to solve a harder inverse problem and invert for both the target scatterer and the random medium. Not surprisingly, the best reconstruction results are obtained when we have the highest quality of information, spectrally separated priors for both $q$ and $\eta$ and a rich dataset from multiple realizations of $\eta$. This scenario could require the need for multiple inversions although just averaging the data seems to produce good results. Finally, if $q$ and $\eta$ have similar priors it will be hard to tell them apart. In this scenario, we need data from multiple realizations of $\eta$ and multiple inversion to be able to disentangle $\eta$ and $q$.  Our simple analysis of the interplay between the priors of $\eta$ and $q$ using their spectrum could be formalized using  Kullback-Leibler divergence between the priors. 

In the future, we intend to extend the study to the case of limited aperture data and to the case when only the magnitude of the scattered field can be measured (and not its phase). We also intend to study the case when dissipation is allowed in the unknown scatterer function.

\section*{Acknowledgments}
This material is based upon work supported by AFOSR grants
FA9550-17-1-0190; and by NSF grant CCF-1337393.  Any opinions,
findings, and conclusions or recommendations expressed herein are
those of the authors and do not necessarily reflect the views of the
AFOSR, and NSF.  The authors would also like to thank Andreas Mang
and Kui Ren for several useful conversations.

\section*{Appendix A - HPS fast solver}
Consider the following forward scattering problem for $u^{\emph{scat}}$ with a source function $f(\xb)$:
\begin{equation}
\Delta u^{\emph{scat}}(\xb) + k^2 (1+q(\xb)) u^{\emph{scat}}(\xb) = f(\xb), 
\label{helm_forc_term}
\end{equation}
where $\xb \in \mathbb{R}^2$, $k$ is wavenumber, the functions $q(\xb)$ and $f(\xb)$ both have compact support in $\Omega$, and where $u^{\emph{scat}}(\xb)$ satisfies the Sommerfeld radiation condition.

We break problem \eqref{helm_forc_term} into two problems, one in the interior and other in the exterior of $\Omega$. In the interior of the domain $\Omega$, we have the problem
\begin{eqnarray}
\Delta u^{\emph{scat}}(\xb ) +k^2(1+q(\xb)) u^{\emph{scat}}(\xb ) &=&f(\xb) \quad  \mbox{in} \quad \Omega, \nonumber \\
u^{\emph{scat}}(\xb )&=&s(\xb ) \quad  \mbox{on} \quad \partial\Omega. \nonumber
\end{eqnarray}
In the exterior of $\Omega$, $u^{\emph{scat}}(\xb)$ must satisfy the constant-coefficient problem
\begin{eqnarray}
\Delta u^{\emph{scat}}(\xb ) +k^2 u^{\emph{scat}}(\xb ) &=&0 \quad\quad\quad \mbox{  in} \quad \mathbb{R}^2\setminus\Omega, \nonumber \\
u^{\emph{scat}}(\xb )&=&s(\xb ) \quad \quad \mbox{on} \quad \partial\Omega,\nonumber \\
 \frac{\partial v}{\partial r} -ikv &=& o(r^{-1/2})\quad  r=\|\xb \|\rightarrow\infty. \nonumber
\end{eqnarray}
We will assume that the interior Dirichlet problem does not have a resonance at $k$. To obtain a coupling condition for the two problems, we write $u^{\emph{scat}}(\xb) = u^{\emph{scat}}_h(\xb)+u^{\emph{scat}}_p(\xb)$, where $u^{\emph{scat}}_h(\xb)$ is the solution of the homogeneous problem and $u^{\emph{scat}}_p(x)$ is the particular solution of the problem. The particular solution can be found via partial differential discretization techniques.

It is straightforward to determine $u^{\emph{scat}}_h(\xb)$ on $\partial\Omega$ by solving the problem
\begin{equation}
\left(T^{int}-T^{ext}\right) u^{\emph{scat}}_h\vert_{\partial\Omega} = T^{ext} u^{\emph{scat}}_p - \frac{\partial u^{\emph{scat}}_p}{\partial n}\vert_{\partial\Omega},
\label{eq:fem-bem}
\end{equation}
where $T^{int}$ and $T^{ext}$ are, respectively, the interior and exterior ``Dirichlet-to-Neumann'' maps.

The construction of $T^{ext}$ has been extensively covered in the literature \cite{Colton}. Using the standard layer potentials, the scattered field $u^{\emph{scat}}(\xb)$ satisfies
\begin{equation*}
u^{\emph{scat}}(\xb)=Du^{\emph{scat}}(\xb)-S \frac{\partial u^{\emph{scat}}}{\partial n} (\xb)
\end{equation*}
for $\xb$ in the exterior of $\Omega$, where $D$ and $S$ are the double and single-layer operators, respectively. Using standard jump relations \cite{Colton}, we have
\begin{equation*}
T^{ext}=S^{-1} \left(D-\frac{I}{2} \right).
\end{equation*}

The construction of $T^{int}$ is rather complicated to fully describe and is not the objective of this article. Summarizing, the solver begins by constructing a hierarchically refined quad-tree superimposed on $\Omega$, in which, within each leaf node, a $K \times K$ tensor product Chebyshev grid is used. The Impedance-to-Impedance (ItI) operator, and an operator mapping the particular solution to the corresponding outgoing impedance data, are constructed in each leaf node. Using a bottom-up procedure, the interior ItI map for each parent node is constructed by merging its four child nodes until the root node is reached. Finally, we obtain $T^{int}$ at the root using its ItI operator.

\section*{Appendix B - Regularization parameter calculation}
To deal with the ill conditioning of the system in \eqref{eq:SISDP:system}, we need to provide a way to chose the regularization parameters $\alpha$ and $\beta$. In this appendix, we provide a heuristic to obtain $\alpha$ in the case that $\beta=0$, which is the case for Example \ref{example:SDP}, parts (a) and (b), and the reconstruction of $q_p$ in Example \ref{example:MDP}.
The same procedure is used for $\beta$. The two regularization parameters are determined independently.

To determine $\alpha$ we solve a synthetic problem and test the quality of the reconstruction for $q$. One complication is that in {\bf RLA} we solve a sequence of inverse problems with different wavenumbers $k$. Since the matrices of the system become better conditioned with increasing maximum wavenumber $k_Q$, it makes sense to look for values of the regularization parameter $\alpha$ that decrease with $k_Q$.  The scheme to find the appropriate regularization $\alpha$ for each $k$ is given Algorithm~\ref{alg:find_alpha}. 

\begin{algorithm}
\caption{Algorithm to find $\alpha$}
\label{alg:find_alpha}
\begin{algorithmic}[1]
\STATE{{\bf Input:} initial value $\alpha_0$ for $\alpha$, functions $\eta^\ast_s$ for $s=1,\ldots,10$ and data $\db_s={\bf F}(q^\ast+\eta^\ast_s)$ at $k_1<\cdots<k_Q$.}
\FOR{ $s=1,\ldots,10$}
\STATE{Use the {\bf RLA} to solve $\tilde{q}_s=\argmin_q \|\db_{\eta^\ast_s}(k)-{\bf F}_k(q+\eta^\ast_s)\|$.}
\STATE{Set $\alpha(k_1)=\alpha_0$, $m=1$, $\epsilon_0=0$ and $flag_\epsilon=true$.}
\FOR{$j=1,\ldots,Q$}
\WHILE{$m<10$ and $flag_\epsilon$}
\STATE{Solve using the Gauss-Newton the problem $\argmin_{q,\eta}\|\db_s(k_j)-{\bf F}(q,\eta)\|-\frac{\alpha_s(k_j)}{2}\log(p(\eta))$.}
\STATE{Calculate $\epsilon_1=\frac{\|q-\tilde{q}_s\|}{\|\tilde{q}_s\|}$}.
\IF{$\epsilon_j>\epsilon_{j-1}$}
\STATE{Set $flag_\epsilon=false$ and $\alpha_s(k_{j+1})=\alpha_s(k_j)$.}
\ELSE
\STATE{Set $\alpha_s(k_j)=\alpha_s(k_j)/2$.}
\ENDIF
\ENDWHILE
\ENDFOR
\STATE{Set $\alpha=\frac{1}{10}\sum_{s=1}^{10} \alpha_s$.}
\ENDFOR
\end{algorithmic}
\end{algorithm}

\bibliography{./Bibnew}

\end{document}